\documentclass{article}
\usepackage{amsmath, amssymb}

\newtheorem{thm}{Theorem}[section]
\newtheorem{cor}[thm]{Corollary}
\newtheorem{prop}[thm]{Proposition}
\newtheorem{lem}[thm]{Lemma}
\newenvironment{dfn}{\medskip\refstepcounter{thm}
\noindent{\bf Definition \thesection.\arabic{thm}\ }}{\medskip}
\newenvironment{ex}{\medskip\refstepcounter{thm}
\noindent{\bf Example \thesection.\arabic{thm}\ }}{\medskip}
\newenvironment{proof}[1][,]{\medskip\ifcat,#1
\noindent{{\it Proof}:\ }\else\noindent{\it Proof of #1.\ }\fi}
{\hfill$\square$\medskip}
\newenvironment{remark}[1][Remark]{\begin{trivlist}
\item[\hskip \labelsep {\bfseries #1}]}{\end{trivlist}}
\newenvironment{remarks}[1][Remarks]{\begin{trivlist}
\item[\hskip \labelsep {\bfseries #1}]}{\end{trivlist}}
\newenvironment{note}[1][Note]{\begin{trivlist}
\item[\hskip \labelsep {\bfseries #1}]}{\end{trivlist}}
\newenvironment{notes}[1][Notes]{\begin{trivlist}
\item[\hskip \labelsep {\bfseries #1}]}{\end{trivlist}}
\newenvironment{ack}[1][Acknowledgements]{\begin{trivlist}
\item[\hskip \labelsep {\bfseries #1}]}{\end{trivlist}}

\def\eq#1{{\rm(\ref{#1})}}
\def\N{{\mathbb N}}
\def\Z{{\mathbb Z}}
\DeclareMathOperator\vol{vol}
\def\R{{\mathbb R}}
\DeclareMathOperator\supp{supp}
\def\d{{\rm d}}
\DeclareMathOperator\loc{loc}
\def\w{\wedge}
\DeclareMathOperator\GG{G}
\DeclareMathOperator\GL{GL}
\def\C{{\mathbb C}}
\def\Re{\mathop{\rm Re}\nolimits}
\def\Im{\mathop{\rm Im}\nolimits}
\def\P{{\mathbb P}}
\DeclareMathOperator\SU{SU}
\def\Nhat{\hat{N}}
\DeclareMathOperator\Ker{Ker}
\DeclareMathOperator\Image{Image}
\DeclareMathOperator\ind{ind}
\DeclareMathOperator\SO{SO}
\def\O{{\mathbb O}}
\def\eps{\varepsilon}

\def\bff{{\bf f}}
\DeclareMathOperator\U{U}
\DeclareMathOperator\Sp{Sp}
\def\bfx{{\bf x}}
\def\bfo{{\bf 1}}
\def\bfi{{\bf i}}
\def\bfj{{\bf j}}
\def\bfk{{\bf k}}
\def\H{{\mathbb H}}
\DeclareMathOperator\DD{D}
\DeclareMathOperator\SSS{S}
\DeclareMathOperator\AAA{A}
\def\bfu{{\bf u}}

\begin{document}

\title{Stability of Coassociative Conical Singularities}
\author{\textsc{Jason D.~Lotay\footnote{Address for correspondence: Department of Mathematics, University College London, Gower Street, London WC1E 6BT, England, U.K. Email: j.lotay@ucl.ac.uk.}}\\ Imperial College London}
\date{}

\maketitle

\begin{abstract}
\noindent We study the stability of coassociative 4-folds with conical singularities under perturbations of the ambient $\GG_2$ structure by 
defining an integer invariant of a coassociative cone which we call the stability index.  The stability index of a coassociative 
cone is determined by the spectrum of the curl operator acting on its link.  We explicitly calculate the stability index for cones 
on group orbits.  We also describe the stability index for cones fibered by 2-planes over algebraic curves using the degree and
 genus of the curve and the spectrum of the Laplacian on the link.  Finally we apply our results to construct the first 
 known examples of coassociative 4-folds with conical singularities in compact manifolds with $\GG_2$ holonomy.
\end{abstract}

\section{Introduction}

Coassociative 4-folds are \emph{calibrated}, hence \emph{minimal}, submanifolds of 7-manifolds with $\GG_2$ structures,   
 first defined in \cite{HarLaw}.
  Of particular interest are coassociative 4-folds in manifolds with $\GG_2$ holonomy, of which few examples are known. 
Coassociative submanifolds with conical singularities have been previously studied by the author in \cite{Lotaycs} and 
\cite{Lotaydesing}, building upon the work on special Lagrangian submanifolds with conical singularities by Joyce in 
\cite{JoyceSLsing5}-\cite{JoyceSLsing4}.  

We continue to generalize the work of Joyce to the coassociative setting by defining
the notion of \emph{stability index} for coassociative cones.  The stability index is a non-negative integer invariant for 
a coassociative cone, the vanishing of which guarantees that coassociative 4-folds with a singularity modelled on that cone 
 will have a smooth moduli space of deformations:  in particular, they are \emph{stable} under small perturbations  
of the $\GG_2$ structure on the 7-manifold.  

We calculate the stability index for certain types of 
coassociative cones and 
use our results to construct the first examples of coassociative 4-folds with conical singularities in compact manifolds 
with $\GG_2$ holonomy.  This is an essential step in the proposed construction given 
in \cite{Kovalevfib} 
 of coassociative fibrations of compact $\GG_2$ manifolds.

\subsection{Motivation}

Even though there has been a wealth of research devoted to calibrated submanifolds with conical singularities 
in manifolds with special holonomy, particularly in \cite{JoyceSLsing1}-\cite{JoyceSLsing5} and \cite{Lotaycs}-\cite{Lotaydesing},
 there were no known examples of such submanifolds.  Given this well-developed theory, 
we were motivated to construct coassociative examples.  As far as the author is aware, it is unknown whether there are special Lagrangian $m$-folds with conical
 singularities in Calabi--Yau $m$-folds for $m\geq 3$.

The stability index of a special Lagrangian cone is defined in \cite{JoyceSLsing1} in terms of the spectrum of the Laplacian on 
the link, and is calculated for cones over flat tori originally given in \cite[$\S$III]{HarLaw}.  
The stability index for special Lagrangian cones over 
certain homogeneous spaces is calculated in \cite{Ohnita} and Haskins \cite{Haskins} shows that the only stable 
special Lagrangian $T^2$-cone is the cone over the flat torus.   

In contrast, the stability index for a coassociative cone is determined 
by the spectrum of the \emph{curl operator} acting on 1-forms on the 3-dimensional link of the cone.  
Though still a natural geometric object, 
there is relatively little material on the spectrum of the curl operator in the literature, and it is certainly 
less straightforward to analyse than the spectrum of the Laplacian which has received so much attention.   
We are thus required to undertake fundamental elementary calculations to describe the stability index even in simple cases. 

 Eigenforms for the curl operator naturally define contact structures and 
are dual to certain Beltrami fields, which are important in hydrodynamics.  
Beltrami fields are also equivalent to Reeb vector fields by the work in \cite{Etnyre}, so are of particular 
interest in contact geometry.  We hope therefore that our spectral calculations may be of wider benefit.

Motivated by the SYZ conjecture (see \cite{SYZ}), one would hope to construct coassociative fibrations of compact manifolds with 
$\GG_2$ holonomy.  An elementary argument shows that some of the fibres must necessarily be singular, so it is natural to assume that 
they have the simplest type of singularity, namely conical singularities.  As already mentioned, the stability results in this paper
 are essential for the proposed construction of a coassociative $K3$ fibration in \cite{Kovalevfib}.  Moreover, our results  
will almost certainly be useful for any other construction of a coassociative fibration with conically singular fibres.

\subsection{Summary}

We begin in $\S$\ref{coasssec} with the basic definitions we need, including concepts from calibrated geometry and 
Geometric Measure Theory.  We also discuss the relationship between complex and coassociative geometry.

In $\S$\ref{cssec} we mainly review the material in \cite{Lotaycs}, defining coassociative 4-folds with conical singularities
and describing their deformation theory.  Our definition of coassociative conical singularities (Definition \ref{csdfn}) on the 
face of it looks rather strong.  However, we prove the following regularity result, generalizing results in \cite{Simon} 
and \cite{JoyceSLsing1}, which shows that the definition is more applicable than it appears.

\begin{thm}\label{integralthm} If a coassociative integral current has a multiplicity one, Jacobi integrable tangent cone with
 isolated singularity at an interior point $p$, then it has a conical singularity at $p$. 
\end{thm}

\begin{remarks} The link of a coassociative cone is a \emph{Lagrangian} (or totally real) submanifold of the nearly K\"ahler 6-sphere $\mathcal{S}^6$.
A coassociative cone is \emph{Jacobi integrable} if every infinitesimal variation of its link 
as a Lagrangian in $\mathcal{S}^6$ is integrable.
 \end{remarks}

In $\S$\ref{stabindsec} we begin by briefly reviewing the known examples of coassociative cones. Then, by refining some of the results 
in \cite{Lotaycs}, we define the notion of $\mathcal{C}$-stability index for 
a coassociative cone $C$ in a family $\mathcal{C}$.  We simply call the $\mathcal{C}$-stability index of $C$ the 
stability index when $\mathcal{C}$ is the family of cones generated by $\GG_2$ transformations and translations of $C$.

Sections \ref{homconessec} and \ref{algcurvessec} are devoted to calculating the stability index for certain coassociative cones.  
In $\S$\ref{homconessec} we analyse the curl operator acting on Berger 3-spheres and their quotients by finite groups, and thus 
determine the stability index for all homogeneous coassociative cones.  We deduce the following.

\begin{thm}
The only stable homogeneous coassociative cones are coassociative 4-planes and the $\Sp(1)$-invariant coassociative cone 
given in Example \ref{u2ex}.
\end{thm}
\begin{remark}
The cone in Example \ref{u2ex} was originally constructed in \cite[$\S$7]{LO} and identified as coassociative in \cite[$\S$IV]{HarLaw}.
\end{remark}

In $\S$\ref{algcurvessec} we study coassociative cones which 
are fibered by 2-planes over algebraic curves: either \emph{holomorphic curves} in $\C\P^2$, where the corresponding 
coassociative cone is complex, or \emph{null-torsion pseudoholomorphic curves} in $\mathcal{S}^6$.  In both cases we express 
the stability index in terms of the spectrum of the Laplacian acting on functions on the link and algebro-geometric data from the 
curve.  These results allow us to examine the behaviour of the stability index under deformations of these types of coassociative
 cones.

Finally, in $\S$\ref{exsec}, after discussing the construction in \cite{Kovalevsums} of compact manifolds with $\GG_2$ holonomy, 
we apply the results of $\S$\ref{homconessec} and $\S$\ref{algcurvessec} to prove the following.

\begin{thm}\label{mainthm}
Given a pair of maximal deformation families of Fano 3-folds, one can construct a one-parameter family of compact manifolds 
with $\GG_2$ holonomy, which 
contain coassociative $K3$ surfaces with conical singularities.
\end{thm}

\begin{remark}
The key ingredients in the proof are Theorem \ref{integralthm}, the $\mathcal{C}$-stability index for a homogeneous complex cone $C$
in a deformation family $\mathcal{C}$ and the invariance of the $\mathcal{C}$-stability index under deformations of $C$ 
in $\mathcal{C}$.
\end{remark}

Given the proof of Theorem \ref{mainthm}, it is clear that this result will naturally extend to give coassociative 4-folds with conical singularities in the compact holonomy $\GG_2$ manifolds constructed by 
Corti, Haskins, Nordstr\"om and Pacini \cite{CHNP}, where one replaces Fano 3-folds by a more general class of complex 3-folds.

\begin{notes}\begin{itemize}\item[]
\item[(a)] Manifolds are taken to be nonsingular and submanifolds to
be embedded, for convenience, unless stated otherwise.
\item[(b)] We use the convention that $\N=\{0,1,2,\ldots\}$ and $\Z^+=\N\setminus\{0\}$.
\end{itemize}\end{notes}

\section{Coassociative 4-folds}\label{coasssec}

In this section we cover all of the basic definitions and theory we need.

\subsection{Calibrated geometry and Geometric Measure Theory}

We will need some general theory from \emph{calibrated geometry} later, so we start with the following definition.

\begin{dfn}\label{caldfn}
Let $(M,g)$ be a Riemannian manifold.  An \emph{oriented tangent $m$-plane} $W$ on $M$ is an oriented $m$-dimensional vector subspace
$W$ of $T_xM$, for some $x$ in $M$. Given an oriented tangent $m$-plane $W$ on $M$, $g|_W$ is a Euclidean metric on $W$ and hence,
using $g|_W$ and the orientation on $W$, there is a natural volume
form, $\vol_W$, which is an $m$-form on $W$.

A closed $m$-form $\phi$ on $M$ is a \emph{calibration} on $M$ if $\phi|_W\leq\vol_W$ for all oriented tangent $m$-planes $W$ on $M$,
where $\phi|_W = \kappa\cdot\vol_W$ for some $\kappa\in\R$, so $\phi|_W\leq\vol_W$ if $\kappa\leq 1$.
An oriented $m$-dimensional submanifold $S$ of $M$ is a \emph{calibrated submanifold} or \emph{$\phi$-submanifold} 
if $\phi|_{S} =\vol_{S}$.
\end{dfn}

We shall need some ideas from Geometric Measure Theory.  A good introduction to this theory, which we only need in a 
superficial way, is given in \cite{Morgan}.

\begin{dfn}\label{GMTdfn}  Let $(M,g)$ be a complete Riemannian manifold and let $\mathcal{H}^m$ be $m$-dimensional Hausdorff measure 
on $M$.  By \cite[Proposition 3.11]{Morgan}, one can define an $\mathcal{H}^m$-measurable
subset $S$ of $M$ to be an \emph{$m$-dimensional rectifiable set} if $S$ has finite $\mathcal{H}^m$-measure and 
$\mathcal{H}^m$-almost all of $S$ is covered by the disjoint union of a countable number of compact 
$C^1$-submanifolds 
of $M$.  Let $\vol(S)$ be the $\mathcal{H}^m$-measure of $S$.  
 By \cite[Proposition 3.12]{Morgan}, an $m$-dimensional rectifiable set $S$ in $M$ has a well-defined $m$-dimensional tangent plane  
$\mathcal{H}^m$-almost everywhere, and so is orientable almost everywhere.  If $S$ is an oriented rectifiable set, let 
$\mathbf{s}(x)$ be the unit $m$-vector to $S$ at $x$ given by the choice of orientation, when this is well-defined.
 
Let $\mathcal{D}^m(M)$ be the space of smooth compactly supported $m$-forms on $M$.  An \emph{$m$-dimensional current} on $M$ is an
 element of the dual space $\mathcal{C}^m(M)=\mathcal{D}^m(M)^*$ and we define the support of $T\in\mathcal{C}^m(M)$,
 $\supp T$, to be the smallest closed set in $M$ such that, for any $\xi\in\mathcal{D}^m(M)$, $\supp\xi\cap\supp T=\emptyset$ implies that 
 $T(\xi)=0$.  We define the \emph{boundary} $\partial T\in\mathcal{C}^{m-1}(M)$ of $T\in\mathcal{C}^m(M)$ via the formula $\partial
 T(\xi)=T(\d\xi)$ for $\xi\in\mathcal{D}^{m-1}(M)$.  We also define the \emph{interior} $T^{\circ}$ of $T$ to be the set 
$\supp T\setminus\supp \partial T$.

Given an $m$-dimensional oriented rectifiable set $S$ and a function $\nu:S\rightarrow\mathbb{Z}^+$ such that
 $\int_S\nu(x)\,\d\mathcal{H}^m<\infty$, we define an associated element $T_S$ of $\mathcal{C}^m(M)$ via
$$T_S(\xi)=\int_S \mathbf{s}(x)\cdot\xi(x)\, \nu(x)\,\d\mathcal{H}^m.$$
If $\supp T_S$ is compact, we say that $T_S$ (or simply $S$, since $T_S$ is defined by $S$) is an $m$-dimensional 
\emph{rectifiable current} and we denote the set 
of $m$-dimensional rectifiable currents by $\mathcal{R}^m(M)$.  We also let 
$\mathcal{I}^m(M)=\{S\in\mathcal{R}^m(M)\,:\,\partial S\in\mathcal{R}^{m-1}(M)\}$ be the set of 
$m$-dimensional \emph{integral currents}.

Finally, we define the set of $m$-dimensional \emph{locally integral currents} by $$\mathcal{I}^m_{\loc}(M)=\{T\in\mathcal{C}^m(M)\,:\,\forall x\in M\, \exists S\in\mathcal{I}^m(M)\,\text{with}\,x\notin
\supp(T-S)\}$$
 and similarly define $\mathcal{R}^m_{\loc}(M)$.
\end{dfn}

\noindent The idea behind a rectifiable current is to generalize the notion of a compact $C^1$-submanifold with boundary to include
 multiplicities (given by the function $\nu$) and to allow for singular behaviour.  Since the boundary of a rectifiable current 
may be very badly behaved, we often require that the rectifiable current be integral.  As we shall need to deal with planes and cones 
which are definitely not compact, we expand our notation to include the ``local'' versions of the integral and 
rectifiable currents.

The majority of this paper will be dedicated to the study of cones, so we make some formal definitions for convenience.

\begin{dfn}\label{conedfn}  Recall the notation of Definition \ref{GMTdfn}, let $W$ be a normed vector space 
and let $\mathcal{S}(W)$ be the unit sphere in $W$ with respect to the norm.  
An element $C\in\mathcal{R}^m_{\loc}(W)$ is a \emph{cone} 
 in $W$ if $tC=C$ for all $t>0$, and we call  $C\cap \mathcal{S}(W)$ the \emph{link} of $C$.  
\end{dfn}

We formally define convergence in the space of currents as follows.

\begin{dfn}\label{currentcvgedfn}  Recall the notation of Definition \ref{GMTdfn}.  
We say that a sequence $(S_j)$ in $\mathcal{R}^m_{\loc}(M)$ converges to $S\in\mathcal{R}^m_{\loc}(M)$ if
 $S_j\rightarrow S$ in the weak topology in $\mathcal{C}^m(M)$; that is,
$$S_j\rightarrow S\quad\text{if and only if}\quad \int_{S_j}\xi\rightarrow\int_S\xi\quad\text{as $j\rightarrow\infty$}$$
for all $\xi\in\mathcal{D}^m(M)$, where integration is carried out with respect to $\mathcal{H}^m$ and includes multiplicities.
\end{dfn}

In the seminal work on calibrated geometry \cite{HarLaw}, the relationship between calibrated geometry and Geometric Measure 
Theory is discussed at length.  We note some of the observations originally presented there.

\begin{dfn}\label{calGMTdfn}  Let $(M,g)$ be a complete Riemannian manifold, let $\phi$ be an $m$-form which is a calibration on $M$ 
and recall the notation of Definition \ref{GMTdfn}.  
By the work in \cite{HarLaw}, we can define $S\in\mathcal{I}^m(M)$ to be an \emph{integral $\phi$-current} if $S$ is
 calibrated with respect to $\phi$; that is, $\int_S\phi=\vol(S)$.  
  Then integral $\phi$-currents are volume-minimizing in their homology class.  We can also define a \emph{locally integral 
  $\phi$-current} in $M$ in the obvious manner.
\end{dfn}

One of the key ideas in Geometric Measure Theory is the concept of a \emph{tangent cone}, which we now define.  

\begin{dfn}\label{tgtconedfn}  Let $(M,g)$ be a complete Riemannian $n$-manifold and recall the notation of Definition \ref{GMTdfn}.  
Let $S\in\mathcal{R}^m_{\loc}(M)$ and let $x\in S^\circ$.  Choose a diffeomorphism $\upsilon:V\rightarrow
B$, where $B$ is an open neighbourhood of the origin in $\R^n$ and $V$ is an open neighbourhood of $x$ in $M$.  Let $U=V\cap S$ and let $\Upsilon=\d\upsilon|_x$, which 
is an isomorphism between $T_xM$ and $\R^n$.    
A \emph{tangent cone} for $S$ at $x$ is a cone $C$ in $T_xM$ such that there exists a strictly decreasing positive sequence
 $(r_j)$, converging to zero as $j\rightarrow\infty$, such  that 
$$r_j^{-1}\upsilon(U\setminus\{x\})\rightarrow \Upsilon(C)\qquad\text{as $j\rightarrow\infty$}$$ in the sense 
of Definition \ref{currentcvgedfn}. 
\end{dfn}

\noindent In Geometric Measure Theory there are two notions of `tangent cone': one is a set and the other is a current.  We have 
defined a current in Definition \ref{tgtconedfn} and so it is, strictly speaking, an \emph{oriented tangent cone} in the sense of 
Geometric Measure Theory.  Oriented tangent cones can be defined for more general currents than locally rectifiable ones, but then 
one only requires a weaker form of convergence in the definition. 

We conclude with an important result that follows from \cite[Theorem 4.4.4]{JoyceRiem}. 

\begin{prop}\label{caltgtconeprop}
 Let $\phi$ be a calibration on a complete Riemannian manifold $M$ and let $S$ be an integral $\phi$-current in $M$.   
There exists a tangent cone to $S$ at each $x\in S^\circ$, and it is a locally integral $\phi|_x$-current in $T_xM$.
\end{prop}

\subsection{Calibrated geometry in $\R^7$}

We define coassociative 4-folds in $\R^7$ by introducing a distinguished 3-form.

\begin{dfn}\label{phisdfn} Let $(x_1,\ldots,x_7)$ be coordinates on $\R^7$ and write
$\d{\bf x}_{ij\ldots k}$ for the form $\d x_i\w \d x_j\w\ldots\w \d x_k$.
Define a 3-form $\varphi_0$ on $\R^7$ by:
\begin{equation}\label{phieq}
\varphi_0 = \d{\bf x}_{123}+\d{\bf x}_{145}+\d{\bf x}_{167}+\d{\bf
x}_{246}- \d{\bf x}_{257}-\d{\bf x}_{347}-\d{\bf x}_{356}.
\end{equation}
The Hodge dual of $\varphi_0$ is a 4-form given by:
\begin{equation}\label{starphieq}
\ast\varphi_0 = \d{\bf x}_{4567}+\d{\bf x}_{2367}+\d{\bf
x}_{2345}+\d{\bf x}_{1357}-\d{\bf x}_{1346}-\d{\bf x}_{1256}-\d{\bf
x}_{1247}.
\end{equation}
The forms $\varphi_0$ and $*\varphi_0$ are calibrations by \cite[Theorems IV.1.4 \& IV.1.16]{HarLaw}.  Submanifolds calibrated with 
respect to $\varphi_0$ and $*\varphi_0$ are called \emph{associative 3-folds} and \emph{coassociative 4-folds} 
respectively. We can also characterize the coassociative 4-folds as the oriented 4-dimensional submanifolds $N$ in
$\R^7$ satisfying $\varphi_0|_N\equiv0$, oriented such that $*\varphi_0|_N>0$, by \cite[Proposition IV.4.5 \&
Theorem IV.4.6]{HarLaw}.
\end{dfn}

\begin{remark}
 The form $\varphi_0$ is sometimes called the `$\GG_2$ 3-form' because the exceptional Lie group $\GG_2$ is the stabilizer of 
$\varphi_0$ in $\GL(7,\R)$.  
\end{remark}

A straightforward calculation yields the following lemma.

\begin{lem}\label{phisplitlem}
Identify $\R^7$ with $\R\oplus\C^3$ so that $x_1$ is the coordinate on $\R$ and $z_1=x_2+ix_3$, 
$z_2=x_4+ix_5$ and $z_3=x_6+ix_7$ are coordinates on $\C^3$.  If $\omega_0$ and $\Omega_0$ are the 
standard K\"ahler and holomorphic forms on $\C^3$, then:
\begin{align}
\varphi_0&=\d x_1\w\omega_0+\Re\Omega_0;\label{phispliteq}\\
*\varphi_0&=\frac{1}{2}\,\omega_0\w\omega_0-\d x_1\w\Im\Omega_0\label{starphispliteq},
\end{align}
where $\varphi_0$ and $*\varphi_0$ are given in \eq{phieq}-\eq{starphieq}.
\end{lem}

Since $\Re\Omega_0$ and $\Im\Omega_0$ are both calibrations on $\C^3$ we have the following definition, again due to 
Harvey and Lawson \cite{HarLaw}.

\begin{dfn}\label{SLdfn}
Let $(z_1,\ldots,z_m)$ be coordinates on $\C^m$ and let $\omega_0$ and $\Omega_0$ be the K\"ahler and holomorphic volume forms
 on $\C^m$.  Then $\cos\theta\Re\Omega_0+\sin\theta\Im\Omega_0$ is a calibration on $\C^m$ for all real constants $\theta$, and 
 its corresponding calibrated submanifolds are real $m$-dimensional submanifolds of $\C^m$ called \emph{special Lagrangian $m$-folds} 
 (with $\emph{phase}$ $e^{i\theta}$).   Moreover, special Lagrangian $m$-folds 
 with phase $e^{i\theta}$ are the oriented real $m$-dimensional submanifolds $L$ of $\C^m$ such that $\omega_0|_L\equiv 0$ and
 $(\sin\theta\Re\Omega_0-\cos\theta\Im\Omega_0)|_L\equiv 0$, up to a choice of orientation.
\end{dfn}

Examination of \eq{phispliteq}-\eq{starphispliteq} immediately yields the following elementary result.

\begin{cor}\label{coasscor}
In the notation of Lemma \ref{phisplitlem}, $\R\times L\subseteq\R\oplus\C^3$ and $N\subseteq\C^3$ are 
coassociative in $\R^7$ if and only if $L$ is a special Lagrangian 3-fold with phase $-i$ and 
$N$ is a complex surface in $\C^3$ respectively.
\end{cor}

Since we will be concerned with coassociative cones, we make the following convenient definition.

\begin{dfn}\label{Lagdef}
The 6-sphere $\mathcal{S}^6$ inherits a \emph{nearly K\"ahler} structure from the standard $\GG_2$ structure on $\R^7$.  In 
particular, if $r$ is the radial coordinate and $\mathbf{e}_r$ is the radial vector field on $\R^7$,
 $\omega=(\mathbf{e}_r\cdot\varphi_0)|_{r=1}$ is a non-degenerate 2-form on $\mathcal{S}^6$ which is \emph{not} closed.   Using this 
 2-form and the round metric $g$ on $\mathcal{S}^6$ we can define an almost complex structure $J$ by $\omega(x,y)=g(Jx,y)$ 
 for tangent vectors $x,y$.  The almost complex structure $J$ is \emph{not} integrable.
 
An oriented 3-dimensional submanifold $L\subseteq\mathcal{S}^6$ is the link of a coassociative cone 
if and only if $\omega|_{L}\equiv 0$.  Thus, 
we say that the link of a coassociative cone in $\R^7$ is a \emph{Lagrangian} submanifold of $\mathcal{S}^6$.

An oriented surface $\Sigma\subseteq\mathcal{S}^6$ is a \emph{pseudoholomorphic curve} if and only if 
  $\omega|_{\Sigma}=\vol_\Sigma$ or, equivalently, if $J(T_\sigma\Sigma)=T_\sigma\Sigma$ for all $\sigma\in\Sigma$.  
Note that $\Sigma$ is the link of an associative cone if and only if 
 $\Sigma$ is a pseudoholomorphic curve.
\end{dfn}

\noindent By Corollary \ref{coasscor}, any complex 2-dimensional cone in $\C^3$ is coassociative in $\R\oplus\C^3\cong\R^7$.  Thus, 
the Hopf lift of any holomorphic curve in $\C\P^2$ to a totally geodesic $\mathcal{S}^5$ in $\mathcal{S}^6$ is Lagrangian.
 Since special Lagrangian 3-folds in $\C^3$ are associative in $\R\oplus\C^3\cong\R^7$, minimal Legendrian surfaces in 
 $\mathcal{S}^5$, which are the links of special Lagrangian cones in $\C^3$, give examples of pseudoholomorphic curves in
 $\mathcal{S}^6$.

\subsection{{\boldmath $\GG_2$} structures}

So that we may define coassociative submanifolds of more general 7-manifolds, we make the following definition.

\begin{dfn}\label{pos3formdfn}
Let $M$ be an oriented 7-manifold and recall the 3-form $\varphi_0$ on $\R^7$ given in \eq{phieq}.  For each $x\in M$ there exists
an orientation preserving isomorphism $\iota_x:T_xM\rightarrow\R^7$.
Since $\dim\GG_2=14$, $\dim\GL_+(T_xM)=49$ and
$\dim\Lambda^3T^*_xM=35$, the $\GL_+(T_xM)$ orbit of
$\iota_x^*(\varphi_0)$ in $\Lambda^3T^*_xM$, denoted
$\Lambda^3_+T^*_xM$, is open.  A 3-form $\varphi$ on $M$ is
\emph{positive} 
if
$\varphi|_{x}\in\Lambda^3_+T^*_xM$ for all $x\in M$.  Denote the
bundle of positive 3-forms by $\Lambda^3_+T^*M$.  
\end{dfn}

\noindent A positive 3-form is identified with the
$\GG_2$ 3-form $\varphi_0$ on $\R^7$ at each point in $M$.  Therefore, to each
positive 3-form $\varphi$ we can uniquely associate a 4-form
$*_{\varphi}\varphi$ and a metric $g_\varphi$ on $M$ such that the triple
$(\varphi,*_{\varphi}\varphi,g_{\varphi})$ corresponds to $(\varphi_0,*\varphi_0,g_0)$
at each point.  Notice that since the metric $g_\varphi$ depends on $\varphi$, the Hodge star $*_\varphi$ depends on $\varphi$ also.  This leads us to our next definition.

\begin{dfn}\label{G2structdfn}
Let $M$ be an oriented 7-manifold and let $\varphi\in C^{\infty}(\Lambda^3_+T^*M)$.  If $g_{\varphi}$ is the metric associated with $\varphi$, 
we call $(\varphi,g_{\varphi})$ a $\GG_2$ \emph{structure} on $M$. If $\varphi$ is
closed (or  $*_{\varphi}\varphi$ is closed) then $(\varphi,g_{\varphi})$ is a \emph{closed} (or
\emph{coclosed}) $\GG_2$ structure.  A closed and coclosed $\GG_2$
structure is called \emph{torsion-free}.
\end{dfn}

\noindent Our choice of notation here agrees with \cite{BryantG2struct}.

\begin{remark}
By \cite[Lemma 11.5]{Salamon}, $(\varphi,g_{\varphi})$ is a torsion-free $\GG_2$ structure on $M$ if and only if 
the holonomy of $g_{\varphi}$ is contained in $\GG_2$.
\end{remark}

\begin{dfn}\label{G2mflddfn} Let $M$ be an oriented
7-manifold with a $\GG_2$ structure $(\varphi,g_{\varphi})$, denoted
$(M,\varphi,g_{\varphi})$. If $(\varphi,g_{\varphi})$ is closed, we say that 
$(M,\varphi,g_{\varphi})$ is an \emph{almost
$\GG_2$ manifold}. If $(\varphi,g_{\varphi})$ is torsion-free, we call $(M,\varphi,g_{\varphi})$ a 
$\GG_2$ \emph{manifold}.
\end{dfn}

\begin{note}
By \cite[Proposition 11.2.1]{JoyceRiem}, the metric $g_{\varphi}$ 
on a compact $\GG_2$ manifold $M$ 
has $\GG_2$ holonomy
if and only if the fundamental group $\pi_1(M)$ 
is finite.
\end{note}

We are now able to complete our definitions regarding coassociative
4-folds.

\begin{dfn}\label{coassdfn}
A 4-dimensional submanifold $N$ of\/ $(M,\varphi,g_{\varphi})$ is
coassociative if and only if\/ $\varphi|_N\equiv 0$ and
$*_\varphi\varphi|_N>0$.
\end{dfn}

\begin{note} Though we may define coassociative 4-folds with respect to any $\GG_2$ structure, for 
deformation theory and related results to hold we need it to be closed. 
Therefore, we shall work with almost $\GG_2$ manifolds for greatest useful generality.
\end{note}

The next result, \cite[cf.~Proposition 4.2]{McLean},
is invaluable in describing the deformation theory of coassociative
4-folds.

\begin{prop}\label{jmathprop}
Let $N$ be a coassociative 4-fold in an almost $\GG_2$ manifold $(M,\varphi,g_{\varphi})$. There is an
isometric isomorphism between the normal bundle $\nu(N)$ of $N$ in $M$ and
$(\Lambda^2_+)_{g_{\varphi}|_N}T^*N$ given by $\jmath_N:v\mapsto(v\lrcorner\varphi)|_{TN}$.  Thus, infinitesimal
coassociative deformations of $N$ are governed by closed self-dual 2-forms on $N$.
\end{prop}

\begin{remarks}
From Proposition \ref{jmathprop} and some further analysis, one may deduce as in \cite[Theorem 12.3.4]{JoyceRiem}, by following 
\cite[Theorem 4.5]{McLean}, that 
the moduli space of deformations of a compact coassociative 4-fold $N$ in an almost $\GG_2$ manifold 
is a  manifold of dimension $b^2_+(N)$.  The author \cite{Lotaycs} adapted this deformation theory result to the situation 
where $N$ has \emph{conical singularities}, which will be invaluable for the study in this article.  
The relevant material in \cite{Lotaycs} will be the focus of $\S$\ref{cssec}.
\end{remarks}

We shall also be briefly concerned with $\SU(3)$ structures, so they form the subject of the next definition.

\begin{dfn}\label{CYdfn} On $\C^3$, let $g_0$, $\omega_0$ and $\Omega_0$ denote the standard Euclidean metric, K\"ahler form 
and holomorphic volume form respectively.  Let $(Y,J,g,\omega)$ be an almost Hermitian 6-manifold; that is, $g$ is a 
Riemannian metric on the almost complex 6-manifold $Y$, $J$ is an almost complex structure preserved by $g$ and $\omega$ is 
the associated (non-degenerate) $(1,1)$-form on $M$.  An \emph{$\SU(3)$ structure} on $Y$ is a choice of nowhere 
vanishing $(3,0)$-form $\Omega$ on $Y$ such that, for 
all $y\in Y$, there exists an orientation preserving isomorphism $\iota_y:T_yY\rightarrow\C^3$ satisfying $\iota_y^*(g_0)=g|_y$, 
$\iota_y^*(\omega_0)=\omega|_y$ and $\iota_y^*(\Omega_0)=\Omega|_y$.  

If $(Y,J,g,\omega,\Omega)$ is an almost Hermitian 6-manifold endowed with an $\SU(3)$ structure, the 
 product 7-manifold $M=\R\times Y$ (or $\mathcal{S}^1\times Y$) has a `product' $\GG_2$ structure given by 
$\varphi=\d x\w\omega+\Re\Omega$ and $g_{\varphi}=\d x^2+g$, where $x$ is the coordinate on $\R$ or $\mathcal{S}^1$, 
by \cite[Proposition 11.1.9]{JoyceRiem}.  
Moreover, $*_\varphi\varphi=\frac{1}{2}\,\omega\w\omega-\d x\w\Im\Omega$.
\end{dfn}

\begin{notes}
If $(Y,g,J,\omega)$ is a compact K\"ahler 3-manifold and $\Omega$ is a nowhere vanishing \emph{holomorphic} $(3,0)$-form on 
$Y$, then $(Y,g,J,\omega,\Omega)$ is called an \emph{almost Calabi--Yau 3-fold}.  If, in addition, 
$\omega^3=\frac{3i}{4}\Omega\w\bar{\Omega}$ then the $\SU(3)$ structure is \emph{torsion-free} and $(Y,J,g,\omega,\Omega)$ is called a \emph{Calabi--Yau 3-fold}.  These are the 
natural `$\SU(3)$ analogues' of the manifolds defined in Definition \ref{G2mflddfn}.  In particular, the Calabi--Yau condition is 
equivalent to saying that the compact K\"ahler manifold has metric $g$ with holonomy contained in $\SU(3)$.
\end{notes}

\section{Conical singularities}\label{cssec}

In this section we review some of the theory of conical singularities of coassociative 4-folds as described in \cite{Lotaycs}.   
We also prove an important new result which shows that singular coassociative integral currents with particularly 
``nice'' tangent cones have conical singularities. 

\subsection{Coassociative 4-folds with conical singularities}

We first define a preferred choice of local coordinates on an almost $\GG_2$ 
manifold near a finite set of points, which is an analogue of one given for
almost Calabi--Yau manifolds in \cite[Definition 3.6]{JoyceSLsing1}.  We let $B(0;\delta)\subseteq\R^7$ denote the open ball 
about $0$ with radius $\delta>0$.

\begin{dfn}\label{G2coordsdfn}
Let $(M,\varphi,g_{\varphi})$ be an almost $\GG_2$ manifold and let
$z_1,\ldots,z_s\in M$ be distinct points. There exist a constant $\epsilon_M\in(0,1)$,
an open set $V_i\ni z_i$ in $M$ with $V_i\cap V_j=\emptyset$ for
$j\neq i$ and a diffeomorphism
$\chi_i:B(0;\epsilon_M)\rightarrow V_i$ with $\chi_i(0)=z_i$,
for $i=1,\ldots,s$, such that $\zeta_i=\d\chi_i|_0:\R^7\rightarrow
T_{z_i}M$ is an isomorphism identifying the standard $\GG_2$
structure $(\varphi_0,g_0)$ on $\R^7$ with the pair
$(\varphi|_{z_i},g_{\varphi}|_{z_i})$. We call the set
$\{\chi_i:B(0;\epsilon_M)\rightarrow V_i:i=1,\ldots,s\}$ a
\emph{$\GG_2$ coordinate system near
${z_1,\ldots,z_s}$}.

We say that two $\GG_2$ coordinate systems near
$z_1,\ldots,z_s$, with maps $\chi_i$ and $\tilde{\chi}_i$ for
$i=1,\ldots,s$ respectively, are \emph{equivalent} if
$\d\tilde{\chi}_i|_0=\d\chi_i|_0=\zeta_i$ for all $i$.
\end{dfn}

\begin{dfn}\label{csdfn}
Let $(M,\varphi,g_{\varphi})$ be an almost $\GG_2$ manifold, let $N\subseteq M$
be compact and connected and let $z_1,\ldots,z_s\in N$ be distinct. 
Let $\{\chi_i:B(0;\epsilon_M)\rightarrow
V_i:i=1,\ldots,s\}$ be a $\GG_2$ coordinate system near
$z_1,\ldots,z_s$, as in Definition \ref{G2coordsdfn}.
We say that
$N$ is a 4-fold in $M$ with \emph{conical singularities at
$z_1,\ldots,z_s$ with rate $\mu$}, denoted a \emph{CS 4-fold},
if $\Nhat=N\setminus\{z_1,\ldots,z_s\}$ is a (nonsingular)
4-dimensional submanifold of $M$ and there exist constants
$0<\epsilon<\epsilon_M$ and $\mu\in(1,2)$, a compact 3-dimensional 
Riemannian submanifold $(L_i,h_i)$ of
$\mathcal{S}^6\subseteq\R^7$, where $h_i$ is the restriction of the
round metric on $\mathcal{S}^6$ to $L_i$, an open set $U_i\ni
z_i$ in $N$ with $U_i\subseteq V_i$ and a smooth map
$\Phi_i:(0,\epsilon)\times L_i\rightarrow
B(0;\epsilon_M)\subseteq\R^7$, for $i=1,\ldots,s$, such that
 $\Psi_{i}=\chi_i\circ\Phi_i:(0,\epsilon)\times L_i\rightarrow
 U_i\setminus\{z_i\}$ is a diffeomorphism, and
 $\Phi_i$ satisfies 
\begin{equation}\label{csnormaleq}
\Phi_i(r_i,x_i)-\iota_i(r_i,x_i)\in \big(T_{r_ix_i}\iota_i(C_i)\big)^{\perp}\qquad 
\text{for all $(r_i,x_i)\in (0,\epsilon)\times L_i$}
\end{equation} and
\begin{equation}\label{csasymeq}
\big|\nabla^j_i\big(\Phi_i(r_i,x_i)-\iota_i(r_i,x_i)\big)\big|
=O\big(r_i^{\mu-j}\big)\qquad \text{for $j\in\N$ as $r_i\rightarrow 0$,}
\end{equation}
where $\iota_i(r_i,x_i)=r_i x_i\in B(0;\epsilon_M)$, $\nabla_i$
is the Levi-Civita connection of the cone metric
$g_i=dr_i^2+r_i^2h_i$ on $C_i=(0,\infty)\times L_i$ coupled with
partial differentiation on $\R^7$, and
 $|.|$ is calculated with respect to $g_i$.

 We call $C_i$ the \emph{cone} at the singularity $z_i$ and
$L_i$ the \emph{link} of the cone $C_i$.  We may write $N$ as
the disjoint union
$N=K_N\sqcup\bigsqcup_{i=1}^sU_i,$
where $K_N$ is compact.

If $\Nhat$ is coassociative in $M$, we say that $N$ is a \emph{CS
coassociative 4-fold}.
\end{dfn}

\begin{remark} If $N$ is a CS 4-fold, $\Nhat$ is \emph{non-compact}.
\end{remark}

Suppose $N$ is a CS 4-fold at $z_1,\ldots,z_s$ with rate $\mu$
in $(M,\varphi,g_{\varphi})$ and use the notation of Definition
\ref{csdfn}. The induced metric on $\Nhat$, $g_{\varphi}|_{\Nhat}$, makes
$\Nhat$ into a Riemannian manifold. Moreover, it is clear from
\eq{csasymeq} that, as long as $\mu<2$, the maps $\Psi_i$ satisfy
 \begin{equation}\label{csmetriceq}
\big|\nabla_i^j\big(\Psi_i^*(g_{\varphi}|_{\Nhat})-g_i\big)\big|=O\big(r_i^{\mu-1-j}\big)\qquad
\text{for $j\in\N$ as $r_i\rightarrow 0$.}
\end{equation}
Consequently, the condition $\mu>1$ guarantees that the induced
metric on $\hat{N}$ genuinely converges to the conical metric on
$C_i$.

\begin{note} As shown on \cite[p.~6]{Lotaycs}, since $\mu\in(1,2)$, Definition \ref{csdfn} is
independent of the choice of $\GG_2$ coordinate system near the
singularities, up to equivalence.\end{note}

\begin{dfn}\label{radiusfndfn}
Let $N$ be a CS coassociative 4-fold in an almost $\GG_2$ manifold and use the notation of Definition 
\ref{csdfn}. 
A \emph{radius function} on $\hat{N}$ is
a smooth map $\rho:\hat{N}\rightarrow (0,1]$ such that there exist positive constants
$c_1<1$ and $c_2>1$ with $c_1r_i<\Psi_i^*(\rho)<c_2r_i$ on
$(0,\epsilon)\times L_i$ for $i=1,\ldots,s$.  
\end{dfn}

\noindent It is clear how we may construct such a function.

We now make a definition which also depends only on equivalence
classes of $\GG_2$ coordinate systems near the singularities.

\begin{dfn}\label{cstgtconedfn}
Let $N$ be a CS coassociative 4-fold at $z_1,\ldots,z_s$ in an almost $\GG_2$
manifold.  Use the notation of Definitions
\ref{G2coordsdfn} and \ref{csdfn}. For $i=1,\ldots,s$, define a
cone $\hat{C}_i$ in $T_{z_i}M$ by
$\hat{C}_i=(\zeta_i\circ\iota_i)(C_i)$. We call $\hat{C}_i$ the
\emph{tangent cone} at $z_i$.
\end{dfn}

\noindent  Using \eq{csasymeq}, one sees 
 that $\hat{C}_i$ is a tangent cone at $z_i$ in the sense of 
Definition \ref{tgtconedfn}.  Since the tangent cone has multiplicity one, \cite[Theorem 5.7]{Simon} implies that $\hat{C}_i$ is 
the unique tangent cone to $N$ at $z_i$.  It is still an open question whether a general calibrated integral current has a unique 
tangent cone at each point.

   We conclude with a straightforward result which follows from Proposition \ref{caltgtconeprop} or, by more elementary means, from 
\cite[Proposition 3.6]{Lotaycs}.

\begin{prop}\label{coasstgtconeprop}
Let $N$ be a CS coassociative 4-fold at $z_1,\ldots,z_s$ in an almost
$\GG_2$ manifold.  The tangent cones at
$z_1,\ldots,z_s$ are coassociative.
\end{prop}

\subsection{Weighted Banach spaces}

For this subsection we let $(P,g)$ denote the nonsingular part of a CS coassociative 4-fold in an almost $\GG_2$ manifold with the induced metric
 and let $\rho$ be a radius function on $P$, 
as given in Definitions \ref{csdfn} and \ref{radiusfndfn}.  
We also let $\nabla$ denote the Levi-Civita connection of $g$.  We define \emph{weighted} Banach spaces of forms on $P$ as 
in \cite[$\S$1]{Bartnik} so as to implement the analytic framework of \cite{Lotaycs}.

We shall use the notation and definition of the usual `unweighted' Banach spaces of forms as  
in \cite[$\S$1.2]{Joyexcept}; that is, Sobolev and H\"older spaces are denoted by $L^p_k$ and $C^{k,\,a}$ respectively,
where $p\geq 1$, $k\in\N$ and $a\in(0,1)$.  
Recall that, by the Sobolev Embedding Theorem, $L^p_{k}$ embeds continuously in $L^q_l$ if $l\leq k$ and
 $l-\frac{4}{q}\leq k-\frac{4}{p}$,
and $L^p_k$ embeds continuously in $C^{l,\,a}$ if $k-\frac{4}{p}\geq l+a$.
We also introduce the notation $C^k_{\text{loc}}$ for the space of forms $\xi$ such that $f\xi$ lies in $C^k$ for every
smooth compactly supported function $f$, and similarly define spaces $L^p_{k,\,\text{loc}}$ and $C^{k,\,a}_{\text{loc}}$.

\begin{dfn}\label{wSobdfn}
Let $p\geq 1$, $k\in\N$ and $\lambda\in\R$. The \emph{weighted Sobolev space}
$L_{k,\,\lambda}^p(\Lambda^mT^*P)$ of $m$-forms $\xi$ on $P$ is the
subspace of
$L^p_{k,\,\text{loc}}(\Lambda^mT^*P)$ such that the norm
\begin{equation*}\label{ch6s2eq2}
\|\xi\|_{L_{k,\,\lambda}^p}=\sum_{j=0}^k\left(\int_{P}
|\rho^{j-\lambda}\nabla^j\xi|^p\rho^{-4} \,dV_g\right)^\frac{1}{p}
\end{equation*} is finite. 
Then $L_{k,\,\lambda}^p(\Lambda^mT^*P)$ is a
Banach space.  
\end{dfn}

We now define the \emph{dual} weighted Sobolev space which shall be useful later.

\begin{dfn}\label{dualdfn}  Use the notation from Definition \ref{wSobdfn}. 
Let $p,q>1$ be such that $\frac{1}{p}+\frac{1}{q}=1$, let $k,l\in\N$ and let $\lambda\in\R$.  Define
a pairing $\langle\,.\,,\,.\,\rangle:L^p_{k,\,\lambda}(\Lambda^mT^*P)\times L^q_{l,\,-4-\lambda}(\Lambda^mT^*P)\rightarrow\R$ by
$$\langle\xi,\eta\rangle=\int_P\xi\w\ast\eta.$$
We shall refer to this as the \emph{dual pairing}.  For our purposes, we take the dual space of
 $L^p_{k,\,\lambda}(\Lambda^mT^*P)$ to be $L^q_{l,\,-4-\lambda}(\Lambda^mT^*P)$, with linear functionals represented by dual pairings.
\end{dfn}

\begin{dfn}\label{wCkdfn} Let $\lambda\in\R$ and $k\in\N$. The
\emph{weighted $C^k$-space} $C_{\lambda}^{k}(\Lambda^mT^*P)$ of
$m$-forms $\xi$ on $P$ is the subspace of
$C^k_{\text{loc}}(\Lambda^mT^*P)$ such that the norm
$$\|\xi\|_{C_{\lambda}^{k}}=\sum_{j=0}^k
\sup_{P}|\rho^{j-\lambda}\nabla^j\xi|$$ is finite. 
We also define $C_{\lambda}^{\infty}(\Lambda^mT^*P)=\bigcap_{k\geq
0}C_{\lambda}^{k}(\Lambda^mT^*P)$.  
Then $C_{\lambda}^{k}(\Lambda^mT^*P)$
is a Banach space, but $C_{\lambda}^{\infty}(\Lambda^mT^*P)$ is not in general. 
\end{dfn}

\begin{dfn}\label{wHolderdfn}
Let $E$ be a vector bundle on $P$ endowed with Euclidean metrics on its fibres and a connection preserving these 
metrics.  Let $d(x,y)$ be the geodesic distance between points $x,y\in
P$, let $a\in
(0,1)$, let $k\in\N$  and let $\lambda\in\R$. Let
\begin{align*}H=\{&(x,y)\in P\times P\,:\,x\neq
y,\,c_1\rho(x)\leq\rho(y)\leq c_2\rho(x)\,\;\text{and}\;\,\\
&\text{there exists a geodesic in $P$ of length $d(x,y)$ from
$x$ to $y$}\},\end{align*} where $0<c_1<1<c_2$ are constant. A
section $s$ of $E$ is \emph{H\"older continuous} (with \emph{exponent $a$})
if
$$[s]^a=
\sup_{(x,y)\in H}\frac{|s(x)-s(y)|_{E}}{d(x,y)^a}<\infty.$$ We
understand the quantity $|s(x)-s(y)|_E$ as follows.  Given $(x,y)\in
H$, there exists a geodesic $\gamma$ of length $d(x,y)$ connecting
$x$ and $y$. Parallel translation along $\gamma$ using the
connection on $E$ identifies the fibres over $x$ and $y$ and the
metrics on them. Thus, with this identification, $|s(x)-s(y)|_E$ is
well-defined.

The \emph{weighted H\"older space}
$C_{\lambda}^{k,\,a}(\Lambda^mT^*P)$ is the subspace of
$C^{k,\,a}_{\text{loc}}(\Lambda^mT^*P)$ such that the norm
$$\|\xi\|_{C^{k,\,a}_{\lambda}}=\|\xi\|_{C^{k}_{\lambda}}+[\xi]^{k,\,a}_{\lambda}$$
is finite, where
$$[\xi]^{k,\,a}_{\lambda}=[\rho^{k+a-\lambda}\nabla^k\xi]^{a}.
$$  Then $C_{\lambda}^{k,\,a}(\Lambda^mT^*P)$ is a
Banach space.  It is clear that we have embeddings 
$C_{\lambda}^{k,\,a}(\Lambda^mT^*P)\hookrightarrow
C_{\lambda}^l(\Lambda^mT^*P)$ and $C_{\lambda}^{k+1}(\Lambda^mT^*P)\hookrightarrow C_{\lambda}^{l,\,a}(\Lambda^mT^*P)$ if $l\leq k$.
\end{dfn}

\begin{remark}
The set $H$ in Definition \ref{wHolderdfn} is introduced so that $[\xi]^{k,\,a}_{\lambda}$ is well-defined.
\end{remark}

Finally, we shall need the analogue of the Sobolev Embedding Theorem for
weighted spaces, which is adapted from \cite[Lemma
7.2]{LockhartMcOwen} and \cite[Theorem 1.2]{Bartnik}.

\begin{thm}[\textbf{Weighted Sobolev Embedding Theorem}]
\label{wembedthm} Let $p,\,q\geq 1$,\\ $a\in (0,1)$,
$\lambda,\nu\in\R$ and $k,l\in\N$.
\begin{itemize}
\item[{\rm (a)}] If $k\geq l$, $k-\frac{4}{p}\geq l-\frac{4}{q}$,
and either $p\leq q$ and $\lambda\geq\nu$, or $p>q$ and $\lambda>\nu$,
there is a continuous embedding
$L_{k,\,\lambda}^p(\Lambda^mT^*P)\hookrightarrow
L_{l,\,\nu}^q(\Lambda^mT^*P)$. \item[{\rm (b)}] If
$k-\frac{4}{p}\geq l+a$, there is a continuous embedding
$L_{k,\,\lambda}^p(\Lambda^mT^*P)\hookrightarrow
C_{\lambda}^{l,\,a}(\Lambda^mT^*P)$.
\end{itemize}
\end{thm}

\subsection{Deformation theory}

We now review and discuss the key deformation theory results for CS coassociative 4-folds from \cite{Lotaycs}.  We begin by recalling
 the linear differential operator governing infinitesimal deformations of CS coassociative 4-folds.

\begin{dfn}\label{dd*dfn}  
Let $N$ be a CS coassociative 4-fold in an almost $\GG_2$ manifold.  Let $p>4$, $k\geq 2$ and let $\lambda\in\R$.  Define 
\begin{equation}\label{dd*eq}
(\d_++\d^*)_{\lambda}:L_{k+1,\,\lambda}^{p}(\Lambda_+^2T^*\hat{N}\oplus\Lambda^4T^*\hat{N})\rightarrow 
L_{k,\,\lambda-1}^{p}(\Lambda^3T^*\hat{N})
\end{equation}
by $(\d_++\d^*)_{\lambda}(\alpha,\beta)=\d\alpha+\d^*\beta$.
\end{dfn}

\begin{notes}  We use the operator \eq{dd*eq} rather than simply the exterior derivative on self-dual 2-forms since the former 
operator is elliptic whereas the latter is not.
The choice of $p>4$ and $k\geq 2$ ensures that $L^p_{k+1,\,\lambda}$-solutions to the (nonlinear) deformation problem 
are in fact smooth.
\end{notes}

\begin{dfn}\label{Ddfn}
Let $N$ be a CS coassociative 4-fold in an almost $\GG_2$ manifold and use the notation of Definition \ref{csdfn}.  Let
\begin{align*}
D(\lambda,i)=\{(\alpha_i,\beta_i)\in C^\infty(\Lambda^2T^*L_i&\oplus
\Lambda^3T^*L_i)\,:\,\\&\d_i\alpha_i=-\lambda\beta_i,\;\d_i\!*_i\!\alpha_i+\d^*_i\beta_i=
-(\lambda+2)\alpha_i\},
\end{align*}
where $*_i$, $\d_i$ and $\d^*_i$ are the Hodge star, the exterior derivative and its formal adjoint on $L_i$.
By \cite[Propositions 5.1 \& 5.2]{Lotaycs}, the set $\mathcal{D}$ of real numbers 
 such that \eq{dd*eq} is not Fredholm is countable and discrete and given by
\begin{equation*}
\mathcal{D}=\bigcup_{i=1}^s\{\lambda\,:\,D(\lambda,i)\neq 0\}.
\end{equation*}
We also set $d(\lambda)=\sum_{i=1}^s\dim D(\lambda,i)$.
\end{dfn}

\begin{note}
The forms $(\alpha_i,\beta_i)\in D(\lambda,i)$ correspond to homogeneous forms
 $(\alpha,\beta)\in C^{\infty}(\Lambda^2_+T^*C_i\oplus\Lambda^4T^*C_i)$ of order $\lambda$ which satisfy 
$\d\alpha+\d^*\beta=0$ on the cone $C_i$. 
\end{note}

In \cite{Lotaycs} the author studied deformations of CS coassociative 4-folds, allowing the singularities and 
tangent cones to vary, and permitting changes in the ambient $\GG_2$ structure.  However, for our purposes, we require 
a slightly more general theory, where we allow the underlying cones on which the singularities are modelled to vary
as well.  We now define the moduli space of deformations.

\begin{dfn}\label{modulidfn}
Let $N$ be a CS coassociative 4-fold in an almost $\GG_2$ manifold $(M,\varphi,g_{\varphi})$ and let
 $\mathcal{D}$ be given by Definition \ref{Ddfn}.  Suppose further that $N$
has singularities at $z_1,\ldots,z_s$ with rate $\mu\in(1,2)\setminus\mathcal{D}$, having 
cone $C_i$ and tangent cone $\hat{C}_i$ at $z_i$ for all $i$.  
Let $\mathfrak{C}=\prod_{i=1}^s\mathcal{C}_i$ where 
 $\mathcal{C}_i$ is a smooth, connected family of coassociative cones in $\R^7$, 
 closed under the natural action of 
$\GG_2\ltimes\R^7$, such that $C_i\in\mathcal{C}_i$ for all $i$.  

The moduli space of deformations $\mathcal{M}(N,\mu,\mathfrak{C})$ is the set of CS coassociative 
4-folds $N^{\prime}$ in $M$ such that:
\begin{itemize}
\item[(a)] $N^{\prime}$ has a singularity at $z_i^{\prime}$ with rate $\mu$ and cone 
in $\mathcal{C}_i$ for $i=1,\ldots,s$; and 
\item[(b)] there exists a diffeomorphism $h:M\rightarrow M$, isotopic to the identity, such that $h(z_i)=z_i^{\prime}$ for all $i$, $h|_N:N\rightarrow N^{\prime}$ is a homeomorphism and $h|_{\hat{N}}:\hat{N}\rightarrow N^{\prime}\setminus\{z_1^{\prime},\ldots,z_s^{\prime}\}$ is a diffeomorphism.
\end{itemize}
\end{dfn}

We state the deformation theory result we require.

\begin{thm}\label{maindefthm} Use the notation of Definition \ref{modulidfn}, let $p>4$ and $k\geq 2$.  
There exist finite-dimensional vector spaces $\mathcal{I}(N,\mu,\mathfrak{C})$ and 
 $\mathcal{O}(N,\mu,\mathfrak{C})$, with $\mathcal{O}(N,\mu,\mathfrak{C})$ contained in $L^p_{k,\,\mu-1}(\Lambda^3T^*\Nhat)$,
 and there exist 
\begin{itemize}
\item[\emph{(a)}] a smooth manifold $\hat{\mathcal{M}}(N,\mu,\mathfrak{C})$,
which is an open neighbourhood of $0$ in
$\mathcal{I}(N,\mu,\mathfrak{C})$, and 
\item[\emph{(b)}] a smooth map
$\pi:\hat{\mathcal{M}}(N,\mu,\mathfrak{C})\rightarrow\mathcal{O}(N,\mu,\mathfrak{C})$,
with $\pi(0)=0$, 
\end{itemize}
such that an open neighbourhood of zero in
$\Ker\pi$ is homeomorphic to an open neighbourhood
of $N$ in $\mathcal{M}(N,\mu,\mathfrak{C})$.

Furthermore, if $\mathcal{O}(N,\mu,\mathfrak{C})=\{0\}$, then $\mathcal{M}(N,\mu,\mathfrak{C})$ is a smooth manifold of dimension equal to 
that of $\mathcal{I}(N,\mu,\mathfrak{C})$.
\end{thm}

We can actually say more about the spaces $\mathcal{I}(N,\mu,\mathfrak{C})$ and $\mathcal{O}(N,\mu,\mathfrak{C})$.

\begin{prop}\label{maindefprop}
Use the notation of Definitions \ref{dd*dfn}-\ref{modulidfn} and Theorem \ref{maindefthm}.
\begin{itemize}
\item[\emph{(a)}] $\mathcal{I}(N,\mu,\mathfrak{C})$ contains a subspace isomorphic to $\Ker(\d_++\d^*)_{\mu}$.
\item[\emph{(b)}] $\mathcal{O}(N,\mu,\mathfrak{C})$ is transverse to the space $\Image (\d_++\d^*)_{\mu}$, it is contained in 
$\overline{\d(L^p_{k+1,\,\mu}(\Lambda^2T^*\Nhat))}\subseteq L^p_{k,\,\mu-1}(\Lambda^3T^*\Nhat)$ and satisfies
\begin{equation}\label{obseq1}\dim\mathcal{O}(N,\mu,\mathfrak{C})\leq \!\!\!\!\!\sum_{\lambda\in(-2,\mu)\cap\mathcal{D}}
\!\!\!\! d(\lambda)-\sum_{i=1}^s\dim\mathcal{C}_i.
\end{equation}
\end{itemize}
\end{prop}

These results essentially follow from \cite[Theorem 7.9 \& Proposition 8.10]{Lotaycs}, the only difference being that we  
allow the cones $C_i$ on which the singularities are modelled to deform in the families $\mathcal{C}_i$,  
 which may be larger than the families given simply by translations and $\GG_2$ transformations of $C_i$.  
Rather than repeating the entire analysis in \cite{Lotaycs} with this change, we appeal to the similar discussion in \cite[$\S$8.3]{JoyceSLsing2}; that
 is, the \emph{infinitesimal deformation space} $\mathcal{I}(N,\mu,\mathfrak{C})$ is unchanged but the dimension of the 
\emph{obstruction space} $\mathcal{O}(N,\mu,\mathfrak{C})$ is reduced by the dimension of the families in which the cones vary.  

\begin{notes}
 In \eq{obseq1}, the sum is over $(-2,\mu)\cap\mathcal{D}$, rather than $(-1,\mu)\cap\mathcal{D}$ as in \cite[Proposition 8.10]{Lotaycs}. This ability to improve to the latter smaller set was based on the claim in 
 \cite[Proposition 5.3]{Lotaycs}, which followed from applying \cite[Theorem 10.2]{LockhartMcOwen}, that $(-2,-1]\cap\mathcal{D}=
 \emptyset$. Though $(-2,-1)\cap\mathcal{D}=\emptyset$ for almost all examples of interest, often $-1\in\mathcal{D}$ so 
the claim is erroneous.    
However, we shall show in Proposition \ref{improvdimprop} that not only is \cite[Proposition 8.10]{Lotaycs} valid, but the estimate in \eq{obseq1} can be improved even further. 
\end{notes}

In \cite[Theorem 7.13]{Lotaycs} we proved a result for deformations of $N$ where the ambient $\GG_2$ structure on $M$ also varies. 
We see that we have to restrict our choice of perturbations of the $\GG_2$ structure as follows.

\begin{dfn}\label{G2familydfn}  Use the notation of Definition \ref{modulidfn}.  
By \cite[Proposition 6.19]{Lotaycs} there exists a neighbourhood $T_N$ of $N$ in $M$ which retracts onto $N$ and $H^3_{\text{dR}}(T_N)\cong H^3_{\text{cs}}(\hat{N})$, the 
third compactly supported cohomology group of $\hat{N}$.  Thus, any closed positive 3-form $\varphi^{\prime}$ on $M$ defines a cohomology class $[\varphi^{\prime}|_{\hat{N}}]\in H^3_{\text{cs}}(\hat{N})$ and, moreover, if $[\varphi^{\prime}|_{\hat{N}}]\neq 0$ in $H^3_{\text{cs}}(\hat{N})$ there are \emph{no} nearby coassociative deformations of $N$.  Therefore, we let
$\mathcal{F}=\{(\varphi^f,g_{\varphi^f})\,:\,f\in B(0;\delta)\subseteq\R^r\}$ be a smooth $r$-dimensional family of closed $\GG_2$ structures on $M$, with $\varphi^{0}=\varphi$ and 
$[\varphi^{f}|_{\hat{N}}]=0$ for all $f\in B(0;\delta)$.
\end{dfn}

The obstruction space for the deformation problem where the ambient $\GG_2$ structure deforms in the family $\mathcal{F}$, suitably generalized to include the possibility that the cones at the singularities 
vary in the family $\mathfrak{C}$, is contained in $\mathcal{O}(N,\mu,\mathfrak{C})$.  Therefore, if
 $\mathcal{O}(N,\mu,\mathfrak{C})=\{0\}$, for any sufficiently small perturbation of the $\GG_2$ structure $(\varphi,g_{\varphi})$ 
in $\mathcal{F}$ we obtain a corresponding CS deformation of $N$ which is coassociative with respect to the new $\GG_2$ structure.  
Thus, we have the following result using \cite[Theorem 7.13]{Lotaycs}.
 
\begin{thm}\label{G2deformthm}  Use the notation of Definitions \ref{modulidfn} and \ref{G2familydfn} and Theorem \ref{maindefthm}.  
If\/ $\mathcal{O}(N,\mu,\mathfrak{C})=\{0\}$, there exists $\delta_N\in(0,\delta)$ such that, if $f\in B(0;\delta_N)$, 
then there exists a CS deformation $N^f$ 
of $N$ in $T_N$ with $s$ singularities with rate $\mu$, modelled on cones in $\mathfrak{C}$, which is coassociative with respect to $\varphi^f$.
\end{thm}

\begin{note}
 By the method of proof of \cite[Theorem 7.13]{Lotaycs}, a `sufficiently small' perturbation $(\varphi^{\prime},
g_{\varphi^{\prime}})$ of $(\varphi,g_\varphi)$ is one for which $\|\varphi^{\prime}-\varphi\|_{C^1}<\epsilon_{C^1}$ and 
$\|\varphi^{\prime}-\varphi\|_{L^p_2}<\epsilon_{L^p_2}$ for some $p>4$, where the norms are calculated in $T_N$ with respect to 
$g_{\varphi}$, and $\epsilon_{C^1}$ and $\epsilon_{L^p_2}$ are small constants determined by the geometry of $(M,\varphi,g_{\varphi})$
 near $N$ (i.e.~in $T_N$).  
Thus $\delta_N$ is chosen so that these $C^1$ and $L^p_2$ norms are smaller than the appropriate constants.  
\end{note}

\subsection{Multiplicity one tangent cones}

In this subsection we study coassociative integral currents with multiplicity one tangent cones, motivated by the work on 
special Lagrangian integral currents in \cite[$\S$6]{JoyceSLsing1}.  The key condition is that 
the underlying cone is \emph{Jacobi integrable}, which we define for a coassociative cone following \cite[Definition 6.7]{JoyceSLsing1}.

\begin{dfn}\label{Jacdfn}
Let $C$ be a coassociative cone in $\R^7$ with compact link $L$ in $\mathcal{S}^6$ such that $C\setminus\{0\}$ is nonsingular.  We say that
 $v\in C^{\infty}(\nu_{\mathcal{S}^6}(L))$ is a \emph{Lagrangian Jacobi field} if 
$\alpha_v=*(v\lrcorner\varphi_0)|_{TL}\in C^{\infty}(T^*L)$ satisfies $\d\alpha_v=-3*\alpha_v$, where $*$ is the Hodge 
star on $L$.  Notice that $v$ is a Lagrangian  
Jacobi field if and only if $r^3*\alpha_v+r^2\alpha_v\w\d r$ is a 
closed self-dual 2-form on $C$ of order $O(r)$, and so defines an infinitesimal deformation of $C$ as a coassociative cone; i.e.~an
infinitesimal deformation of $L$ as a Lagrangian in $\mathcal{S}^6$ in the sense of Definition \ref{Lagdef}.  Since Lagrangians in $\mathcal{S}^6$ are minimal, Lagrangian Jacobi fields are Jacobi fields in the usual sense.

We say that $C$ is \emph{Jacobi integrable} if every Lagrangian 
Jacobi field $v$ on $L$ defines a smooth 
one-parameter family $\{L_t=\exp_{tv}(L)\subseteq\mathcal{S}^6\,:\,t\in(-\tau,\tau)\}$, for some $\tau>0$, of Lagrangian submanifolds 
of $\mathcal{S}^6$.
\end{dfn}

Our next result proves Theorem \ref{integralthm}; namely, that interior singular points of 
 coassociative integral currents, with multiplicity one tangent cones modelled on Jacobi integrable cones, are conical 
 singularities in the sense of Definition \ref{csdfn}.

\begin{thm}\label{cosingthm}
Let $N$ be a coassociative integral current in an almost $\GG_2$ manifold $M$ and let $z\in N^\circ$ be a singular point of $N$.  Let $\{\chi:B(0;\epsilon_M)
\rightarrow V\}$ be a $\GG_2$ coordinate system near $z$ in the sense of Definition \ref{G2coordsdfn}, with $\zeta=\d\chi|_0$.  Suppose that $C$ is a cone in $\R^7$ with
 compact link $L$ such that $C\setminus\{0\}$ is nonsingular and $\hat{C}=\zeta(C)\subseteq T_zM$ is a multiplicity one tangent cone for $N$ at $z$.   Then $C$ is coassociative and
 $\hat{C}$ is the unique tangent cone for $N$ at $z$.  
 
Suppose further that $C$ is Jacobi integrable in the sense of Definition \ref{Jacdfn} and let $U=N\cap V$.  Then there exist
 $\epsilon\in(0,\epsilon_M)$ and an embedding $\Phi:(0,\epsilon)\times L\rightarrow B(0;\epsilon_M)$ such that 
 $U\setminus\{z\}= \chi\circ\Phi((0,\epsilon)\times L)$ as an embedded submanifold of $M$ and $\Phi$ satisfies 
 \eq{csnormaleq} and \eq{csasymeq} for some $\mu\in(1,2)$.
\end{thm}

\begin{proof}
The fact that $C$ is coassociative follows from Proposition \ref{caltgtconeprop}.  By \cite[Theorem 5.7]{Simon}, the tangent cone 
$\hat{C}$ is unique and $U\setminus\{z\}$ can be realized as a $C^2$-embedding $\Psi$ of $(0,\epsilon)\times L$, for some $\epsilon>0$.  Moreover, $\Psi=\chi\circ\Phi$ where $\Phi$ is a $C^2$-embedding of $(0,\epsilon)\times L$ into $B(0;\epsilon_M)$ which satisfies \eq{csnormaleq}. 
 This last point is not explicit in the statement of the theorem, but does follow from the proof.  Moreover, if ~$\iota:C\setminus\{0\}\cong(0,\infty)\times L\rightarrow\R^7$ is the inclusion map, $\Phi$ also satisfies:
\begin{align}
\big|\Phi(r,x)-\iota(r,x)\big|&=o(r);\label{asymeq0}\\
\big|\nabla\big(\Phi(r,x)-\iota(r,x)\big)\big|&=o(1);\label{asymeq1}\\
\big|\nabla^2\big(\Phi(r,x)-\iota(r,x)\big)\big|&=O(1)\label{asymeq2}
\end{align}
as $r\rightarrow 0$.  Equation \eq{asymeq2} again is not stated explicitly in \cite[Theorem 5.7]{Simon} but follows from the minimality of $C$ and 
$N$ as discussed in the proof of \cite[Theorem 1]{AdamsSimons}.  The aforementioned theorem \cite[Theorem 1]{AdamsSimons} and its
 proof imply that there exists $\mu\in(1,2)$ such that the estimates \eq{asymeq0} and \eq{asymeq1} can be improved to:
\begin{align}
\big|\Phi(r,x)-\iota(r,x)\big|&=O(r^{\mu});\label{asymeq0a}\\
\big|\nabla\big(\Phi(r,x)-\iota(r,x)\big)\big|&=O(r^{\mu-1})\label{asymeq1a}
\end{align}
as $r\rightarrow 0$, since $C$ is Jacobi integrable.  We now only need a regularity argument to complete the proof, 
similar to \cite[Proposition 4.17]{Lotayac}, which we now detail.

Let $P=\iota\big((0,\epsilon)\times L\big)\subseteq\R^7$ and recall that, since $P$ is coassociative, we have an isomorphism $\jmath_P:\nu(P)\rightarrow\Lambda^2_+T^*P$ by 
Proposition \ref{jmathprop}.  For $\alpha\in C^2_{\text{loc}}(\Lambda^2_+T^*P)$ 
and $v=\jmath_P^{-1}(\alpha)$, define 
\begin{align*}
F(\alpha)&=\exp_v^*\left(\varphi_0|_{\exp_v(P)}\right)\in C^1_{\text{loc}}(\Lambda^3T^*P)\\ 
\intertext{and}
G(\alpha)&=\pi_{\Lambda^2_+}\big(\d^*F(\alpha)\big)\in C^0_{\text{loc}}(\Lambda^2_+T^*P),
\end{align*}
where $\pi_{\Lambda^2_+}$ is the projection from 2-forms to self-dual 2-forms on $P$.  Clearly, if $F(\alpha)=0$, the 
deformation $\exp_v(P)$ of $P$ is coassociative and $G(\alpha)=0$.
Moreover, $\d F|_0(\alpha)=\d\alpha$ by \cite[p.~731]{McLean} and thus $\d G|_0(\alpha)=\d^*_+\d\alpha$, where $\d^*_+=\pi_{\Lambda^2_+}\circ\d^*$ acting on 3-forms on $P$.  We deduce 
that $G(\alpha)=0$ is a nonlinear \emph{elliptic} equation for $\alpha$ at $0$; that is, its linearisation at $0$ is elliptic.

  Since $\Phi$ is a $C^2$-coassociative embedding satisfying \eq{csnormaleq}, it defines a $C^2$-normal vector field $v_{\Phi}$ on 
$P$ and hence a $C^2$-self-dual 2-form $\alpha_{\Phi}=\jmath_P(v_\Phi)$ on $P$ which satisfies $F(\alpha_\Phi)=0$.  Moreover, 
$\alpha_\Phi\in C^2_{\mu}(\Lambda^2_+T^*P)$ by \eq{asymeq0a}-\eq{asymeq1a} and \eq{asymeq2}, since $\mu-2<0$.  Thus, 
$\alpha_{\Phi}\in C^{1,\,a}_{\mu}(\Lambda^2_+T^*P)$ for any $a\in(0,1)$.

Let $\nabla$ be the Levi-Civita connection of the conical metric on $P$.  We can write, for $\alpha\in C^2_{\text{loc}}(\Lambda^2_+T^*P)$ and $p\in P$, 
$$G(\alpha)(p)=R\big(p,\alpha(p),\nabla\alpha(p)\big)\nabla^2\alpha(p)+E\big(p,\alpha(p),\nabla\alpha(p)\big),$$
where $R$ and $E$ are smooth functions of their arguments, since $G(\alpha)$ is linear in $\nabla^2\alpha$ and a smooth function of
$\alpha$.  This leads us to define a new operator on $\beta\in C^2_{\text{loc}}(\Lambda^2_+T^*P)$ by
$$S_\Phi(\beta)(p)=R\big(p,\alpha_\Phi(p),\nabla\alpha_\Phi(p)\big)\nabla^2\beta(p).$$
Note that $S_\Phi$ is \emph{not} the linearisation of $G$ of $0$, but it is a \emph{linear second order elliptic} operator on $\beta$ with 
coefficients in $C^{0,\,a}_{\text{loc}}$.  

Since $G(\alpha_\Phi)=0$, $S_\Phi(\alpha_\Phi)=-E(\alpha_\Phi)$, where we set  
 $$E(\beta)(p)=E\big(p,\beta(p),\nabla\beta(p)\big)$$ 
for $p\in P$ and $\beta\in C^2_{\text{loc}}(\Lambda^2_+T^*P)$.  
As argued in \cite[Proposition 4.17]{Lotayac}, $E(\beta)\in C^{k,\,a}_{\mu-1}$ if $\beta\in C^{k+1,\,a}_{\mu}$ since $E$ is no worse
 than quadratic in $\beta$ and $\nabla\beta$.  Thus, using the Schauder regularity estimates as given in \cite{Mazya} and the 
fact that $\alpha_\Phi\in C^{1,\,a}_{\mu}$, we deduce that $\alpha_\Phi\in C^{2,\,a}_{\mu}(\Lambda^2_+T^*P)$.  A standard inductive argument leads us to the conclusion that $\alpha_\Phi\in C^{k,\,a}_{\mu}$ for all $k\in\N$.  Thus the corresponding 
map $\Phi$ satisfies \eq{csasymeq} as claimed.
\end{proof}

From Definition \ref{csdfn} we have an immediate corollary to Theorem \ref{cosingthm}.

\begin{cor}\label{cosingcor}
Let $N$ be a connected coassociative integral current in an almost $\GG_2$ manifold with $\partial N=\emptyset$.  If 
$N$ has multiplicity one tangent cones at its singular points modelled on Jacobi integrable coassociative cones, then $N$ is a CS coassociative 4-fold in the sense of Definition \ref{csdfn}.
\end{cor}

\section{Coassociative cones and stability index}\label{stabindsec}

In this section we define the notion of \emph{stability} of conical singularities of coassociative 4-folds using a 
numerical invariant for a coassociative cone which we call the \emph{stability index}.  This stability index is calculated 
using the spectrum of the curl operator acting on 1-forms on the link of the cone.  We begin by giving a brief 
survey of the known examples of coassociative cones.  

\subsection{Examples}\label{orbitsec}

For this subsection it is convenient to identify $\R^7$ with the imaginary octonions $\Im\O$.  Let $\{\eps_1,\ldots,\eps_7\}$ 
be a basis for $\Im\O$ satisfying the multiplication law below.
\begin{equation*}
\begin{gathered}
\begin{array}{rrrrrrrrrr}
  & \vline &  \eps_1 &  \eps_2 &  \eps_3 &  \eps_4 &  \eps_5 &  \eps_6 &  \eps_7 \\
\hline
 \eps_1 & \vline  &  -1  &  \eps_3 & -\eps_2 &  \eps_5 & -\eps_4 &  \eps_7 & -\eps_6 \\
\eps_2 & \vline  & -\eps_3 &  -1  &  \eps_1 &  \eps_6 & -\eps_7 & -\eps_4 &  \eps_5 \\
\eps_3 & \vline  &  \eps_2 & -\eps_1 &  -1  & -\eps_7 & -\eps_6 &  \eps_5 &  \eps_4 \\
\eps_4 & \vline  & -\eps_5 & -\eps_6 &  \eps_7 &  -1  &  \eps_1 &  \eps_2 & -\eps_3 \\
\eps_5 & \vline  &  \eps_4 &  \eps_7 &  \eps_6 & -\eps_1 &  -1  & -\eps_3 & -\eps_2 \\
\eps_6 & \vline  & -\eps_7 &  \eps_4 & -\eps_5 & -\eps_2 &  \eps_3 &  -1  &  \eps_1 \\
\eps_7 & \vline  &  \eps_6 & -\eps_5 & -\eps_4 &  \eps_3 &  \eps_2 & -\eps_1 &  -1
\end{array}
\end{gathered}
\end{equation*}
If we identify $\Im\O\cong\R^7$ so that $\{\eps_1,\ldots,\eps_7\}$ corresponds to the standard oriented orthonormal 
basis in $\R^7$, then the multiplication law we have chosen is consistent with the $\GG_2$ structure on $\R^7$ given by 
$\varphi_0$ in \eq{phieq}.

\medskip

 The coassociative cones invariant under a closed 3-dimensional Lie subgroup of $\GG_2$ were 
classified in \cite{Mashimo}.  We now review these examples starting with the degenerate example of a coassociative cone, namely a coassociative 4-plane.

\begin{ex}\textbf{(Planes)} \label{planeex}
By \cite[Proposition 12.1.2]{JoyceRiem}, $\GG_2$ acts transitively on the set of coassociative 4-planes with isotropy $\SO(4)$.  Therefore, any 
coassociative 4-plane is equivalent up to $\GG_2$ transformation to
\begin{equation*}
C_0=\{x_0\eps_1+x_1\eps_3+x_2\eps_5+x_3\eps_7\in\Im\O\,:\,(x_0,x_1,x_2,x_3)\in\R^4\}.
\end{equation*}
The link $L_0$ of $C_0$ is trivially a totally geodesic $\SO(4)$-invariant $\mathcal{S}^3$ in $\mathcal{S}^6$.
\end{ex}

Our next example was first introduced in \cite[$\S$7]{LO} and shown to be coassociative in \cite[Theorem IV.3.2]{HarLaw}.

\begin{ex}\textbf{({\boldmath $\U(2)$} symmetry)} \label{u2ex} The cone $C_1$ in $\Im\O$, with link
$$L_1=\left\{\textstyle\frac{\sqrt{5}}{3}\,\bar{q}\eps_1 q+\textstyle\frac{2}{3}\,q\eps_5\,:\,\text{$q\in\langle 1,\eps_1,
\eps_2,\eps_3\rangle_{\R}$ with $|q|=1$}\right\},$$
is coassociative and invariant under a $\U(2)$ subgroup of $\GG_2$.  The Lagrangian $L_1$ is realized as an $\Sp(1)$-orbit in 
$\mathcal{S}^6$.  
\end{ex}

Recall that, by Corollary \ref{coasscor}, any complex 2-dimensional cone in $\C^3$ embedded in
 $\R^7$ is coassociative.

\begin{ex}\textbf{(Complex {\boldmath $\SO(3)$} symmetry)}\label{cxso3ex} Identify $\R^7\cong\R\oplus\C^3$ as 
in Lemma \ref{phisplitlem}.  
The cone
$$C_2=\{(0,z_1,z_2,z_3)\in\R\oplus\C^3\,:\,z_1^2+z_2^2+z_3^2=0\}$$
is coassociative and invariant under the standard $\SO(3)$ action on $\C^3$.  The real link $L_2$ of $C_2$ in $\mathcal{S}^6$
 is diffeomorphic to $\SO(3)$ and is the Hopf lift of the constant curvature degree 2 $\C\P^1$ in $\C\P^2$.
\end{ex}

Our final symmetric examples are most easily described using homogeneous harmonic cubics on $\R^3$.

\begin{ex}\label{so3ex}\textbf{({\boldmath $\SO(3)$} symmetry)} 
Identify $\Im\O$ with the homogeneous harmonic cubics $\mathcal{H}^3(\R^3)$ on $\R^3$ by:
\begin{align*}
\eps_1&\mapsto \textstyle\frac{\sqrt{10}}{10}\,x(2x^2-3y^2-3z^2);&&\displaybreak[0]\\
\eps_2&\mapsto -\sqrt{6}xyz;&
\eps_3&\mapsto \textstyle\frac{\sqrt{6}}{2}\,x(y^2-z^2);\displaybreak[0]\\
\eps_4&\mapsto -\textstyle\frac{\sqrt{15}}{10}\, y(4x^2-y^2-z^2);&
\eps_5&\mapsto -\textstyle\frac{\sqrt{15}}{10}\,z(4x^2-y^2-z^2);\displaybreak[0]\\
\eps_6&\mapsto \textstyle\frac{1}{2}\,y(y^2-3z^2);&
\eps_7&\mapsto -\textstyle\frac{1}{2}\, z(z^2-3y^2).\displaybreak[0]
\end{align*}
The standard $\SO(3)$ action on $\R^3$ induces an action on $\mathcal{H}^3(\R^3)$, hence on $\Im\O$.

Let $L_3$ be the orbit through $\eps_6$ of this $\SO(3)$ action on $\Im\O$ and let $L_4$ be the orbit through $\eps_2$.
By \cite[Theorem 4.3]{Mashimo} and observations in \cite[Examples 6.6 \& 6.15]{LotayLag}, $L_3\cong\SO(3)/\SSS_3$ and 
$L_4\cong\SO(3)/\AAA_4$ are Lagrangian.  Thus the cones $C_3$ and $C_4$ on $L_3$ and $L_4$ respectively are coassociative 
and $\SO(3)$-invariant.
\end{ex}

\begin{note}
The $\SO(3)$-orbit through $\epsilon_1$ is a constant curvature $\frac{1}{6}$ pseudoholomorphic curve in $\mathcal{S}^6$ often called the  \emph{Bor\r{u}vka sphere} in $\mathcal{S}^6$. 
\end{note}

All of the coassociative cones introduced so far have links which are fibered by oriented circles of constant radius over a surface.  
 Lagrangians in $\mathcal{S}^6$ with this property were classified in \cite{LotayLag}.  This classification 
 describes all of the known links of coassociative cones and it turns out that the circles have radius either $\frac{2}{3}$ or $1$.  
To describe these examples we make two definitions.
 
\begin{dfn}\label{tubedef}
Let $\mathbf{u}:\Sigma\rightarrow\mathcal{S}^6$ be a surface, let $\Pi$ be a 2-plane subbundle of $\mathbf{u}^*(T\mathcal{S}^6)$ and 
let $\mathcal{U}(\Pi)=\{\mathbf{v}\in\Pi\,:\,|\mathbf{v}|=1\}$.  For $\gamma\in(0,\frac{\pi}{2}]$ we define 
$\mathbf{x}_{\gamma}:\mathcal{U}(\Pi)\rightarrow\mathcal{S}^6$ by $\mathbf{x}_{\gamma}(\mathbf{v})=\cos\gamma\mathbf{u}+\sin\gamma
\mathbf{v}$.  We say that the image of $\mathbf{x}_\gamma$ is a \emph{tube of radius $\gamma$} (in $\Pi$) about $\Sigma$.
\end{dfn}

\begin{dfn}\label{pholodfn}
If $\mathbf{u}:\Sigma\rightarrow\mathcal{S}^6$ is a non-totally geodesic 
pseudoholomorphic curve in the sense of Definition \ref{Lagdef}, there is an orthogonal decomposition
$\mathbf{u}^*(T^{1,0}\mathcal{S}^6)=T^{1,0}\Sigma\oplus N_1\Sigma\oplus N_2\Sigma$,  
where $N_1\Sigma$ and $N_2\Sigma$ are holomorphic line bundles  
such that the 
second fundamental form of $\mathbf{u}$ takes values in $N_1\Sigma$.  We call $N_1\Sigma$ and $N_2\Sigma$ 
the \emph{first} and \emph{second normal bundle} respectively.

If $(\bff_1,\bff_2,\bff_3)$ is a moving orthonormal frame for $T^{1,0}\Sigma\oplus N_1\Sigma\oplus N_2\Sigma$ and 
$\theta_1$ is the $(1,0)$-form dual to $\bff_1$, then the structure equations for $\Sigma$ are:
\begin{align*}
\d\bfu&=-2i\bff_1\theta_1+2i\bar{\bff}_1\bar{\theta}_1;\displaybreak[0]
&\d\bff_1&=-i\bfu\bar{\theta}_1+\bff_1\kappa_{11}+\bff_2\kappa_{21};\displaybreak[0]\\
\d\bff_2&=-\bff_1\bar{\kappa}_{21}+\bff_2\kappa_{22}+\bff_3\kappa_{32}-\bar{\bff}_3\theta_1;\displaybreak[0]
&\d\bff_3&=-\bff_2\bar{\kappa}_{32}+\bff_3\kappa_{33}+\bar{\bff}_2\theta_1;\displaybreak[0]\\
\d\theta_1&=-\kappa_{11}\w\theta_1;\displaybreak[0]
&\d\kappa_{11}&=-\kappa_{21}\w\bar{\kappa}_{21}+2\theta_1\w\bar{\theta}_1;\displaybreak[0]\\
\d\kappa_{22}&=\kappa_{21}\w\bar{\kappa}_{21}-\kappa_{32}\w\bar{\kappa}_{32}-\theta_1\w\bar{\theta}_1;\displaybreak[0]&
\d\kappa_{33}&=\kappa_{32}\w\bar{\kappa}_{32}-\theta_1\w\bar{\theta}_1;\displaybreak[0]\\
\d\kappa_{21}&=(\kappa_{11}-\kappa_{22})\w\kappa_{21};\displaybreak[0]
&\d\kappa_{32}&=(\kappa_{22}-\kappa_{33})\w\kappa_{32},\displaybreak[0]
\end{align*}
for some imaginary-valued 1-forms $\kappa_{11},\kappa_{22},\kappa_{33}$ such that $\kappa_{11}+\kappa_{22}+\kappa_{33}=0$ and 
complex-valued 1-forms $\kappa_{21}$ and $\kappa_{32}$.  Moreover, by the work in \cite[$\S$4]{BryantOct}, there exist holomorphic 
functions $K$ and $T$ such that $\kappa_{21}=K\theta_1$ and $\kappa_{32}=T\theta_1$.  We identify $K$ with the second fundamental 
form of $\Sigma$, so $K\neq 0$, and we call $T$ the \emph{torsion} of $\Sigma$.  Of particular interest are the pseudoholomorphic 
curves with \emph{null-torsion}, i.e.~with $T\equiv 0$, since they may be viewed as certain algebraic curves in the 5-quadric in
 $\C\P^6$ and every compact Riemann surface can be realized as such a curve in $\mathcal{S}^6$ by the work in \cite[$\S$4]{BryantOct}.
\end{dfn}

\begin{remark}
The Bor\r{u}vka sphere $\Sigma$ in $\mathcal{S}^6$ has null-torsion.  
\end{remark}

The coassociative cones with links which are fibered by circles of radius $\frac{2}{3}$ are locally cones over tubes of radius 
$\sin^{-1}(\frac{2}{3})$ in $N_2\Sigma$ about null-torsion pseudoholomorphic curves $\Sigma\subseteq\mathcal{S}^6$.  This includes 
$C_1$ given in Example \ref{u2ex}, where $\Sigma$ is a totally geodesic $\mathcal{S}^2$ (which is the degenerate case of 
a null-torsion curve), and $C_4$ given in Example \ref{so3ex}, where $\Sigma$ is the Bor\r{u}vka sphere.

We shall be more concerned with the case where the link is fibered by oriented geodesic circles so we make the following 
 definition as in \cite{Lotay2r}.

\begin{dfn} A 4-dimensional submanifold $N$ of $\R^n$ is \emph{2-ruled} 
 if there is a surface $\Sigma$ and a smooth fibration $\pi:N\rightarrow\Sigma$
 whose fibres are affine 2-planes in $\R^n$.  In addition we say that $N$ is 
 \emph{r-framed} if there is a choice of oriented frame for each 2-plane $\pi^{-1}(\sigma)$ 
 which varies smoothly with $\sigma\in\Sigma$.
\end{dfn} 

 The 2-ruled coassociative 4-folds in $\R^7$ were 
 studied in \cite{Lotay2r}, and there was a further focus on the conical case in \cite{Foxcoass}.  
 By \cite[Theorems 1.1-1.2]{LotayLag} we can extend the work in \cite{Foxcoass} and describe all r-framed 2-ruled coassociative
 cones.  We refer the reader to \cite{LotayLag} for full details but give a brief description here.

\begin{ex}\label{2rex} ({\bf The 2-ruled family}) 
Complex 2-dimensional cones in $\C^3$ embedded in $\R^7$ are 2-ruled coassociative cones by 
 Corollary \ref{coasscor}.  

Recall Definition \ref{tubedef}.  The general r-framed 2-ruled coassociative cone $C$ has link $L$ such that, for all $p$ in an open 
dense subset $L^*$ of $L$, there exist an open set $U\ni p$, a non-totally geodesic pseudoholomorphic curve 
$\mathbf{u}:\Sigma\rightarrow\mathcal{S}^6$, and a holomorphic line subbundle $\Pi$ of $\mathbf{u}^*(T^{1,0}\mathcal{S}^6)$  
such that $U\cap L^*$ is a tube of radius $\frac{\pi}{2}$ in $\Pi$ about $\Sigma$.  There are restrictions on the choice of 
line subbundle $\Pi$ (see \cite[Example 7.4]{LotayLag}),  however, in particular, we may always choose $\Pi=N_2\Sigma$.
\end{ex}

We conclude this subsection with the following important examples of 2-ruled coassociative cones.

\begin{ex}\label{Killingex}
Let $C$ be a coassociative cone with link $L$ admitting a Killing vector field whose integral curves are geodesic circles in
 $\mathcal{S}^6$.
By \cite[Theorem 2]{Vrancken}, either $C$ is a complex 2-dimensional cone in $\C^3$ or $L$ is locally a tube of 
radius $\frac{\pi}{2}$ in $N_1\Sigma$ or $N_2\Sigma$ about  a \emph{null-torsion} pseudoholomorphic curve $\Sigma$ in $\mathcal{S}^6$.
\end{ex}

\noindent Cones as in Example \ref{Killingex} are given trivially by $C_2$ in Example \ref{cxso3ex} but also by $C_3$ 
in Example \ref{so3ex}, since $L_3$ is a tube of radius $\frac{\pi}{2}$ in $N_2\Sigma$ about the Bor\r{u}vka sphere $\Sigma$.

\subsection{Stability index and the curl operator}

For this subsection, use the notation of Definitions \ref{csdfn} and \ref{modulidfn} and Theorem \ref{maindefthm}.  In 
particular, $N$ is a CS coassociative 4-fold where 
the $s$ singularities have rate $\mu\in(1,2)\setminus\mathcal{D}$, where $\mathcal{D}$ is given in Definition \ref{Ddfn}, 
and the $s$-tuple of cones at the singularities $(C_1,\ldots,C_s)$ lies in $\mathfrak{C}=\prod_{i=1}^s\mathcal{C}_i$.  

   To get an effective notion of stability for the conical singularities of $N$ 
we need to improve the estimate \eq{obseq1} for the dimension of the obstruction space $\mathcal{O}(N,\mu,\mathfrak{C})$
  for the deformation problem for $N$.  
To understand this we need to compare the maps $(\d_++\d^*)_{\lambda}$ given by \eq{dd*eq} and
\begin{equation}
\begin{split}
\label{dd*eq2}
(\d+\d^*)_{\lambda}:L^p_{k+1,\,\lambda}(\Lambda^2T^*\Nhat\oplus\Lambda^4T^*\Nhat)&\rightarrow L^p_{k,\,\lambda-1}(\Lambda^3T^*\Nhat)\\
(\alpha,\beta)&\mapsto \d\alpha+\d^*\beta
\end{split}
\end{equation}
for $\lambda\in(-2,\mu]$.  
From Theorem \ref{maindefthm}, we see that $\mathcal{O}(N,\mu,\mathfrak{C})$ is a subspace of the cokernel of $(\d_++\d^*)_{\mu}$ which is transverse to the cokernel of $(\d+\d^*)_{\mu}$.  
 Moreover, from the work in \cite[$\S$8]{Lotaycs}, the sum over $\lambda\in(-2,\mu)\cap\mathcal{D}$ in \eq{obseq1} is the upper bound
 of the dimension of the space of forms which adds to the cokernel of $(\d_++\d^*)_{\lambda}$ as the rate $\lambda$ increases from
 $-2$ to $\mu$.   
Hence, we can improve the upper bound by considering the forms which add 
 to the cokernel of $(\d_++\d^*)_{\lambda}$ but \emph{not} the cokernel of $(\d+\d^*)_{\lambda}$.  
 These cokernels are isomorphic to the annihilators of the images of $(\d_++\d^*)_{\lambda}$ and $(\d+\d^*)_{\lambda}$ 
 under the dual pairing given in Definition \ref{dualdfn}.  Thus, by comparing these annihilators we may obtain our estimate.

We begin with the following.

\begin{prop}\label{annprop}
Let $\mathcal{A}_+(\lambda)$ and $\mathcal{A}(\lambda)$ denote the annihilators of the images of \eq{dd*eq} and 
\eq{dd*eq2} via the dual pairing given in Definition \ref{dualdfn}.  
For $\lambda\leq -1$, $\mathcal{A}_+(\lambda)=\mathcal{A}(\lambda)$. 
\end{prop}

\begin{proof}
As observed in the proof of 
\cite[Proposition 8.6]{Lotaycs}, if $\lambda\notin\mathcal{D}$, 
there exist finite-dimensional spaces $\mathcal{C}_+(\lambda)$ and $\mathcal{C}(\lambda)$ of smooth
 compactly supported 3-forms on $\Nhat$ such that 
$$ L^p_{k,\,\lambda-1}(\Lambda^3T^*\Nhat)=\Image(\d_++\d^*)_{\lambda}\oplus\mathcal{C}_+(\lambda)
=\overline{\Image(\d+\d^*)_{\lambda}}\oplus\mathcal{C}(\lambda)$$
 and the dual pairings between $\mathcal{C}_+(\lambda)$ and $\mathcal{A}_+(\lambda)$, and between 
 $\mathcal{C}(\lambda)$ and $\mathcal{A}(\lambda),$ are non-degenerate.  Let 
$\mathcal{C}^{\prime}(\lambda)$ be such that $\mathcal{C}_+(\lambda)=\mathcal{C}(\lambda)\oplus
\mathcal{C}^{\prime}(\lambda)$.  

Notice by Definition \ref{dualdfn} that if $q>1$ is such that $\frac{1}{p}+\frac{1}{q}=1$ and $l\in\N$ then, for any 
$\lambda\in\R$,
\begin{align*}
\mathcal{A}_+(\lambda)&=\big\{\gamma\in L^q_{l+1,-3-\lambda}(\Lambda^3T^*\Nhat)\,:\,\langle \d\alpha+\d^*\beta,\gamma \rangle=0\;\,\\
& \!\!\quad\qquad\qquad\qquad\qquad\qquad\qquad\text{for all $(\alpha,\beta)\in L^p_{k+1,\,\lambda}(\Lambda^2_+T^*\Nhat\oplus\Lambda^4T^*\Nhat)$}\big\}\\
&=\big\{\gamma\in L^q_{l+1,-3-\lambda}(\Lambda^3T^*\Nhat)\,:\langle \alpha,\d^*\gamma \rangle=0,\,\langle 
\beta,\d\gamma\rangle=0\;\,\\
& \!\!\quad\qquad\qquad\qquad\qquad\qquad\qquad\text{for all $(\alpha,\beta)\in L^p_{k+1,\,\lambda}(\Lambda^2_+T^*\Nhat\oplus\Lambda^4T^*\Nhat)$}\big\}\\
&=\big\{\gamma\in L^q_{l+1,-3-\lambda}(\Lambda^3T^*\Nhat)\,:\,\d\gamma=0,\,
\d^*\gamma\in L^q_{l,\,-4-\lambda}(\Lambda^2_-T^*\Nhat)\big\},
\end{align*}
where the integration by parts is justified by the choice of weight for the dual weighted Sobolev space.  Similarly,
\begin{align*}
\mathcal{A}(\lambda)=\big\{\gamma\in L^q_{l+1,-3-\lambda}(\Lambda^3T^*\Nhat)\,:\,
\d\gamma=\d^*\gamma=0\big\}.
\end{align*}

If $\lambda\leq-1$, then $-3-\lambda\geq\lambda-1$ so, by Theorem \ref{wembedthm}, $L^q_{l+1,-3-\lambda}\hookrightarrow
 L^p_{k,\,\lambda-1}$ for sufficiently large $l$.  Since $\mathcal{A}_+(\lambda)$ consists of smooth forms, as it is the kernel 
 of an elliptic operator, we can choose $l$ arbitrarily large and see that 
$\mathcal{A}_+(\lambda)\subseteq L^p_{k,\,\lambda-1}(\Lambda^3T^*\Nhat)$ and the same is clearly true for $\mathcal{A}(\lambda)$.  
  Moreover, if additionally $\lambda\notin\mathcal{D}$, since $\mathcal{A}_+(\lambda)$ and $\mathcal{A}(\lambda)$ are of equal
 dimension to $\mathcal{C}_+(\lambda)$ and $\mathcal{C}(\lambda)$, and the annihilators are by construction orthogonal to 
 the closures of the images of \eq{dd*eq} and \eq{dd*eq2}, we deduce that $\mathcal{A}_+(\lambda)=\mathcal{C}_+(\lambda)$ and 
 $\mathcal{A}(\lambda)=\mathcal{C}(\lambda)$. 
 
 Suppose that $\gamma\in\mathcal{C}^{\prime}(\lambda)$.  Then $\gamma$ is compactly supported and lies in $\mathcal{A}_+(\lambda)$.  
Therefore, since $\d^*\gamma$ is anti-self-dual,  
$$\|\d^*\gamma\|_{L^2}^2=\int_{\Nhat}-\d^*\gamma\w \d^*\gamma =\int_{\Nhat}-\d\!*\!\gamma\w
\d\!*\!\gamma=\int_{\Nhat}-\d(*\gamma\w \d\!*\!\gamma)=0,$$
where the integration by parts is valid since $\gamma$ is compactly supported.  Thus, $\gamma\in\mathcal{A}(\lambda)=\mathcal{C}(\lambda)$, so $\gamma=0$.  We deduce that $\mathcal{C}(\lambda)=\mathcal{C}_+(\lambda)
=\mathcal{A}(\lambda)=\mathcal{A}_+(\lambda)$ for $\lambda\leq -1$, $\lambda\notin\mathcal{D}$.  The annihilators are well-defined for 
$\lambda\in\mathcal{D}$ and hence, since the dimension of $\mathcal{A}_+(\lambda)$ is lower semi-continuous at $\lambda=-1$, as 
remarked in the proof of \cite[Proposition 8.4]{Lotaycs}, we see that $\mathcal{A}(\lambda)=\mathcal{A}_+(\lambda)$ for all 
$\lambda\leq -1$.
\end{proof}

\begin{prop}\label{improvdimprop}
Recall the notation of Definition \ref{Ddfn}.  Let
\begin{align*}
\check{D}(\lambda,i)=\{\gamma_i\in C^\infty(T^*L_i)\,:\,\d_i^*\gamma_i=0,\;\d_i\gamma_i=-(\lambda+2)*_i\!\gamma_i\}
\end{align*}
and let $\check{d}(\lambda)=\sum_{i=1}^s\dim \check{D}(\lambda,i)$.  Then 
\begin{equation}\label{obseq2}\dim\mathcal{O}(N,\mu,\mathfrak{C})\leq \!\!\!\!\!\sum_{\lambda\in(-1,\mu)\cap\mathcal{D}}
\!\!\!\! \check{d}(\lambda)-\sum_{i=1}^s\dim\mathcal{C}_i.
\end{equation}
\end{prop}

\begin{proof}
By Proposition \ref{annprop}, there are no forms which add to the cokernel of \eq{dd*eq} but not 
 the cokernel of \eq{dd*eq2} for $\lambda\leq -1$.  Therefore, we can first improve the estimate 
 \eq{obseq1} by restricting to the range of rates $(-1,\mu)\cap\mathcal{D}$.

The map $(\d+\d^*)_\lambda$, given in \eq{dd*eq2}, is not elliptic, but forms part of the map 
$$\d+\d^*:L_{k+1,\,\lambda}^p(\Lambda^{\text{even}}T^*\hat{N})\rightarrow L_{k,\,\lambda-1}^p(\Lambda^{\text{odd}}T^*\hat{N}),$$
which \emph{is} elliptic.  Therefore, we can apply the theory of \cite{LockhartMcOwen} and deduce that 
there is a countable discrete set of rates $\mathcal{E}\supseteq\mathcal{D}$ for which $\d+\d^*$ is Fredholm and, moreover, we 
can calculate which forms on $L_i$ correspond to forms on $\hat{N}$ which subtract from the kernel or add to the cokernel as  
$\lambda$ increases.   

The cokernel of $\d+\d^*$ is isomorphic to the annihilator of the image of $\d+\d^*$, 
using the dual pairing given in Definition \ref{dualdfn}.   The work in \cite{LockhartMcOwen} identifies the forms which add to 
this annihilator as the rate crosses elements $\lambda\in\mathcal{E}$.  To calculate the changes in the annihilator, we need to 
consider homogeneous forms on the cone $C_i$ of the following type:
$$(\gamma^1,\gamma^3)=(r^{-\lambda-2}\alpha^1+r^{-\lambda-3}\alpha^0\w\d r,r^{-\lambda}\alpha^3+r^{-\lambda-1}\alpha^2\w\d r)$$ for
$(\alpha^0,\alpha^1,\alpha^2,\alpha^3)\in C^{\infty}(\oplus_{m=0}^3\Lambda^mT^*L_i)$, which 
satisfy $\d^*\gamma^1=\d\gamma^1+\d^*\gamma^3=\d\gamma^3=0$, where $\d^*$ is calculated using the cone metric.  
We therefore have the following conditions which define $\mathcal{E}$ 
and the changes to the kernel or cokernel of $\d+\d^*$ acting on even forms, using the notation of Definition \ref{Ddfn}:
\begin{gather}
\d_i^*\alpha^1=-\lambda\alpha^0,\quad \d_i\alpha^0+\d_i^*\alpha^2=-(\lambda+2)\alpha^1,\label{dd*eq3a}\\
\d_i\alpha^2=-\lambda\alpha^3,\quad\d_i\alpha^1+\d_i^*\alpha^3=-(\lambda+2)\alpha^2.\label{dd*eq3b}
\end{gather}

Thus, by \eq{dd*eq3a}-\eq{dd*eq3b}, forms $(\alpha,\beta)\in D(\lambda,i)$ giving rise to cokernel forms for $(\d_++\d^*)_{\lambda}$, 
$\lambda\in\mathcal{D}$, which can lie in $\mathcal{O}(N,\mu,\mathfrak{C})$ must be transverse to forms
 $(\alpha^{\prime},\beta^{\prime})\in C^{\infty}(\Lambda^2T^*L_i\oplus\Lambda^3T^*L_i)$ satisfying:
\begin{equation}\label{dd*eq4}
\d_i\alpha^{\prime}=-\lambda\beta^{\prime},\quad \d_i^*\beta^{\prime}=-(\lambda+2)\alpha^{\prime}, \quad \d_i^*\alpha^{\prime}=0.
\end{equation}
For $\lambda\neq -2,0$, solutions to \eq{dd*eq4} are equivalent to giving exact $\beta^{\prime}\in C^{\infty}(\Lambda^3T^*L_i)$ with
 $\Delta_i\beta^{\prime}=\lambda(\lambda+2)\beta^{\prime}$, where $\Delta_i$ is the Laplacian on $L_i$, and setting
 $\alpha^{\prime}=-(\lambda+2)^{-1}\d_i^*\beta^{\prime}$.  For $\lambda=0$, $\beta^{\prime}$ is locally constant and
 $\alpha^{\prime}=0$.
 Since $\beta$ is also an exact eigenform of $\Delta_i$ 
with eigenvalue $\lambda(\lambda+2)$ for $\lambda\neq 0$ and is constant if $\lambda=0$, we quickly deduce that $\beta=0$.   
 Notice that for the case $\lambda=-2$, $\beta$ must be zero since it is harmonic and exact.
 Observe that $\alpha$ is automatically orthogonal to coexact forms since $\d_i\alpha=0$ and $-(\lambda+2)\alpha=\d_i\!*_i\!\alpha$.  
Setting $\gamma_i=*_i\alpha$, we see that we can replace the quantities $d(\lambda)$ by $\check{d}(\lambda)$ in our estimate \eq{obseq2} for the 
dimension of $\mathcal{O}(N,\mu,\mathfrak{C})$ as claimed.
\end{proof}

Theorem \ref{maindefthm} and Proposition \ref{improvdimprop} invite us to define the stability index 
of a coassociative cone in a similar manner to \cite[Definition 3.6]{JoyceSLsing2}.

\begin{dfn}\label{stabdfn}
Let $C$ be a coassociative cone in $\R^7$ with compact link $L\subseteq\mathcal{S}^6$ such that $C\setminus\{0\}$ 
is non-singular. 
Let $\mathcal{C}$ be a smooth, connected family of coassociative cones in $\R^7$  
 which contains $C$ and is closed under the natural action of $\GG_2\ltimes\R^7$.  Finally, for $\lambda\in\R$, let
\begin{equation}\label{checkDeq}
\check{D}(\lambda)=\{\gamma\in C^\infty(T^*L)\,:\,\d^*\gamma=0,\;\d\gamma=-(\lambda+2)*\gamma\}.
\end{equation}

We define the $\mathcal{C}$-\emph{stability index} of $C$ by
\begin{equation}\label{sindeq}
\ind_{\mathcal{C}}(C)=\sum_{\lambda\in(-1,1]}\!\!\!\! \dim\check{D}(\lambda)-\dim\mathcal{C}.
\end{equation}
If the family $\mathcal{C}$ consists solely of the $\GG_2\ltimes\R^7$
 transformations of $C$, we simply write $\ind_{\mathcal{C}}(C)=\ind(C)$ and call $\ind(C)$ the \emph{stability 
index} of $C$.
 
Notice that the sum in \eq{sindeq} is well-defined because the set of $\lambda$ for which
 $\dim\check{D}(\lambda)\neq 0$ is countable and discrete by the observations in Definition \ref{Ddfn}.
 Moreover, $\dim\check{D}(1)$ is the dimension of the space of Lagrangian Jacobi fields 
  on $L$ by Definition \ref{Jacdfn}, so the space of all infinitesimal deformations of $C$ as a coassociative cone, and 
$\check{D}(0)$ corresponds to the $O(1)$ closed self-dual 2-forms on $C$, so $\dim\check{D}(0)$ is at least as large as the 
space of translations of $C$.  Thus, $\ind_{\mathcal{C}}(C)\geq 0$.

We say that the cone $C$ is $\mathcal{C}$-\emph{stable} if $\ind_{\mathcal{C}}(C)=0$.
 If $\ind(C)=0$ we say that the cone $C$ is \emph{stable}.  We also say that $C$ is \emph{rigid} if 
 $\dim\check{D}(1)=14-\dim\GG$, where $\GG$ is the Lie subgroup of $\GG_2$ preserving $C$.
\end{dfn}

\begin{notes}
\begin{itemize}\item[]\item[(a)] 
It is clear that stability of $C$ implies rigidity. It is also clear that if $C$ is rigid then it is Jacobi integrable, since 
then every Lagrangian deformation of $L$ in $\mathcal{S}^6$ comes from $\GG_2$ transformations.  However, we shall see that one may have cones 
which are Jacobi integrable but neither stable nor rigid.
\item[(b)] By Theorem \ref{maindefthm} and Proposition \ref{improvdimprop}, if all the cones $C_i$ at the singularities of a CS
 coassociative 4-fold $N$ are $\mathcal{C}_i$-stable and $\mu\in(1,2)$ is such that $(1,\mu]\cap\mathcal{D}=\emptyset$, then the obstruction space $\mathcal{O}(N,\mu,\mathfrak{C})=\{0\}$ and so 
the moduli space of deformations $\mathcal{M}(N,\mu,\mathfrak{C})$ is a smooth manifold near $N$.  Moreover, $N$ is `stable' under
deformations of the ambient closed $\GG_2$ structure by Theorem \ref{G2deformthm}.
\end{itemize}
\end{notes}

We see from \eq{checkDeq} that to determine the stability index for a coassociative cone we need to study the equation  
$$*\d\gamma=-(\lambda+2)\gamma$$
for 1-forms $\gamma$ on a compact Riemannian 3-manifold $L$ for $\lambda\in(-1,1]$.  
The operator $*\d:C^{\infty}(T^*L)\rightarrow C^{\infty}(T^*L)$ is 
a natural self-adjoint operator on $L$ which we call the \emph{curl operator}.  Thus, our problem is to calculate 
the negative eigenvalues of the curl operator and their multiplicities.  For convenience we make the following 
definition.

\begin{dfn}\label{cLdfn}
Let $(L,g_L)$ be a compact 
Riemannian 3-manifold and let $c_L=-*\d:C^{\infty}(T^*L)\rightarrow C^{\infty}(T^*L)$.  Denote by $\sigma_L(\nu)$ the
 multiplicity of a non-zero eigenvalue $\nu$ of $c_L$.  Note that $\dim\check{D}(\lambda)=\sigma_L(\lambda+2)$ for $\lambda\neq -2$.
\end{dfn}

Finding the positive spectrum of $c_L$ and the multiplicities is an extremely complicated problem and 
in general there is no hope to solve it.  However, since we need only consider eigenvalues in $(1,3)$, 
the problem is tractable in special cases.  For possible further applications to coassociative 
geometry it is of greatest practical use to study eigenvalues in the range $(0,4)$.

\section{Homogeneous cones}
\label{homconessec}

In this section  
we explicitly determine the stability index for coassociative cones whose links are orbits of closed 
3-dimensional subgroups of $\GG_2$.  We achieve this by calculating the small eigenvalues of the curl operator on 
Berger 3-spheres and their quotients using elementary methods.
 
\subsection{Berger 3-spheres}

\begin{dfn}
\label{Bergerdfn}
Let $\H$ denote the quaternions with standard basis $\{\bfo,\bfi,\bfj,\bfk\}$.  
Identify $\mathcal{S}^3\cong\Sp(1)$ and let $\bfx:\Sp(1)\rightarrow\H$ denote the inclusion map of $\Sp(1)$ as unit quaternions.
  Then $\d\bfx=\bfx\omega$ for a 1-form $\omega$ taking values in the Lie algebra of $\Sp(1)$, which here is represented by 
$\Im\H$.  Therefore 
  we can write 
 $$\bfx=x_0\bfo+x_1\bfi+x_2\bfj+x_3\bfk\quad\text{and}\quad\omega=\omega_1\bfi+\omega_2\bfj+\omega_3\bfk.$$
 Since the Maurer--Cartan form $\omega$ satisfies the structure equation $\d\omega+\omega\w
 \omega=0$ and $\d\bfx=\bfx\omega$, we have that
\begin{align}
\d x_0&=-x_1\omega_1-x_2\omega_2-x_3\omega_3, & \d x_1&=x_0\omega_1-x_3\omega_2+x_2\omega_3,\label{structeq1a}\\ 
\d x_2&=x_3\omega_1+x_0\omega_2-x_1\omega_3, & \d x_3&=-x_2\omega_1+x_1\omega_2+x_0\omega_3\label{structeq1b}
\end{align}
and
\begin{gather}
\d\omega_1=-2\omega_2\w\omega_3,\qquad\d\omega_2=-2\omega_3\w\omega_1,\qquad\d\omega_3=-2\omega_1\w\omega_2.\label{structeq2}
\end{gather}
 We define a 1-parameter family of metrics on $\mathcal{S}^3$ by $g_{\tau^2}=\tau^2\omega_1^2+\omega_2^2+\omega_3^2$ 
 for $\tau>0$.  The Riemannian manifolds $(\mathcal{S}^3,g_{\tau^2})$ are the \emph{Berger 3-spheres}.
\end{dfn}

It is immediately clear that finding the eigenvalues of $*\d$ on a Berger 3-sphere will involve the Laplacian acting on 
$\mathcal{S}^3$, thus homogeneous harmonic polynomials on $\R^4$.  We are thus lead to make the following definitions.

\begin{dfn}\label{weightdfn}
Use the notation of Definition \ref{Bergerdfn}.  For $m=1,2,3$, we define operators $\partial_m$ on 
$f\in C^{\infty}(\mathcal{S}^3)$ by the expression: $\d f=\partial_1 f\omega_1+\partial_2f\omega_2+\partial_3f\omega_3$. 
For a unit imaginary quaternion $\mathbf{q}=q_1\mathbf{i}+q_2\mathbf{j}+q_3\mathbf{k}$ we define 
$\partial_{\mathbf{q}}=q_1\partial_1+q_2\partial_2+q_3\partial_3$.  We say that $f\in C^{\infty}(\mathcal{S}^3)$ has 
\emph{$\mathbf{q}$-weight} $w\geq 0$ if $\partial_{\mathbf{q}}^2f=-w^2f$.

Observe that the Laplacian $\Delta_{\tau^2}$ on $(\mathcal{S}^3,g_{\tau^2})$ acting on functions is given by 
$\Delta_{\tau^2}=-\frac{1}{\tau^2}\partial_1^2-\partial_2^2-\partial_3^2$.
\end{dfn}

\noindent Of course, $\partial_m$ is just Lie derivative along the vector field dual to $\omega_m$.

\begin{dfn}\label{polydfn}
Use the notation of Definition \ref{Bergerdfn}. 
For $k\in\N$, let $\mathcal{P}_k$ be the space of homogeneous polynomials in $x_0,x_1,x_2,x_3$ of degree $k$ on $\mathcal{S}^3$.
Let $\mathcal{Q}_k$ be the subspace of $\mathcal{P}_k$ consisting of polynomials which are eigenfunctions of the 
standard Laplacian $\Delta_{1}$; that is, restrictions of homogeneous harmonic polynomials in 4 real variables to $\mathcal{S}^3$.  
We let $\mathcal{A}_k=\{p_1\omega_1+p_2\omega_2+p_3\omega_3\,:\,p_1,p_2,p_3\in\mathcal{P}_k\}$ and 
$\mathcal{B}_k=\{p_1\omega_1+p_2\omega_2+p_3\omega_3\,:\,p_1,p_2,p_3\in\mathcal{Q}_k\}$.
\end{dfn}

It is often clearer to work with the representation of $\Sp(1)=\SU(2)$ on $\C^2$ rather than $\H$, so we make the following 
useful definition.

\begin{dfn}\label{cxpolydfn}
Recall the notation of Definition \ref{Bergerdfn}.  We define complex coordinates on $\H\cong\C^2$ by $z_1=x_0+ix_1$ and 
$z_2=x_2+ix_3$.  For $\mathbf{q}=q_0\bfo+q_1\bfi+q_2\bfj+q_3\bfk\in\Sp(1)$, its action on $\C^2$ is given by
$$\mathbf{q}:\left(\begin{array}{c}z_1 \\ z_2\end{array}\right)\mapsto\left(\begin{array}{rr} q_0+iq_1 & q_2+iq_3\\ -q_2+iq_3 & q_0-iq_1
 \end{array}\right)\left(\begin{array}{c} z_1\\ z_2\end{array}\right).$$
For $k\in\N$, let $\mathcal{P}_k^{\C}$ denote the space of homogeneous polynomials in $z_1,z_2,\bar{z}_1,\bar{z}_2$ of 
degree $k$ restricted to $\mathcal{S}^3$ and let 
$\mathcal{Q}_k^{\C}$ be the subspace of $\mathcal{P}_k^{\C}$ consisting of polynomials which are harmonic on $\C^2$.  We also let 
$\mathcal{R}_m^{\C}$, for $m\in\Z$, be the space of homogeneous polynomials $p$ in $z_1,z_2,\bar{z}_1,\bar{z}_2$, restricted 
to $\mathcal{S}^3$, such that under the action of $\cos\theta\bfo+\sin\theta\bfi$, $p$ maps to $e^{im\theta}p$.  Notice that 
$p\in\mathcal{R}_m^{\C}$ if and only if $\bar{p}\in\mathcal{R}_{-m}^{\C}$.
\end{dfn}

We now recall the following well-known facts concerning eigenfunctions of the Laplacian on $\mathcal{S}^3$.

\begin{thm}\label{harmthm}
Use the notation of Definitions \ref{Bergerdfn}-\ref{polydfn}.  
\begin{itemize}
\item[\emph{(a)}] For each $k\in\N$, there is a direct sum decomposition 
$\mathcal{Q}_k=\oplus_{l=0}^{[k/2]}\mathcal{Q}_{k,k-2l}$ such that the elements of $\mathcal{Q}_{k,k-2l}$ have $\mathbf{i}$-weight 
$k-2l$.  Moreover, 
$$\dim\mathcal{Q}_{k,k-2l}=2k+2\;\,\text{if\/ $l<[\textstyle\frac{k}{2}]$,}\quad
\dim\mathcal{Q}_{2l,0}=2l+1$$
and $\dim\mathcal{Q}_k=(k+1)^2$. 
\item[\emph{(b)}] The eigenvalues of $\Delta_{\tau^2}$ are of the form $k(k+2)+(\frac{k-2l}{\tau})^2(1-\tau^2)$ for $k\in\N$ and $l\leq[\frac{k}{2}]$, and the corresponding eigenspace is $\mathcal{Q}_{k,k-2l}$.
\end{itemize}
\end{thm}

\begin{proof}
The results of \cite[Lemmas 3.1 \& 4.1]{Tanno} give the decomposition in (a) and all of (b).  We can determine the dimensions of 
the $\mathcal{Q}_{k,k-2l}$ by explicitly identifying the functions as the real and imaginary parts of elements in 
 $\mathcal{Q}_k^{\C}\cap\mathcal{R}_{k-2l}^{\C}$, in the notation of Definition \ref{cxpolydfn}.  Clearly $\mathcal{Q}_{2l,0}$ is 
 the space of lifts of eigenfunctions of the Laplacian on $\C\P^1$ with eigenvalue $4l(l+1)$, which has dimension $2l+1$.  For 
 the remaining spaces, it is straightforward to see that each $\mathcal{Q}_{k,k-2l}$ has the same dimension for fixed $k$, and 
 since $\dim\mathcal{Q}_k=(k+1)^2$, as it is the multiplicity of the eigenvalue $k(k+2)$ for $\Delta_1$ on $\mathcal{S}^3$, 
 it is an elementary calculation to find that $\dim\mathcal{Q}_{k,k-2l}=2k+2$ if $l\neq[\frac{k}{2}]$.
\end{proof}

\noindent The proofs of the quoted results from \cite{Tanno} rest 
on the fact that $\partial_1$ commutes with $\Delta_{\tau^2}$ for any 
$\tau>0$.  This is certainly not true of $\partial_2$ and $\partial_3$ if $\tau\neq 1$.  

\begin{prop}\label{Bergerprop}  
Use the notation of Definitions \ref{cLdfn} and \ref{Bergerdfn} and Theorem \ref{harmthm}.  
The positive eigenvalues of $c_L$ 
on $(L,g_L)=(\mathcal{S}^3,g_{\tau^2})$ are 
\begin{equation*}
\nu_{k,k-2l}=\tau+
\sqrt{\tau^2+k(k+2)+(\textstyle\frac{k-2l}{\tau})^2(1-\tau^2)} 
\quad\text{and}\quad\nu_k=\frac{k+2}{\tau}
\end{equation*}
for $k,l\in\N$ with $l\leq[\frac{k}{2}]$.  Moreover, 
$$\sigma_L(\nu_{k,k-2l})=\dim\mathcal{Q}_{k,k-2l}\quad\text{and}\quad\sigma_L(\nu_k)=2k+2.$$
\end{prop}

\begin{note}
For the multiplicity count here we regard the $\nu_{k,k-2l}$ and $\nu_k$ as distinct.  If they agree then we add the multiplicities.
\end{note}

\begin{proof}
From \eq{structeq1a}-\eq{structeq2} 
we see that $c_L=-\!*\d$ sends $\mathcal{A}_k$, given in Definition \ref{polydfn}, to itself.  Since $\cup_{k\in\N}\mathcal{A}_k$ is
 dense in $C^{\infty}(T^*L)$ we 
need only consider $c_L\alpha=\nu\alpha$ for $\alpha\in\mathcal{A}_k$ to determine the eigenvalues $\nu$ of $c_L$.

The equation $c_L\alpha=\nu\alpha$ for 
$\alpha=p_1\omega_1+p_2\omega_2+p_3\omega_3$ and $\nu>0$ is equivalent to the following system, using the notation of Definition 
\ref{weightdfn}:
\begin{align}
(2-\textstyle\frac{\nu}{\tau})p_1&=\partial_2p_3-\partial_3p_2,\label{curleq1}\\
(2-\nu\tau)p_2&=\partial_3p_1-\partial_1p_3,\label{curleq2}\\ (2-\nu\tau)p_3&=\partial_1p_2-\partial_2p_1.\label{curleq3}
\end{align}
Moreover, since $*\d\alpha=-\nu\alpha$ for $\nu\neq 0$, we have that $\d^*\alpha=0$, which is equivalent to the condition:
\begin{equation}\label{diveq}
\textstyle\frac{1}{\tau}\partial_1p_1+\tau\partial_2p_2+\tau\partial_3p_3=0.
\end{equation}
From \eq{structeq1a}-\eq{structeq1b} we see that, if $\epsilon_{abc}$ is the standard permutation symbol,
\begin{equation}\label{partialeq}
[\partial_a,\partial_b]=2\epsilon_{abc}\partial_c.
\end{equation}

Using \eq{curleq1}-\eq{partialeq} we calculate:
\begin{align}
\Delta_{\tau^2}p_1&=-(\textstyle\frac{1}{\tau^2}\partial_1^2+\partial_2^2+\partial_3^2)p_1\nonumber\\
&=
\big([\partial_1,\partial_2]-(2-\nu\tau)\partial_3\big)p_2+\big([\partial_1,\partial_3]+(2-\nu\tau)\partial_2\big)p_3\nonumber\\
&=-\nu\tau(\partial_2p_3-\partial_3p_2)\nonumber\\
&=\nu(\nu-2\tau)p_1.\label{nonzerop1eq1}
\end{align}
Thus $p_1$ is a $\nu(\nu-2\tau)$-eigenfunction of $\Delta_{\tau^2}$ if $p_1\neq0$ and hence  
\begin{equation}\label{nonzerop1eq2}
\nu(\nu-2\tau)=k(k+2)+(\textstyle\frac{k-2l}{\tau})^2(1-\tau^2),
\end{equation}
from which the formula for $\nu_{k,k-2l}$ follows.  

If $p_1=0$ then \eq{curleq1}-\eq{partialeq} imply that
\begin{equation}\label{zerop1eq}
-\partial_1^2p_j=(\nu\tau-2)^2p_j\quad\text{and}\quad -(\partial_2^2+\partial_3^2)p_j=2(\nu\tau-2)p_j
\end{equation}
for $j=2,3$.  Thus, $\Delta_1 p_j= \nu\tau(\nu\tau-2)p_j$, so $\nu\tau(\nu\tau-2)=k(k+2)$ from which the formula for $\nu_k$ follows.

To determine the multiplicity of $\nu_{k,k-2l}$ we make the following observations.  First, using \eq{curleq2}-\eq{curleq3}
we see that 
\begin{align}
\big(\partial_1^2+(\nu\tau-2)^2\big)p_2&=
\big(\partial_1\partial_2-(\nu\tau-2)\partial_3\big)p_1,\label{p2eq1}\\
\big(\partial_1^2+(\nu\tau-2)^2\big)p_3&=\big(\partial_1\partial_3+(\nu\tau-2)\partial_2\big)p_1.\label{p3eq1}
\end{align}
Second, from \eq{curleq1}-\eq{diveq}, we calculate
\begin{align}
\big(\partial_2^2+\partial_3^2+2(\nu\tau-2)\big)p_2&=-\big(\textstyle\frac{1}{\tau^2}\partial_2\partial_1+
(4-\textstyle\frac{\nu}{\tau})\partial_3)p_1,\label{p2eq2}\\
\big(\partial_2^2+\partial_3^2+2(\nu\tau-2)\big)p_3&=-\big(\textstyle\frac{1}{\tau^2}\partial_3\partial_1-(4-\textstyle\frac{\nu}{\tau})\partial_2)p_1.\label{p3eq2}
\end{align}
Combining \eq{p2eq1}-\eq{p3eq2} we see that $p_1$ determines $p_2$ and $p_3$ unless they satisfy \eq{zerop1eq}, which 
happens if and only if $\nu_{k,k-2l}=\nu_k$.  Thus, the multiplicity of $\sigma_L(\nu_{k,k-2l})$ is determined by the number 
of choices for $p_1$, so $\sigma_L(\nu_{k,k-2l})=\dim\mathcal{Q}_{k,k-2l}$.

For the multiplicity of $\nu_k$, we may take $p_1=0$ and see that \eq{curleq3} determines $p_3$ given $p_2$ unless $k=0$.  
  If $k\neq 0$, $\sigma_L(\nu_k)$ is the number of choices for $p_2$.  By \eq{zerop1eq}, we see that $p_2\in\mathcal{Q}_{k,k}$, so  $\sigma_L(\nu_k)=\dim\mathcal{Q}_{k,k}=2k+2$ by Theorem \ref{harmthm}.  
  For $k=0$, $p_2$ and $p_3$ are arbitrary constants so $\sigma_L(\nu_0)=2$.
\end{proof} 

\begin{notes}
We have some important elementary observations from the proof of Proposition \ref{Bergerprop} which will be important later.
\begin{itemize}\item[(a)] There is a basis for the $\nu_{k,k-2l}$-eigenforms on $(\mathcal{S}^3,g_{\tau^2})$ consisting of forms 
$p_1\omega_1+p_2\omega_2+p_3\omega_3$ with $p_1\in\mathcal{Q}_{k,k-2l}$ and $p_2,p_3\in\mathcal{Q}_k$ determined by 
$p_1$.
\item[(b)] There is a basis for the $\nu_k$-eigenforms on $(\mathcal{S}^3,g_{\tau^2})$ consisting of forms $p_2\omega_2+p_3\omega_3$ such that
 $p_2,p_3\in\mathcal{Q}_{k,k}$ and $p_3$ is determined by $p_2$ if $k>0$.
\end{itemize}
\end{notes}

We have an immediate corollary to Proposition \ref{Bergerprop} which builds on an already well-known result.

\begin{cor}
\label{roundcor}
Use the notation of Definitions \ref{cLdfn} and \ref{Bergerdfn}.  The positive eigenvalues of $c_{L_0}$ on 
$(L_0,g_{L_0})=(\mathcal{S}^3,g_1)$ are $k+2$ for $k\in\N$ and $\sigma_{L_0}(k+2)=(k+2)^2-1$.
  Hence, a coassociative 4-plane in $\R^7$ is stable.  
\end{cor}

\begin{proof}
Clearly if $\tau=1$ then $\nu_{k,k-2l}=1+\sqrt{(k+1)^2}=k+2=\nu_k$.  If $k=2n$, then the multiplicity of the eigenvalue 
$2n+2$ is 
$$\sum_{l=0}^{n-1}\sigma_{L_0}(\nu_{2n,l})+\sigma_{L_0}(\nu_{2n,n})+\sigma_{L_0}(\nu_{2n})=n(4n+2)+2n+1+4n+2=(2n+2)^2-1.$$
Similarly, if $k=2n+1$ the multiplicity of the eigenvalue $2n+3$ is
$$\sum_{l=0}^{n}\sigma_{L_0}(\nu_{2n+1,l})+\sigma_{L_0}(\nu_{2n+1})=(n+1)(4n+4)+4n+4=(2n+3)^2-1.$$
The multiplicities follow.
 
For a coassociative 4-plane $C_0$,
$$\sum_{\lambda\in(-1,1]}\!\!\!\! \dim\check{D}(\lambda)=\dim\check{D}(0)+\dim\check{D}(1)=\sigma_{L_0}(2)+\sigma_{L_0}(3)=3+8=11.$$
The number of non-trivial translations of $C_0$ is 3, and the stabilizer of 
$C_0$ in $\GG_2$ is isomorphic to $\SO(4)$ as observed in Example \ref{planeex}.  Thus, the dimension of the space of 
$\GG_2\ltimes\R^7$ transformations of $C_0$ is $3+14-6=11$.  Therefore 
$\ind(C_0)=0$.  
\end{proof}

By \cite[Example 5.1 \& Theorem 5.1]{Dillen}, the Lagrangian $L_1$ given in Example \ref{u2ex} is isometric to 
$\mathcal{S}^3$ with the metric $\frac{8}{3}g_{\frac{1}{6}}$.  We may thus calculate the stability index of $C_1$ as follows. 

\begin{cor}
\label{squashedcor}
Use the notation of Definitions \ref{cLdfn} and \ref{Bergerdfn} and Example \ref{u2ex}.  The eigenvalues of $c_{L_1}$ on  
$(L_1,g_{L_1})=(\mathcal{S}^3,\frac{8}{3}g_{\frac{1}{6}})$ which lie in $(0,4)$ are given by
$\{\textstyle\frac{1}{2},2,3,\textstyle\frac{1+\sqrt{145}}{4},\textstyle\frac{7}{2}\}$.   
Moreover, 
$$\sigma_{L_1}(\textstyle\frac{1}{2})=1,\quad\sigma_{L_1}(2)=7,\quad\sigma_{L_1}(3)=10,\quad\sigma_{L_1}(\textstyle\frac{1+\sqrt{145}}{4})=5,
\quad \sigma_{L_1}(\textstyle\frac{7}{2})=6.$$
Hence, the coassociative cone $C_1$ is stable. 
\end{cor}

\begin{proof}
Let $\lambda_{k,k-2l}=\sqrt{\frac{3}{8}}\nu_{k,k-2l}$ and let $\lambda_k=\sqrt{\frac{3}{8}}\nu_{k}$, 
using Proposition \ref{Bergerprop} with 
$\tau=\frac{1}{\sqrt{6}}$.  Then $\lambda_{k,k-2l}$ and $\lambda_k$ are the eigenvalues of $c_{L_1}$.   

We notice that $\lambda_{5,1}=\frac{1}{4}(1+\sqrt{241})>4$, so we need only consider 
$\lambda_{k,k-2l}$ for $k\leq 4$ for eigenvalues in $(0,4)$.  We may calculate:
\begin{gather*}
\lambda_{0,0}=\textstyle\frac{1}{2},\quad\lambda_{1,1}=2,\quad\lambda_{2,2}=\textstyle\frac{7}{2},\quad\lambda_{2,0}=2,\quad
\lambda_{3,3}=5,\quad\lambda_{3,1}=3,\\ \quad\lambda_{4,4}=\textstyle\frac{13}{2},\quad\lambda_{4,2}=\textstyle\frac{1}{4}(1+\sqrt{265})>4,\quad
\lambda_{4,0}=\textstyle\frac{1}{4}(1+\sqrt{145})<4.
\end{gather*}
Since $\lambda_k=\frac{3k}{2}+3$, only $\lambda_0=3$ is relevant here.  By 
Proposition \ref{Bergerprop},
\begin{gather*}
\sigma_{L_1}(\textstyle\frac{1}{2})=\dim\mathcal{Q}_{0,0},\quad\sigma_{L_1}(2)=\dim\mathcal{Q}_{1,1}+\dim\mathcal{Q}_{2,0},
\quad\sigma_{L_1}(3)=\dim\mathcal{Q}_{3,1}+2,\\
\sigma_{L_1}(\textstyle\frac{1+\sqrt{145}}{4})=\dim\mathcal{Q}_{4,0},\quad\sigma_{L_1}(\textstyle\frac{7}{2})=\dim\mathcal{Q}_{2,2}.
\end{gather*}
The multiplicities now follow from Theorem \ref{harmthm}.

Since $C_1$ is non-planar and the stabilizer of $C_1$ under $\GG_2$ transformations is $\U(2)$, the dimension of the family 
 of $\GG_2\ltimes\R^7$ transformations of $C_1$ is $7+14-4=17$.  Further,  
$$\sum_{\lambda\in(-1,1]}\!\!\!\! \dim\check{D}(\lambda)=\dim\check{D}(0)+\dim\check{D}(1)=\sigma_{L_1}(2)+\sigma_{L_1}(3)=7+10=17,$$
so $\ind(C_1)=0$.
\end{proof}

\begin{remark}
We may observe, as in \cite[Theorem 5]{Urakawa}, that 
$\tau^2=\frac{1}{6}$ is 
the critical value at which the multiplicity 
of the first eigenvalue of $\Delta_{\tau^2}$ on $\mathcal{S}^3$ ``jumps''.
\end{remark}

\subsection{Quotients of the 3-sphere}

We now wish to consider quotients of $\mathcal{S}^3$ by finite groups.  The possible finite groups are the cyclic groups $\Z_n$, 
the binary dihedral groups $\DD_{n}^*$, the binary tetrahedral group $\AAA_4^*$, the binary octahedral group $\SSS_4^*$ and the 
binary icosahedral group $\AAA_5^*$.  We know, from $\S$\ref{orbitsec}, that the only groups we must consider are $\Z_2$, 
$\DD_{3}^*$ and $\AAA_4^*$. 
 We describe the actions of the finite subgroups we need explicitly.
 
\begin{dfn}\label{finitegpdfn}
Let unit quaternions act on $\H$ by left-multiplication. 
\begin{itemize}
\item[(a)]  The cyclic group of order $n\geq 2$ acts as
$\Z_n=\{\cos(\frac{2k\pi}{n})\bfo+\sin(\frac{2k\pi}{n})\bfi\,:\,k=0,1,\ldots,n-1\}$.
\item[(b)] The binary dihedral group of order $4n$, for $n\geq 1$, acts as 
$\DD_{n}^*=\Z_{2n}\cup\bfj\Z_{2n}$. 
\item[(c)] The binary tetrahedral group acts as  
$\AAA_4^*  
=\DD_{2}^*\cup\{\frac{1}{2}(\pm\bfo\pm\bfi\pm\bfj\pm\bfk)\}$, where any combination of signs is permissible. 
 \end{itemize}
\end{dfn}
 
 To understand the spectrum of $c_L$ on quotients of $\mathcal{S}^3$ we need the following.
 
\begin{prop}\label{commprop}
Use the notation of Definitions \ref{cLdfn} and \ref{Bergerdfn}.  Let 
$\xi_m$ be the vector field on $(L,g_L)=(\mathcal{S}^3,g_{\tau^2})$ dual to $\omega_m$.  Then $c_L$ and 
the Lie derivative $\mathcal{L}_{\xi_1}$ commute for all $\tau>0$, and $c_L$ commutes with $\mathcal{L}_{\xi_m}$ for all $m$ if 
$\tau=1$.
\end{prop} 

\begin{proof}
Using Cartan's formula, we see that 
\begin{align*}
[c_L,\mathcal{L}_{\xi_1}]\alpha&=\xi_1\lrcorner\d *\d\alpha+\d(\xi_1\lrcorner*\d\alpha)-*\d(\xi_1\lrcorner\d\alpha)
\end{align*}
for 1-forms $\alpha$.  Let $\alpha=p_1\omega_1+p_2\omega_2+p_3\omega_3$ and let $\partial_j$ be the operator given in 
Definition \ref{weightdfn} for $j=1,2,3$.  We calculate:
\begin{align*}
*\d\alpha&=\tau(\partial_2p_3-\partial_3p_2-2p_1)\omega_1\\
&\quad+\textstyle\frac{1}{\tau}(\partial_3p_1-\partial_1p_3-2p_2)\omega_2
+\textstyle\frac{1}{\tau}(\partial_1p_2-\partial_2p_1-2p_3)\omega_3,
\intertext{hence} \xi_1\lrcorner\d*\d\alpha&=\Big(\big(\textstyle\frac{1}{\tau}\partial_1\partial_3+2(\tau+\textstyle\frac{1}{\tau})\partial_2\big)p_1+\big(\tau\partial_2\partial_3-\textstyle\frac{4}{\tau}\partial_1\big)p_2\\
&\qquad+\big(\textstyle\frac{4}{\tau}-\textstyle\frac{1}{\tau}\partial_1^2-\tau\partial_2^2\big)p_3\Big)\omega_2\\
&+\Big(\big(-\textstyle\frac{1}{\tau}\partial_1\partial_2+2(\tau+\textstyle\frac{1}{\tau})\partial_3\big)p_1+\big(-\textstyle\frac{4}{\tau}+\textstyle\frac{1}{\tau}\partial_1^2+\tau\partial_3^2\big)p_2\\
&\qquad+\big(-\tau\partial_3\partial_2-\textstyle\frac{4}{\tau}\partial_1\big)p_3\Big)\omega_3
\intertext{and}
\d(\xi_1\lrcorner*\d\alpha)&=\tau(-2\partial_1p_1-\partial_1\partial_3p_2+\partial_1\partial_2p_3)\omega_1+\tau
(-2\partial_2p_1-\partial_2\partial_3p_2+\partial_2^2p_3)\omega_2\\
&\quad+\tau(-2\partial_3p_1-\partial_3^2p_2+\partial_3\partial_2p_3)\omega_3.
\end{align*}
Since
\begin{align*}
\xi_1\lrcorner\d\alpha&=(\partial_1p_2-\partial_2p_1-2p_3)\omega_2-(\partial_3p_1-\partial_1p_3-2p_2)\omega_3,
\intertext{we have that:}
*\d(\xi_1\lrcorner\d\alpha)&=\tau\big([\partial_3,\partial_2]p_1+(-\partial_3\partial_1+2\partial_2)p_2+(\partial_2\partial_1+2\partial_3)p_3\big)\omega_1
\\
&\quad+\textstyle\frac{1}{\tau}\big((\partial_1\partial_3+2\partial_2)p_1-4\partial_1p_2+(4-\partial_1^2)p_3\big)\omega_2\\&\quad+\textstyle\frac{1}{\tau}
\big((-\partial_1\partial_2+2\partial_3)p_1-(4-\partial_1^2)p_2-4\partial_1p_3\big)\omega_3.
\end{align*}
Combining these formulae and using \eq{partialeq} shows that $[c_L,\mathcal{L}_{\xi_1}]\alpha=0$.  Clearly this argument will extend 
to $\mathcal{L}_{\xi_2}$ and $\mathcal{L}_{\xi_3}$ in the case $\tau=1$.
\end{proof}

\begin{prop}\label{cyclicprop}  
Use the notation of Definitions \ref{cLdfn} and \ref{finitegpdfn} and Proposition \ref{Bergerprop}. 
 Let $n\in\Z^+\setminus\{1\}$.  The positive eigenvalues of
 $c_L$ on $(L,g_L)=(\mathcal{S}^3/\Z_n,g_{\tau^2})$ are 
\begin{align*}
\nu_{rn+2l,rn}&=\tau+\sqrt{\tau^2+(rn+2l)(rn+2l+2)+(\textstyle\frac{rn}{\tau})^2(1-\tau^2)}
\quad\text{and}\\
\nu_{sn+n-2}&=\frac{(s+1)n}{\tau}
\end{align*}
for $l,r,s\in\N$. Moreover, 
\begin{gather*}
\sigma_L(\nu_{rn+2l,rn})=2rn+4l+2\quad(r>0), \qquad \sigma_L(\nu_{2l,0})=2l+1,\\ \sigma_L(\nu_{sn+n-2})=2sn+2n-2.
\end{gather*}
\end{prop} 

\begin{proof}
Let $\xi_1$ be the vector field dual to $\omega_1$. 
By Proposition \ref{commprop}, $c_L$ and $\mathcal{L}_{\xi_1}$ commute.  
  Thus, using the notation of Definition \ref{polydfn}, we may restrict attention to $\alpha\in \mathcal{B}_k\otimes\C$ such that 
$$*\d\alpha=-\nu\alpha\quad\text{and}\quad \mathcal{L}_{\xi_1}\alpha=im\alpha$$
for $\nu>0$ and $m\in\Z$.  For $\alpha$ to descend to the quotient $\mathcal{S}^3/\Z_n$, 
we must have that $m\equiv 0$ $(\text{mod }n)$. 
 If we write $\alpha=p_1\omega_1+p_2\omega_2+p_3\omega_3$ then using Cartan's formula and the notation of Definition \ref{weightdfn}
 we find that 
\begin{align*}
im\alpha&=\xi_1\lrcorner\d\alpha+\d(\xi_1\lrcorner\alpha)\\
&=-\nu\xi_1\lrcorner*\alpha+\d\big(\xi_1\lrcorner(p_1\omega_1+p_2\omega_2+p_3\omega_3)\big)\\
&=-\nu\xi_1\lrcorner\big(\textstyle\frac{1}{\tau}p_1\omega_2\w\omega_3+\tau p_2\omega_3\w\omega_1 +\tau p_3\omega_1\w\omega_2\big) +
\d p_1\\
&=\partial_1p_1\omega_1+(\partial_2p_1-\nu\tau p_3)\omega_2+(\partial_3 p_1+\nu\tau p_2)\omega_3.
\end{align*}
For $p_1\neq 0$ we see that $-\partial_1^2p_1=m^2 p_1$, so $\alpha$ is a $\nu_{k,k-2l}$-eigenform for $c_L$ by
 Proposition \ref{Bergerprop}, where $k-2l=rn$ for some $r\in\N$.  For $p_1=0$ we observe that 
 $$-m^2p_2=im(imp_2)=im(-\nu\tau p_3)=-\nu^2\tau^2p_2.$$
Since $\nu\tau=k+2$ by Proposition \ref{Bergerprop}, we have that $k=(s+1)n-2$ for some $s\in\N$, since $k\geq 0$.  
Thus $\alpha$ has eigenvalue 
$\nu_{sn+n-2}$ in these cases.  The eigenvalues and multiplicities now follow from
 Theorem \ref{harmthm} and Proposition \ref{Bergerprop}.
\end{proof}

From the observations in \cite[Example 6.14]{LotayLag} we see that $L_2$ given in Example \ref{cxso3ex} is isometric to 
$\mathcal{S}^3/\Z_2\cong\SO(3)$ with metric $2g_2$.  We therefore calculate the spectrum of $c_{L_2}$ 
 in $(0,4)$ using Proposition \ref{cyclicprop} as follows.

\begin{cor}\label{cxso3cor}
Use the notation of Definitions \ref{cLdfn}, \ref{Bergerdfn} and \ref{finitegpdfn} and Example \ref{cxso3ex}.  The eigenvalues of $c_{L_2}$ on
$(L_2,g_{L_2})=(\mathcal{S}^3/\Z_2,2g_2)$ in $(0,4)$ are: $\{1,2,3,1+\sqrt{5}\}$.
Moreover,
$$\sigma_{L_2}(1)=2,\quad\sigma_{L_2}(2)=7,\quad\sigma_{L_2}(3)=16,\quad\sigma_{L_2}(1+\sqrt{5})=3.$$

Let $\mathcal{C}_{\circ}$ be the family of cones generated by $\GG_2\ltimes\R^7$ transformations of cones 
\begin{equation*}
C_{\mathbf{a}}=\{(0,z_1,z_2,z_3)\in\R\oplus\C^3\cong\R^7\,:\,a_1z_1^2+a_2z_2^2+a_3z_3^2=0\},
\end{equation*}
where $\mathbf{a}\in\mathcal{T}=\{(a_1,a_2,a_3)\in\R^3\,:\,\text{$a_i> 0$ for all $i$ and $a_1+a_2+a_3=1$}\}$.  
Then $C_2=C_{(\frac{1}{3},\frac{1}{3},\frac{1}{3})}$, $\ind(C_2)=5$ and $\ind_{\mathcal{C}_{\circ}}(C_2)=0$, so $C_2$ is $\mathcal{C}_{\circ}$-stable and Jacobi integrable but not stable or rigid.
\end{cor}

\begin{proof}
Using the notation of Proposition \ref{Bergerprop}, let $\lambda_{2r+2l,2r}=\textstyle\frac{1}{\sqrt{2}}\nu_{2r+2l,2r}$ and 
$\lambda_{2s}=\textstyle\frac{1}{\sqrt{2}}\nu_{2s}$, calculated using $\tau=\sqrt{2}$. Then $\lambda_{2r+2l,2r}$ and $\lambda_{2s}$
are the eigenvalues of $c_{L_2}$ by Proposition \ref{cyclicprop}.  

Since $\lambda_{4,4}=4$, we need only consider $\lambda_{2r+2l,2r}$ for $r,l\leq 1$.  We also see that $\lambda_{2s}=s+1$.  
  We therefore calculate 
\begin{gather*}
\lambda_{0,0}=2,\quad\lambda_{2,2}=3,\quad\lambda_{2,0}=1+\sqrt{5},\\
\lambda_0=1,\quad\lambda_2=2,\quad\lambda_4=3.
\end{gather*}
The eigenvalues and multiplicities follow from Proposition \ref{cyclicprop}.

Since $C_2$ is non-planar and the stabilizer of $C_2$ under $\GG_2$ transformations is $\SO(3)$, the dimension of the family 
 of $\GG_2\ltimes\R^7$ transformations of $C_2$ is $7+14-3=18$.  Further,  
 $$\sum_{\lambda\in(-1,1]}\!\!\!\! \dim\check{D}(\lambda)=\dim\check{D}(0)+\dim\check{D}(1)=\sigma_{L_2}(2)+
\sigma_{L_2}(3)=7+16=23,$$
so $\ind(C_2)=5$.  Now, $C_2=C_{(\frac{1}{3},\frac{1}{3},\frac{1}{3})}\in\mathcal{C}_{\circ}$.  Moreover, if 
$\mathbf{a}\neq(\frac{1}{3},\frac{1}{3},\frac{1}{3})$, $C_{\mathbf{a}}$ has trivial stabilizer in $\GG_2\ltimes\R^7$ and has a 
two-parameter family of deformations up to rigid motion given by varying $\mathbf{a}$.  Therefore $\dim\mathcal{C}_{\circ}=7+14+2=23$ and 
$C_2$ is $\mathcal{C}_{\circ}$-stable.  Finally note that $\dim\check{D}(1)=16$, so the space of Lagrangian Jacobi fields on $L_2$ 
is equal to the space of genuine deformations of $L_2$ in the family of links of cones in $\mathcal{C}_{\circ}$, which also shows that $\mathcal{C}_{\circ}$ is a
maximal deformation family for $C_2$.
\end{proof}

We now consider the binary dihedral group. 

\begin{prop}\label{dihedralprop}  
Use the notation of Definitions \ref{cLdfn}, \ref{Bergerdfn} and \ref{finitegpdfn} and Proposition \ref{Bergerprop}.  
For $n\in\Z^+$,   
 the positive eigenvalues of $c_L$ 
on $(L,g_L)=(\mathcal{S}^3/\DD_{n}^*,g_{\tau^2})$ are:
\begin{align*}
\nu_{2rn+2n+2l,2rn+2n}&=\tau\\
&+\sqrt{\tau^2+4(rn+n+l)(rn+n+l+1)+4(\textstyle\frac{rn+n}{\tau})^2(1-\tau^2)};\end{align*}
$$
\nu_{4l,0}=\tau+\sqrt{\tau^2+8l(2l+1)};\quad\text{and}\quad
\nu_{2sn+2n-2}=\frac{2(s+1)n}{\tau}$$
for $l,r,s\in\N$.  Moreover, 
\begin{gather*}
\sigma_L({\nu}_{2rn+2n+2l,2rn+2n})=2rn+2n+2l+1, \qquad \sigma_L(\nu_{4l,0})=4l+1,\\ 
\sigma_L({\nu}_{2sn+2n-2})=4sn+4n-2.
\end{gather*}
\end{prop}

\begin{proof}
Let $\alpha=p_1\omega_1+p_2\omega_2+p_3\omega_3$ be an eigenform of $c_L$
 of positive eigenvalue $\nu$. 
From Definition \ref{finitegpdfn}(b), $\Z_{2n}\subseteq\DD_{n}^*$.  Hence, by Proposition \ref{cyclicprop}, the possible $\nu$ are of the form $\nu_{k,k-2l}$ where $k-2l\equiv 0$ (mod $2n$) or $\nu_{2sn+2n-2}$, using the notation of Proposition \ref{Bergerprop}. 

Recalling Definition \ref{cxpolydfn}, we see that $\bfj$ sends $(z_1,z_2)\in\C^2$ to $(z_2,-z_1)$ and hence $\bfj$ maps
 $\mathcal{R}_m^{\C}$ to $\mathcal{R}_{-m}^{\C}$ for each $m\in\Z$.  
By note (a) after Proposition \ref{Bergerprop} and Definition \ref{finitegpdfn}(b), the multiplicity of $\nu_{2rn+2l,2rn}$
 is the number of choices of $p_1\in\mathcal{Q}_{2rn+2l,2rn}$ which
  are $\bfj$-invariant.  By Theorem \ref{harmthm}, $p_1$ is the sum of real
  and imaginary parts of polynomials in $\mathcal{Q}_{2rn+2l}^{\C}\cap\mathcal{R}_{2rn}^{\C}$.  Since $\bfj^2$ clearly acts as 
  the identity on $\mathcal{Q}_{2rn+2l}^{\C}$, we can decompose $\mathcal{Q}_{2rn+2l,2rn}$ into $\pm1$-eigenspaces for $\bfj$.   
Thus, for $r>0$, the subspace of $\mathcal{Q}_{2rn+2l,2rn}$ which
 is $\bfj$-invariant is half the total dimension.  For $r=0$, it is
 straightforward to see that $\bfj$ acts as $(-1)^l$ on 
 $\mathcal{Q}_{2l}^{\C}\cap\mathcal{R}_0^{\C}$, so $p_1\in\mathcal{Q}_{2l,0}$ is $\bfj$-invariant if and only if $l$ is even.  

By note (b) after Proposition \ref{Bergerprop}, the multiplicity of $\nu_k$, if $k\neq 0$, is determined by the number of
choices for $p_2$.  Now, since $p_2$ is only well-defined up to sign on the quotient of $\mathcal{S}^3$ by $\DD_n^*$, we need to 
calculate the number of polynomials in $\mathcal{Q}_{k,k}$ on which $\bfj$ acts as $\pm1$.  However, since $k=2sn+2n-2>0$ is even,
 we can find a basis for $\mathcal{Q}_{k,k}$ consisting of $\pm1$-eigenfunctions for $\bfj$ by the previous paragraph. 
  For $\nu_0$ the multiplicity is 2 because $p_2$ and $p_3$ are constant. 

The result follows from Theorem \ref{harmthm} and Proposition \ref{Bergerprop}.
\end{proof}

By \cite[Example 6.15]{LotayLag}, $L_3$ given in Example \ref{so3ex} is isometric to $\SO(3)/\SSS_3\cong\mathcal{S}^3/\DD_{3}^*$
 with metric $6g_6$.  We may thus apply Proposition \ref{dihedralprop} as follows.

\begin{cor}\label{so3cor1}
Use the notation of Definitions \ref{cLdfn}, \ref{Bergerdfn} and \ref{finitegpdfn} and Example \ref{so3ex}.  The eigenvalues of $c_{L_3}$ on
$(L_3,g_{L_3})=(\mathcal{S}^3/\DD_{3}^*,6g_6)$ in $(0,4)$ are: $\{1,2,3,1+\sqrt{5}\}$.
Moreover,
$$\sigma_{L_3}(1)=10
,\quad\sigma_{L_3}(2)=23
,\quad\sigma_{L_3}(3)=41
,
\quad\sigma_{L_3}(1+\sqrt{5})=5.$$
Hence, the  
cone 
$C_3$ 
is not stable and $\ind(C_3)=46$.  Moreover, $C_3$ is not rigid.
\end{cor}

\begin{proof}
By Proposition \ref{dihedralprop} applied with $\tau=\sqrt{6}$, the eigenvalues of $c_{L_3}$ are 
$\lambda_{6r+6+2l,6r+6}=\frac{1}{\sqrt{6}}\nu_{6r+6+2l,6r+6}$, $\lambda_{4l,0}=\frac{1}{\sqrt{6}}\nu_{4l,0}$ and $\lambda_{6s+4}=\frac{1}{\sqrt{6}}\nu_{6s+4}$, for $l,r,s\in\N$. 
Since $\lambda_{8,0}=1+\textstyle\frac{1}{3}\sqrt{129}>4$, $\lambda_{8,6}=1+\textstyle\frac{1}{3}\sqrt{84}>4$ and $\lambda_{12,12}=4$, we need only calculate:
\begin{gather*}
\lambda_{0,0}=2,\quad\lambda_{4,0}=1+\sqrt{5},\quad\lambda_{6,6}=3,\quad\lambda_{6s+4}=s+1.
\end{gather*}
The multiplicities now follow from Proposition \ref{dihedralprop}.

Since $C_3$ is non-planar and has $\SO(3)$ stabilizer in $\GG_2$, 
\begin{align*}
\ind(C_3)&=\sum_{\lambda\in(-1,1]}\!\!\!\! \dim\check{D}(\lambda)-7-(14-3)=\dim\check{D}(0)
+\dim\check{D}(1)-18\\&=\sigma_{L_3}(2)+\sigma_{L_3}(3)-18=23+41-18=46
\end{align*}
as claimed.  We also see that $\dim\check{D}(1)=41>11=\dim\GG_2-\dim\SO(3)$, so $C_3$ is not rigid.
\end{proof}

Lastly, we study the constant curvature 1 metric on $\mathcal{S}^3/\AAA_4^*$.

\begin{prop}\label{tetraprop}  Use the notation of Definitions \ref{cLdfn} and \ref{finitegpdfn}.  
Let $(L,g_L)=(\mathcal{S}^3/\AAA_4^*,g_{1})$.  The positive eigenvalues of $c_L$ are 
$$\nu_{2r}=2r+2\quad\text{and}\quad\nu_{6s+4,6s+4}=6s+6$$
 for $r,s\in\N$.  Moreover,
$$\sigma_L(\nu_{2r})=(2r+1)(1+2[\textstyle\frac{r}{3}]+[\textstyle\frac{r}{2}]-r)\quad\text{and}
\quad\sigma_L(\nu_{6s+4,6s+4})=12s+10 
,$$
where $[q]$ denotes the integer part of a non-negative rational number $q$.
\end{prop}

\begin{proof}
From Definition \ref{finitegpdfn}(c), we observe that $\AAA_4^*$ is generated by $\bfi$ and $\bfj$, which are elements of order $4$,
 and
$$\frac{1}{2}(\bfo+\bfi+\bfj+\bfk)=\frac{1}{2}\bfo+\frac{\sqrt{3}}{2}\left(\frac{\bfi+\bfj+\bfk}{\sqrt{3}}\right)$$
 which has order $6$.  By Proposition \ref{Bergerprop} and Corollary \ref{roundcor}, we know that the possible 
 positive eigenvalues of $c_L$ are $k+2$ for $k\in\N$ with corresponding eigenforms $\alpha\in\mathcal{B}_k$, in 
the notation of Definition \ref{polydfn}.

Suppose $\alpha=p_1\omega_1+p_2\omega_2+p_3\omega_3\in\mathcal{B}_k\otimes\C$ is an eigenform of $c_L$.  We may write 
$\alpha=P_1\Omega_1+P_2\Omega_2+P_3\Omega_3$ where 
\begin{align*}
\Omega_1&=\frac{\omega_1+\omega_2+\omega_3}{\sqrt{3}}, & \Omega_2&=\frac{\omega_1-\omega_3}{\sqrt{2}}, 
& \Omega_3&=\frac{-\omega_1+2\omega_2-\omega_3}{\sqrt{6}},\\[4pt]
P_1&=\frac{p_1+p_2+p_3}{\sqrt{3}}, & P_2&=\frac{p_1-p_3}{\sqrt{2}}, & P_3&=\frac{-p_1+2p_2-p_3}{\sqrt{6}}.
\end{align*}

Let $\xi_m$ be the vector field dual to $\omega_m$ and set $\xi=(\xi_1+\xi_2+\xi_3)/\sqrt{3}$. 
Since $\xi$ commutes with $c_L$ by Proposition \ref{commprop},  
we may look for $\alpha$ satisfying 
$$c_L\alpha=(k+2)\alpha\quad\text{and}\quad\mathcal{L}_{\xi}\alpha=im\alpha$$
for $m\in\Z$.  For $\alpha$ to be well-defined on $L$ we must have that $m\equiv 0$ (mod 6).  Using \eq{structeq2} we calculate:
\begin{align*}
\mathcal{L}_{\xi}\Omega_1&=\textstyle\frac{1}{3}(\mathcal{L}_{\xi_2+\xi_3}\omega_1+\mathcal{L}_{\xi_3+\xi_1}\omega_2+
\mathcal{L}_{\xi_1+\xi_2}\omega_3)\\
&=\textstyle\frac{1}{3}(-2\omega_3+2\omega_2-2\omega_1+2\omega_3-2\omega_2+2\omega_1)=0;\\
\mathcal{L}_{\xi}\Omega_2&=2\Omega_3;\\
\mathcal{L}_{\xi}\Omega_3&=-2\Omega_2.
\end{align*}

Recall the notation of Definitions \ref{weightdfn} and \ref{cxpolydfn}.  By similar methods to the proof of 
 Proposition \ref{cyclicprop}, $P_1$ has $\mathbf{q}$-weight $6l$ for some 
$l\in\N$, where $\mathbf{q}=\frac{\bfi+\bfj+\bfk}{\sqrt{3}}$.  Moreover, $P_1$ must be invariant under $\bfi$ and $\bfj$, so $P_1$ is  
$\AAA_4^*$-invariant.  The $\AAA_4^*$-invariant eigenfunctions of the Laplacian on $\mathcal{S}^3$ are effectively determined in 
\cite[Theorem 4.4]{Ikeda}, so we deduce that $k=2r$ and $P_1$ lies in a subspace of $\mathcal{Q}_{2r}^{\C}$ of dimension $(2r+1)(1+2[\frac{r}{3}]+[\frac{r}{2}]-r)$. 
The formulae for $\nu_{2r}$ and $\sigma_L(\nu_{2r})$ follow.

Again following the proof of Proposition \ref{cyclicprop}, we observe that $P_1$ determines $P_2$ and $P_3$ unless $P_2,P_3\in\mathcal{Q}_{6s+4}^{\C}$ have $\mathbf{q}$-weight $6s+4$ for some $s\in\N$.  Thus, 
we also have eigenvalues of the form $6s+6$ with multiplicity equal to the number of choices for $P_2$.  
As in the proof of Proposition \ref{dihedralprop}, $P_2$ is only well-defined up to sign on the quotient of $\mathcal{S}^3$. 
 Since $P_2$ is a polynomial of even degree it has even $\bfi$-weight and $\bfj$-weight, so $P_2$ maps to $\pm P_2$ under 
 the action of $\bfi$ and $\bfj$.  Since  
we may decompose the harmonic polynomials on $\mathcal{S}^3$ by $\mathbf{q}$-weight rather than $\bfi$-weight as in Theorem 
\ref{harmthm}, and the metric on $L$ has constant curvature, we deduce that $\sigma_L(6s+6)=\dim\mathcal{Q}_{6s+4,6s+4}=12s+10$.  
\end{proof}

As observed in \cite[Example 6.6]{LotayLag}, the link $L_4$ of $C_4$, given in Example \ref{so3ex}, is isometric to
 $\SO(3)/\AAA_4\cong\mathcal{S}^3/\AAA_4^*$ with constant curvature $\frac{1}{16}$.  We may therefore deduce the following.

\begin{prop}
Use the notation of Definitions \ref{cLdfn}, \ref{Bergerdfn} and \ref{finitegpdfn} and Example \ref{so3ex}.
  The eigenvalues of $c_{L_4}$ on $(L_4,g_{L_4})=(\mathcal{S}^3/\AAA_4^*,16g_1)$ in $(0,4)$ are: 
  $\{\frac{1}{2},\frac{3}{2},2,\frac{5}{2},3,\frac{7}{2}\}$.  Moreover,
\begin{align*}
\sigma_{L_4}(\textstyle\frac{1}{2})&=1,&\sigma_{L_4}(\textstyle\frac{3}{2})&=10
,&\sigma_{L_4}(2)&=7,\\
\sigma_{L_4}(\textstyle\frac{5}{2})&=9,&\sigma_{L_4}(3)&=22
,&\sigma_{L_4}(\textstyle\frac{7}{2})&=26.
\end{align*}
Hence, the cone $C_4$ is not stable and $\ind(C_4)=30$.  Moreover, $C_4$ is not rigid.
\end{prop}

\begin{proof}
By Proposition \ref{tetraprop}, the positive eigenvalues of $c_{L_4}$ are $\lambda_{2r}=\frac{1}{2}(r+1)$ and 
$\lambda_{6s+4,6s+4}=\frac{3}{2}(s+1)$.  We are restricted to $r=0,3,4,6$ and $s=0,1$ for 
 eigenvalues in $(0,4)$ since the multiplicity of $\lambda_{2r}$ for $r=1,2,5$ is zero.  

Hence, on $C_4$,
\begin{align*}
\sum_{\lambda\in(-1,1]}\!\!\!\! \dim\check{D}(\lambda)&
=\dim\check{D}(-\textstyle\frac{1}{2})
+\dim\check{D}(0)+\dim\check{D}(\textstyle\frac{1}{2})+\dim\check{D}(1)\\ 
&=\sigma_{L_4}(\textstyle\frac{3}{2})+\sigma_{L_4}(2)+\sigma_{L_4}(\textstyle\frac{5}{2})+\sigma_{L_4}(3)=10+7+9+22=48.
\end{align*}
Since $C_4$ is non-planar and $\SO(3)$-invariant, the dimension of the family of $\GG_2\ltimes\R^7$ transformations of $C_4$ 
has dimension $7+14-3=18$.  Thus $\ind(C_4)=30$.  
Observe that $\dim\check{D}(1)=22>11=\dim\GG_2-\dim\SO(3)$, so 
$C_4$ is not rigid.  
\end{proof}

The stability index is a measure of the geometry complexity of the coassociative cone.
We might therefore expect homogeneous cones to have the greatest chance of being stable.  Since we find that the only
stable homogeneous coassociative cones are 4-planes and the cone in Example \ref{u2ex}, we might naively expect that 
these are the only stable coassociative cones.  

In Geometric Measure Theory one defines the $m$-dimensional \emph{density} 
of a set 
$S\subseteq\R^n$ at a point $p$ by 
$$\Theta(S,p)=\lim_{r\rightarrow 0}\frac{\mathcal{H}^m\big(S\cap B_n(p;r)\big)}{\vol\big(B_m(0;r)\big)},$$
where $\mathcal{H}^m$ is $m$-dimensional Hausdorff measure and $B_n(p;r)$ is the ball  in $\R^n$ of radius $r$ about $p$.  
If $S$ is an $m$-dimensional submanifold then $\Theta(S,p)=1$ for all $p\in S$.
For coassociative cones $C$ in $\R^7$ with isolated singularities at the origin and compact links $L$ it is straightforward 
  to calculate the 4-dimensional density as $\Theta(C,0)=\vol(L)/\vol(\mathcal{S}^3)$.  We may therefore easily calculate 
  the density for the homogeneous coassociative cones given in Examples \ref{planeex}-\ref{so3ex} as follows:
$$\Theta(C_0,0)=1,\!\quad\Theta(C_1,0)=\textstyle\frac{16}{9},\!\quad\Theta(C_2,0)=4,\!\quad\Theta(C_3,0)=36,\!\quad\Theta(C_4,0)=64.$$
The fact that $C_0$ and $C_1$ are the homogeneous cones with the lowest density and are also the only stable ones is 
 suggestive.

\section{Algebraic curves and 2-ruled cones}\label{algcurvessec}

Since Example \ref{2rex} gives the largest known family of coassociative cones, we are motivated to analyse the 
stability index for 2-ruled cones.  We therefore need to understand the curl operator on geodesic 
$\mathcal{S}^1$-bundles over surfaces.  We are particularly interested in the case where the surface is an algebraic curve.  

We start with some definitions from the theory of Riemannian submersions.

\begin{dfn}\label{horizontaldfn}
Let $(L,g_L)$ be a compact Riemannian 3-manifold which is an $\mathcal{S}^1$-bundle $\pi:L\rightarrow\Sigma$ over a 
 compact Riemannian surface $(\Sigma,g_\Sigma)$.  
Suppose further that $\pi$ is a Riemannian submersion.  
  Let $\xi$ be a unit vector field spanning the vertical distribution of $\pi$ and let $\theta$ be the 1-form dual to 
  $\xi$.  

We define a form $\alpha$ on $L$ to be \emph{horizontal} if $\xi\lrcorner\alpha=0$ and we denote the bundle of horizontal 
  $m$-forms by $\Lambda^m_hT^*L$.  Trivially all functions are horizontal.  We can identify the horizontal $m$-forms with 
$m$-forms on $\Sigma$.   We define the \emph{horizontal Hodge star} $*_h$ 
on horizontal forms $\alpha$ via the equation: $*\alpha=*_h\alpha\w\theta$.

By \cite[Equation (2.1.2)]{Nagy}, 
$$\Lambda^{m+1}T^*L=\Lambda^{m+1}_hT^*L\oplus(\Lambda^{m}_hT^*L\w\theta),$$ so 
we can define a \emph{horizontal derivative} $\d_h:C^{\infty}(\Lambda^mT^*L)\rightarrow C^{\infty}(\Lambda^{m+1}_hT^*L)$ by sending 
$\alpha$ to the horizontal part of $\d\alpha$.  
Since $\d_h$ sends any section of $\Lambda^{m}_hT^*L\w\theta$ to zero, $\d_h:C^{\infty}(\Lambda^m_hT^*L)\rightarrow C^{\infty}(\Lambda^{m+1}_hT^*L)$ is an antiderivation.  We can therefore define the 
 formal adjoint $\d_h^*$ of $\d_h$ and the \emph{horizontal Laplacian}
 $\Delta_h=\d_h\d_h^*+\d_h^*\d_h$.  The horizontal Laplacian is not necessarily elliptic.

If $\Delta_{L}$ is the ordinary Laplacian on $L$, we call $\Delta_v=\Delta_L-\Delta_h$ the \emph{vertical Laplacian}.  Notice in the case of functions that $\Delta_v=-\mathcal{L}_{\xi}^2$.
\end{dfn}

\noindent The situation above includes the Berger 3-spheres, 
tubes over pseudoholomorphic curves in $\mathcal{S}^6$ and real links of 
complex 2-dimensional cones. 

The key result for understanding the Laplacian on functions in the situation of Definition \ref{horizontaldfn} is the following \cite[Theorem 1.5]{Berard}.

\begin{prop}\label{oldcommprop}
Use the notation of Definition \ref{horizontaldfn}.  If $\pi:L\rightarrow\Sigma$  has totally geodesic fibres then, on functions,  $\Delta_{L}$, $\Delta_h$ and $\Delta_v$ commute.
\end{prop}

For convenience we make the following definition.

\begin{dfn}\label{Lapmultdfn}
Let $L$ be a compact Riemannian 3-manifold.  Let $m_L(\nu)$ denote the multiplicity of the eigenvalue $\nu$ of the 
Laplacian $\Delta_L$ on functions.  Let $\mathcal{E}_L=\{\nu\in\R\,:\,m_L(\nu)\neq 0\}$, which is a countable discrete set.
\end{dfn}

\begin{remarks}
Use the notation of Definition \ref{horizontaldfn} and suppose that $\pi:L\rightarrow\Sigma$
has totally geodesic fibres.  Let $m_L(\nu_h,\nu_v)$ denote the dimension of the space of functions 
$f$ such that $\Delta_h f=\nu_hf$ and $\Delta_v f=\nu_vf$.  Notice that,
since $\Delta_h$ and $\Delta_v$ are non-negative operators, we must have $\nu_h\geq 0$ and $\nu_v\geq 0$.  Moreover,
 by one of the main results in \cite{Watson}, $m_L(\nu_h,0)$ is the multiplicity of the eigenvalue $\nu_h$ of the Laplacian
 $\Delta_{\Sigma}$ acting on functions on $\Sigma$.   Finally, by Proposition \ref{oldcommprop}, every eigenvalue of $\Delta_L$ is of
 the form $\nu_h+\nu_v$, so $m_L(\nu)\leq\sum_{\nu_h+\nu_v=\nu}m_{L}(\nu_h,\nu_v)$.
\end{remarks}  

We now have the following result which is similar to Proposition \ref{commprop}.

\begin{prop}\label{commprop2}
Use the notation of Definitions \ref{cLdfn} and \ref{horizontaldfn}.  If the fibres of $\pi:L\rightarrow\Sigma$ 
are totally geodesic and $\xi$ is a Killing vector field for $g_L$, then $c_L$ and the Lie derivative $\mathcal{L}_{\xi}$ commute.  
\end{prop}

\begin{proof} Recall Definition \ref{horizontaldfn}.  
By the work in \cite[$\S$2.1]{Nagy}, since $\nabla_{\xi}\xi=0$,
$\d\theta
$ is horizontal.  Moreover, if we write a 1-form $\gamma$ on $L$ as $\gamma=\alpha+f\theta$ with $\alpha$ horizontal and $f$ a function, then
$$*\d\gamma=*_h(\mathcal{L}_{\xi}\alpha-\d_hf 
)+(*_h\d_h\alpha+f\!*_h\!\d\theta)\theta\quad\text{and}\quad\mathcal{L}_{\xi}\gamma=
\mathcal{L}_{\xi}\alpha
+(\mathcal{L}_{\xi}f)\theta.$$
Thus, since $\xi$ is Killing and $\d\theta$ is horizontal, 
it is straightforward to compute:
$$*\d(\mathcal{L}_\xi\gamma)=\big(*_h\!\mathcal{L}_\xi^2\alpha-*_h\d_h(\mathcal{L}_\xi f)\big)+\big(*_h\!\d_h(\mathcal{L}_\xi\alpha)+
(\mathcal{L}_{\xi}f)*_h\!\d\theta\big)
\theta$$
and 
$$\mathcal{L}_\xi(*\d\gamma)=\big(*_h\!\mathcal{L}_\xi^2\alpha-*_h\mathcal{L}_\xi(\d_hf)\big)+\big(*_h\!\mathcal{L}_\xi(\d_h\alpha)+
(\mathcal{L}_{\xi}f)*_h\!\d\theta\big)\theta.$$
 Note that, since the fibres of $\pi$ are totally geodesic, 
$[\d_h,\mathcal{L}_{\xi}]=0$ on horizontal forms by \cite[Equation (2.1.6)]{Nagy}.  
  We conclude that $[c_L,\mathcal{L}_{\xi}]=0$ as claimed. 
\end{proof}

Observe that, by Example \ref{Killingex}, we have two types of links of coassociative cones where Definition
 \ref{horizontaldfn} and Propositions \ref{oldcommprop} and \ref{commprop2} apply: namely, Hopf lifts of 
holomorphic curves in $\C\P^2$ and tubes of radius $\frac{\pi}{2}$ in the first or second normal bundle about a null-torsion 
pseudoholomorphic curve in $\mathcal{S}^6$.  We may therefore try to describe the spectrum of the curl operator in these cases.

We begin with the links of complex 2-dimensional cones.

\begin{thm}\label{holcurvethm}
Let $\Sigma$ be a compact, connected, holomorphic curve in $\C\P^2$ with degree $d_\Sigma$.  Let $L$ be the Hopf lift
 of $\Sigma$ in $\mathcal{S}^5$.  Use the  notation of Definitions \ref{cLdfn} and \ref{Lapmultdfn}.  The eigenvalues of $c_L$ in
 $(0,4)$ are $1$, $2$, $3$ and $\lambda+2\in(2,4)\setminus\{3\}$ such that $\lambda(\lambda+2)\in\mathcal{E}_L$.  Moreover, 
\begin{align*}\sigma_L(1)&=d_\Sigma^2-d_\Sigma, &\sigma_L(2)&=d_\Sigma^2+d_\Sigma+1,\\
\sigma_L(3)&=m_L(3)+d_\Sigma^2+3d_\Sigma,&\sigma_L(\lambda+2)&=m_L\big(\lambda(\lambda+2)\big).
\end{align*}
Furthermore, $m_L(3)\geq6$ if $d_\Sigma\geq 2$. 
\end{thm}

\begin{proof}
Notice that $L$ is a $\U(1)$-bundle over $\Sigma$, and that the projection $\pi:L\rightarrow\Sigma$ is a Riemannian submersion with
 totally geodesic fibres when $L$ and $\Sigma$ are given the induced metrics from the standard metrics on $\mathcal{S}^5$ and 
$\C\P^2$.   Therefore, Definition \ref{horizontaldfn} and Proposition \ref{oldcommprop} apply.  If we let $\xi$ be the vector field given by the 
$\U(1)$ action then, as observed in Example \ref{Killingex}, $\xi$ is a Killing vector field for the metric on $L$, so Proposition 
\ref{commprop2} applies.  Hence, to understand the positive eigenvalues $\nu$ of $c_L$, it is enough to study the equations
\begin{align}\label{stabpropeq1}
\d\gamma&=-\nu*\gamma\quad\text{and}\\
 \mathcal{L}_{\xi}\gamma &=im\gamma\label{stabpropeq2}
\end{align}
for $\gamma\in C^{\infty}(T^*L\otimes\C)$ and $m\in\Z$. 

Let $\omega_0$ be the standard K\"ahler form on $\C^3$.  Since $\Sigma$ is holomorphic, $\omega_0|_\Sigma=\vol_\Sigma$ and thus $\omega=\omega_0|_L$ is a nowhere vanishing 2-form on $L$.  Moreover, $*\omega$ is the 1-form on $L$ dual to 
the $\U(1)$ vector field $\xi$, where $*$ is the Hodge star on $L$, so $\theta=*\omega$ in the notation 
of Definition \ref{horizontaldfn}.  

By \eq{phispliteq}, if $x_1$ is the coordinate on $\R$ in the decomposition $\R^7=\R\oplus\C^3$ and $C$ is the cone on $L$ embedded 
in $\R^7$, then $\beta=\frac{\partial}{\partial x_1}\lrcorner\varphi|_C=\omega_0|_C$ is a self-dual 2-form 
which corresponds to the deformation of $C$ by translation in the $x_1$ direction.  By Corollary \ref{coasscor} and 
Proposition \ref{jmathprop}, $\beta$ must be closed.  
Notice that $\omega_0|_C=r^2\omega+r*\omega\w\d r$, where $r$ is 
the radial coordinate, and so $\d\beta=0$ if and only if $\d\!*\!\omega=-2\omega$.  

By Definition \ref{horizontaldfn}, given $\gamma\in C^{\infty}(T^*L\otimes\C)$,  
there exist $f\in C^{\infty}(L)$ and $\alpha\in C^{\infty}(T^*_hL\otimes\C)$ such that
$\gamma=f*\omega+\alpha$.  Thus, 
\begin{align*}
\d\gamma&=\d (f*\omega)+\d\alpha =\d_h f\w*\omega-2f\omega+\d\alpha\quad\text{and}\quad
*\gamma=f\omega+*_h\alpha\w*\omega
\end{align*}
We deduce from looking at the horizontal components of \eq{stabpropeq1} that
\begin{equation}\label{stabpropeq3}
\d_h\alpha=-(\nu-2) f\omega.
\end{equation}
For convenience we set $\d_v=\d-\d_h$.  Hence, using Cartan's formula and the facts that $\xi\lrcorner\alpha=\xi\lrcorner\omega=0$ and 
$\xi\lrcorner*\omega=1$, \eq{stabpropeq2} becomes:
\begin{align*}
im\gamma&=\d(\xi\lrcorner\gamma)+\xi\lrcorner\d\gamma=\d f + \xi\lrcorner\d(f*\omega)+\xi\lrcorner\d\alpha\\
&=\d f+\xi\lrcorner(\d_h f\w *\omega)+\xi\lrcorner\d_v\alpha=\d_vf +\xi\lrcorner\d_v\alpha.
\end{align*}
We can also use \eq{stabpropeq1} in \eq{stabpropeq2} to see that:
\begin{align*}
im\gamma&=\d f-\nu\xi\lrcorner*\gamma=\d f-\nu f\xi\lrcorner \omega - \nu\xi\lrcorner(*_h\alpha\w*\omega)=\d f+\nu*_h\alpha.
\end{align*}
Together, we deduce that 
\begin{align}
\d_h f&=im\alpha-\nu*_h\alpha, \label{stabpropeq4}\\
\d_vf&=imf*\omega,\label{stabpropeq5}\\
\d_v\alpha&=im*\omega\w\alpha\label{stabpropeq6}.
\end{align}
Since $\nu>0$, \eq{stabpropeq1} forces $\d^*\gamma=0$ and thus 
\begin{equation}\label{stabpropeq7}
\d_h*_h\alpha=-imf\omega.
\end{equation}

As in \cite[$\S$2.2]{Nagy}, the complex structure $J$ on $\Sigma$ defines a complex structure map $J_h$ on horizontal forms, 
which agrees with $J$ on horizontal 1-forms (viewed as lifts of 1-forms on $\Sigma$).  This
allows us to define $\Lambda^{p,q}_hT^*L\otimes\C$, where $p,q\in\N$, to be the bundle of horizontal $(p+q)$-forms 
which have eigenvalue $i(p-q)$ under $J_h$, and thus give the decomposition
 $\Lambda^m_hT^*L\otimes\C=\oplus_{p+q=m}\Lambda^{p,q}_hT^*L\otimes\C$.  We may therefore define operators 
\begin{align*}
\partial_h&:C^{\infty}(\Lambda^{p,q}_hT^*L\otimes\C)\rightarrow C^{\infty}(\Lambda^{p+1,q}_hT^*L\otimes\C) 
\quad\text{and}\\
 \bar{\partial}_h&:C^{\infty}(\Lambda^{p,q}_hT^*L\otimes\C)\rightarrow C^{\infty}(\Lambda^{p,q+1}_hT^*L\otimes\C)
 \end{align*}
by taking appropriate components of $\d_h$.

Since $\alpha$ is a horizontal complex 1-form on $L$, it can be decomposed into $(1,0)$ and $(0,1)$ components.   
Thus, $*_h\alpha=*_h(\alpha^{1,0}+\alpha^{0,1})=i\alpha^{1,0}-i\alpha^{0,1}$ and so \eq{stabpropeq4} becomes:
\begin{align}
\partial_h f&=i(m-\nu)\alpha^{1,0},\label{stabpropeq4a}\\
\bar{\partial}_h f&=i(m+\nu)\alpha^{0,1}.\label{stabpropeq4b}
\end{align}
We notice that $f=0$ forces $m^2=\nu^2$ otherwise $\alpha=0$ by \eq{stabpropeq4a}-\eq{stabpropeq4b}.
Furthermore, \eq{stabpropeq3} for $f=0$ becomes
\begin{equation}\label{stabpropeq3a}
\bar{\partial}_h\alpha^{1,0}+\partial_h\alpha^{0,1}=0.
\end{equation}
Thus, $f=0=m-\nu$ is equivalent to $\bar{\partial}_h\alpha^{1,0}=0=\alpha^{0,1}$ and $f=0=m+\nu$ is equivalent to 
$\partial_h\alpha^{0,1}=0=\alpha^{1,0}$.   

The hyperplane bundle over $\C\P^2$ defines a complex line bundle $H$ over $\Sigma$, which is also a real line bundle over $L$.  
We can, of course, identify $H$ over $L$ with the cone $C$, and so $H$ has a global section $s$ over $L$ given simply by 
$s(x)=x$.  It is clear that $\mathcal{L}_\xi s=is$ and $\alpha\otimes s^{-m}$, by \eq{stabpropeq6}, is a $\U(1)$-invariant section
 of $T^*L\otimes\C\otimes H^{-m}$, and so pushes 
down to be a well-defined section of $T^*\Sigma\otimes H^{-m}$ over $\Sigma$.  The condition $\bar{\partial}_h\alpha^{1,0}=0$ given 
by \eq{stabpropeq3a}, when $m=\nu$, is then equivalent to saying that $\alpha^{1,0}\otimes s^{-\nu}$ is a holomorphic section of
 $P_\nu=T^{*{1,0}}\Sigma
\otimes H^{-\nu}$.  Since $H$ is the bundle $\mathcal{O}_\Sigma(-1)$ and $T^{*1,0}\Sigma\cong \mathcal{O}_\Sigma(d_\Sigma-3)$, 
we see that $P_\nu\cong\mathcal{O}_\Sigma(\nu+d_\Sigma-3)$.  Thus, by Riemann--Roch, 
the dimension of the vector space of holomorphic sections of $P_\nu$ is 
$$h^0(P_\nu)=d_\Sigma(\nu+d_\Sigma-3)+1-g_\Sigma.$$
(Notice that $P_3$ is isomorphic to the normal bundle of $\Sigma$ in $\C\P^2$ and so its holomorphic sections correspond precisely to
 the infinitesimal deformations of $\Sigma$ as a holomorphic curve.)  
Similarly, for $m=-\nu$, we have that $\alpha^{0,1}\otimes s^{\nu}$ is an anti-holomorphic section of $\overline{P_{\nu}}$.  
Hence, for $f=0$, we get integer eigenvalues $\nu\in\{1,2,3\}$ for $c_L$ in $(0,4)$ (since $m^2=\nu^2$) and contributions of $2h^0(P_{\nu})$ to $\sigma_L(\nu)$.   
Using the degree-genus formula, we calculate
\begin{equation}\label{holseceq}
2h^0(P_{\nu})=d_{\Sigma}(d_\Sigma-3+2\nu).
\end{equation} 

It follows from \eq{stabpropeq3}-\eq{stabpropeq7} that
\begin{align}
\Delta_L f&=*\d\!*\big(imf*\omega+im\alpha-\nu*_h\alpha\big)\nonumber\\
&
=-*\d\big(imf\omega+im*_h\alpha\w*\omega+\nu\alpha\w*\omega\big)\nonumber\\
&=-*\big(im(imf)\vol_L+im(-imf)\vol_L-\nu(\nu-2)f\vol_L\big)\nonumber\\
&
=\nu(\nu-2)f.\label{stabpropeq8}
\end{align}
Since $\mathcal{L}_{\xi}f=imf$ by \eq{stabpropeq4}-\eq{stabpropeq5}, we deduce from Definition \ref{horizontaldfn} 
and Proposition \ref{oldcommprop} that
\begin{equation}\label{stabpropeq9}
\Delta_h f=\big(\nu(\nu-2)-m^2\big)f.
\end{equation}
Recall that the eigenvalues of $\Delta_L$ and $\Delta_h$ have to be non-negative.

If $\nu\in(0,2)$ then the only solution 
of \eq{stabpropeq8} is $f=0$.  Thus, $m^2=\nu^2$ by \eq{stabpropeq4a}-\eq{stabpropeq4b}, so $\nu=1$ 
and $\sigma_L(1)=d_{\Sigma}(d_\Sigma-1)$ from \eq{holseceq} as claimed. 

If $\nu=2$ and $f\neq 0$ then $m=0$ by \eq{stabpropeq9} and the solutions of \eq{stabpropeq8} consist of constant (non-zero)
functions.  Since $f\neq 0$ is constant, we see from \eq{stabpropeq4a}-\eq{stabpropeq4b} that $\alpha=0$.  Thus, 
we have a 1-dimensional space of solutions to \eq{stabpropeq1}-\eq{stabpropeq2} for $\nu=2$ and $f\neq 0$. 
We deduce the formula for $\sigma_L(2)$ from \eq{holseceq}. 

If $\nu\in(2,4)$ and $f\neq 0$ then the non-negativity of $\nu^2-2\nu-m^2$ by \eq{stabpropeq9} forces $|m|\leq 2$.  Moreover, 
solutions to \eq{stabpropeq8} define $\alpha$ via \eq{stabpropeq4a}-\eq{stabpropeq4b} since $m^2\neq\nu^2$.  We deduce the remaining 
eigenvalues of $c_L$ and multiplicities from \eq{holseceq} as claimed.  

Now suppose $d_{\Sigma}>1$ so that the cone $C$ is non-planar.  We can view the Lie algebra of $\SU(3)$ as a subalgebra of $\mathfrak{g}_2$ and 
decompose $\mathfrak{g}_2=\mathfrak{su}(3)\oplus\mathfrak{su}(3)^{\perp}$.  If $v\in\mathfrak{su}(3)^{\perp}$ is identified with 
a tangent vector on $\R^7$, $v|_C$ defines a normal vector $w$ which in turn defines an infinitesimal deformation of $C$ as a 
coassociative cone.  Moreover, since $v\in\mathfrak{su}(3)^{\perp}$ and $C$ is non-planar, $w$ is a non-trivial infinitesimal
 deformation which is not complex.  Thus $(w\lrcorner\varphi)|_C$ will be an order $O(r)$ self-dual 2-form which does not arise purely
 from a horizontal 1-form on $L$, and will define a $3$-eigenform of $c_L$ with $f\neq 0$ (since, for $f=0$, the 3-eigenforms 
correspond to the infinitesimal deformations of $\Sigma$ as a holomorphic curve).  Since $\dim\mathfrak{su}(3)^{\perp}=6$,
 we deduce that $m_L(3)\geq 6$ as claimed.
\end{proof}

It is easy to check that Theorem \ref{holcurvethm} applied to a totally geodesic $\C\P^1$ in $\C\P^2$ 
implies the stability of complex 2-planes as coassociative 4-planes as we already knew 
from Corollary \ref{roundcor}. However, we can in fact prove the following.

\begin{prop}\label{cxunstabprop}
The only complex 2-dimensional cones in $\C^3$, with compact nonsingular complex links in $\C\P^2$, which are stable 
as coassociative cones in $\R^7$ are complex 2-planes.
\end{prop}

\begin{proof}
Suppose $C$ is a counterexample to the statement of the proposition and suppose, for simplicity, that the complex link $\Sigma$ of $C$ 
is connected.  Embed $C$ as a coassociative cone in $\R^7$ and let $\GG$ be the Lie subgroup of $\GG_2$ preserving $C$.  
Since $C$ is non-planar and supposed to be stable, in the notation of Proposition \ref{stabdfn} 
we should have $\dim\check{D}(0)=7$, since every translation of $C$ will define an order $O(1)$ coassociative deformation, and $\dim\check{D}(1)=14-\dim\GG$.  By Theorem \ref{holcurvethm}, this occurs if and only if
\begin{align*}
d_\Sigma^2+d_\Sigma&=6\quad\text{and}\quad
d_\Sigma^2+3d_\Sigma=14-\dim\GG-m_L(3).
\end{align*}
Subtracting these equations and using the fact that $m_L(3)\geq 6$ shows that $2d_\Sigma\leq 2-\dim\GG$, but this contradicts $d_\Sigma>1$ as $C$ is non-planar.
\end{proof}

We now wish to mimic our result for holomorphic curves in $\C\P^2$ for the case of null-torsion pseudoholomorphic curves in $\mathcal{S}^6$.

\begin{thm}\label{nullthm}  Recall Definitions \ref{tubedef}, \ref{pholodfn}, \ref{cLdfn} and \ref{Lapmultdfn}.
Let $\Sigma\subseteq\mathcal{S}^6$ be a compact, connected,
 null-torsion pseudoholomorphic curve of genus $g_\Sigma$ and let $c_1(N_2\Sigma)$ be the first Chern number of $N_2\Sigma$. Let $L$
 be the tube of radius $\frac{\pi}{2}$ in $N_2\Sigma$ about $\Sigma$.  The eigenvalues of $c_L$ in
 $(0,4)$ are $1$, $2$, $3$ and $\lambda+2\in(2,4)\setminus\{3\}$ such that $\lambda(\lambda+2)\in\mathcal{E}_L$.  Moreover,
\begin{align*}\sigma_L(1)&=2g_\Sigma-2c_1(N_2\Sigma)-2, &\sigma_L(2)&=2g_\Sigma-4c_1(N_2\Sigma)-1,\\
\sigma_L(3)&=m_L(3)+2g_\Sigma-6c_1(N_2\Sigma)-2,&\sigma_L(\lambda+2)&=m_L\big(\lambda(\lambda+2)\big).
\end{align*}
\end{thm}

\begin{proof}
It is straightforward to see that Definition \ref{horizontaldfn}  
applies to the natural projection $\pi:L\rightarrow\Sigma$ and that the fibres of $\pi$ are totally geodesic since 
$L$ is a tube of radius $\frac{\pi}{2}$.  
Moreover, by Example \ref{Killingex}, if $\xi$ is the unit vector field given by the $\mathcal{S}^1$-fibration of $L$ over $\Sigma$ then $\xi$ is Killing.  By Proposition \ref{commprop2} it is therefore enough to 
 study \eq{stabpropeq1}-\eq{stabpropeq2} for $\gamma\in C^{\infty}(T^*L\otimes\C)$ and $m\in\Z$ to find the positive eigenvalues 
$\nu$ of $c_L$. 

If we let $\theta$ be the 1-form on $L$ dual to $\xi$ then, since $\Sigma$ has null-torsion, the structure equations for $L$ given in \cite[$\S$6.4]{LotayLag} imply that $\d\theta=-2*\theta$.  This formula, together with the fact that $\Sigma$ is endowed with a
 complex structure, enables us to follow the proof of Theorem \ref{holcurvethm} and see that we essentially have two possible contributions to eigenvalues of $c_L$:
 either eigenfunctions of the Laplacian on $L$ or holomorphic sections of $P_{\nu}=T^{*1,0}\Sigma\otimes H^{-\nu}$ for $\nu\in\Z^+$,
 where $H=N_2\Sigma$ since we can identify the cone over $L$ with $N_2\Sigma$ over $\Sigma$. 
 
Recall the structure equations for $\Sigma$ given in Definition \ref{pholodfn}.  Since $\kappa_{32}=0$, the curvature form of $N_2\Sigma$ is
 $$\d\kappa_{33}=-\theta_1\w\bar{\theta}_1,$$
 so $c_1(N_2\Sigma)<0$ (in fact, it is proportional to $-\vol(\Sigma)$).  Thus, by Riemann--Roch, we have that
 $$h^0(P_{\nu})=c_1(T^{*1,0}\Sigma)-\nu c_1(N_2\Sigma)+1-g_\Sigma=g_\Sigma-1-\nu c_1(N_2\Sigma).$$
 The result now follows from the proof of Theorem \ref{holcurvethm}. 
\end{proof}

\begin{remarks}
\begin{itemize}
\item[]
\item[(a)]
The first Chern number $c_1(N_2\Sigma)$ is the negative of the degree of $\Sigma$ as a holomorphic curve in the 5-quadric in $\C\P^6$, 
so is just $-\frac{1}{4\pi}\vol(\Sigma)$.
\item[(b)]   Theorem \ref{nullthm} will not immediately generalize to the other possible case given by 
Example \ref{Killingex}, namely tubes $L$ of radius $\frac{\pi}{2}$ in $N_1\Sigma$. 
 First, the
 relationship between the derivative of the vertical 1-form $\theta$ and its dual $*\theta$ is not as straightforward.  Second,
 $c_1(N_1\Sigma)$ is often positive, so it is not as easy to calculate $\sigma_L(\nu)$.
\end{itemize}
\end{remarks}

\noindent One can apply Theorem \ref{nullthm} to the Bor\r{u}vka sphere $\Sigma$, for which 
$c_1(N_2\Sigma)=-6$, and recover the result of Corollary \ref{so3cor1}.  Theorem \ref{nullthm} can also be 
applied to a totally geodesic 2-sphere, which is the degenerate case of a null-torsion curve, to prove Corollary 
\ref{roundcor}.

We now prove the analogue of Proposition \ref{cxunstabprop}.

\begin{prop}  Recall Definitions \ref{tubedef} and \ref{pholodfn}.
Let $C$ be a coassociative cone whose link is a tube of radius $\frac{\pi}{2}$ in $N_2\Sigma$ about a 
compact (non-totally geodesic) null-torsion pseudoholomorphic curve $\Sigma$ in $\mathcal{S}^6$.  Then $C$ is not stable.
\end{prop}

\begin{proof}
Suppose for simplicity that $\Sigma$ is connected and let $\GG$ be the subgroup of $\GG_2$ preserving $C$.  By Theorem \ref{nullthm},
$C$ is stable only if
$$2g_\Sigma-4c_1(N_2\Sigma)=8\quad\text{and}\quad m_L(3)+2g_\Sigma-6c_1(N_2\Sigma)-2=14-\dim\GG,$$
since $C$ is non-planar.
These equations force
$$g_\Sigma=m_L(3)+\dim\GG-4\quad\text{and}\quad 2c_1(N_2\Sigma)=m_L(3)+\dim\GG-8.$$
Therefore $m_L(3)+\dim\GG\geq 4$ as $g_\Sigma\geq 0$.  Hence, $c_1(N_2\Sigma)\geq -2$ so 
the degree of $\Sigma$ as a holomorphic curve in $\C\P^6$ must be $1$ or $2$.  We deduce that $\Sigma$ must have genus zero and lie 
in some $\C\P^2$ in $\C\P^6$.  
We now show that $\Sigma$ cannot be a plane curve by
 the structure equations in Definition \ref{pholodfn}, giving our required contradiction.     
 
The embedding of $\Sigma$ as a holomorphic curve in $\C\P^6$ is given by $\bff_3$.  Since 
$\Sigma$ is non-totally geodesic, $\d\bff_3$ depends on $\bff_2$ and $\d\bff_2$ depends on $\bff_1$.  However,
$\d\bff_1$ has a non-zero component in the direction of $\bfu$, which is independent of $\bff_1,\bff_2,\bff_3$, 
since $\theta_1$ is nowhere vanishing on the curve.  It therefore follows that $\Sigma$ cannot lie in some $\C\P^2$ in $\C\P^6$.  
\end{proof}

Theorems \ref{holcurvethm} and \ref{nullthm} invite us to make the following definition.

\begin{dfn}\label{etadfn}
For a compact Riemannian 3-manifold $L$ let $$\eta(L)=\sum_{\lambda\in(0,1]}m_L\big(\lambda(\lambda+2)\big),$$ 
using the notation of Definition \ref{Lapmultdfn}.
\end{dfn}

\begin{remark}
By Theorems \ref{holcurvethm} and \ref{nullthm}, the stability index of certain 2-ruled cones with link $L$ fibered over an algebraic
 curve $\Sigma$ will be the sum of a topological term, determined by the degree and genus of $\Sigma$, and an analytic piece given
 by $\eta(L)$.  We can therefore think of $\eta(L)$ as a sort of ``$\eta$-invariant''.
\end{remark}

We now make an elementary observation, which we state for complex cones though it is equally valid for 
cones whose links are as in Theorem \ref{nullthm}.

\begin{prop}\label{cxstabindprop}
Let $C$ be a complex 2-dimensional cone in $\C^3$ with compact real link $L$ such that $C\setminus\{0\}$ is nonsingular.  
Let $\mathcal{C}$ consist of $\GG_2\ltimes\R^7$ transformations of a deformation family for $C$ as 
a complex cone and recall Definitions \ref{stabdfn} and \ref{etadfn}.  Let $L^{\prime}$ be the link of $C^{\prime}\in\mathcal{C}$. 
 Then
$$\ind_{\mathcal{C}}(C)-\ind_{\mathcal{C}}(C^{\prime})=\eta(L)-\eta(L^{\prime})$$ for all
 $C^{\prime}\in\mathcal{C}$.  
Moreover, there exists an open neighbourhood $\mathcal{C}^{\prime}$ of $C$ in $\mathcal{C}$ such that, for all
 $C^{\prime}\in\mathcal{C}^{\prime}$, $\eta(L)-\eta(L^{\prime})\geq 0$. 
Thus, if $C$ is $\mathcal{C}$-stable then $C^{\prime}$ is $\mathcal{C}$-stable for all $C^{\prime}\in\mathcal{C}^{\prime}$.
\end{prop}

\begin{proof}
By Definitions \ref{stabdfn} and \ref{cLdfn} and Theorem \ref{holcurvethm}, the difference in the $\mathcal{C}$-stability indices
of $C$ and $C^\prime$ is determined by the spectra of $L$ and $L^{\prime}$ as claimed because the degree of the complex link is the same for $C$ and $C^{\prime}$.  

Deformations of $L$  will change the spectrum, but the only way in which $\eta(L)\neq\eta(L^{\prime})$, for a sufficiently small perturbation $L^{\prime}$ of $L$, is if a new element of the spectrum is created strictly above 3 under the deformation and the number of the elements of the spectrum in $(0,3]$ decreases. Thus, $\eta(L)\geq\eta(L^{\prime})$ for all $L^{\prime}$ in some open neighbourhood of $L$, proving the existence of $\mathcal{C}^{\prime}$.

If $C$ is $\mathcal{C}$-stable then $\ind_{\mathcal{C}}(C^{\prime})\leq \ind_{\mathcal{C}}(C)=0$ for all $C^{\prime}\in\mathcal{C}^{\prime}$.  We deduce that $\ind_{\mathcal{C}}(C^{\prime})=0$ from the non-negativity of the stability index.  
\end{proof}

We conclude this section  with an application of Proposition \ref{cxstabindprop}.

\begin{cor}\label{Cacor}
Recall $C_{\mathbf{a}}$, $\mathcal{C}_{\circ}$ and $\mathcal{T}$ defined in Corollary \ref{cxso3cor}.  The cone $C_{\mathbf{a}}$ is
Jacobi integrable and $\mathcal{C}_{\circ}$-stable for all $\mathbf{a}\in\mathcal{T}$.
\end{cor}

\begin{proof}
Recall that $C_2$, with link $L_2$, given by Example \ref{cxso3ex} satisfies $C_2=C_{(\frac{1}{3},\frac{1}{3},\frac{1}{3})}$ and 
$\ind_{\mathcal{C}_{\circ}}(C_2)=0$ by Corollary \ref{cxso3cor}.
Let $L_{\mathbf{a}}$ be the link of $C_{\mathbf{a}}$.  

Let $$\mathcal{T}_{\circ}=\{\mathbf{a}\in\mathcal{T}\,:\,\ind_{\mathcal{C}_{\circ}}(C_{\mathbf{a}})=0\},$$
which is non-empty since $(\frac{1}{3},\frac{1}{3},\frac{1}{3})\in\mathcal{T}_{\circ}$.  By Proposition \ref{cxstabindprop} $\mathcal{T}_{\circ}$ is open.  
Our aim is to show that $\mathcal{T}_{\circ}$ is also closed in $\mathcal{T}$, since then $\mathcal{T}_{\circ}=\mathcal{T}$ by connectedness.
  
Recall that
 $\dim\mathcal{C}_{\circ}=23$ and that $m_{L_{\mathbf{a}}}(3)\geq 6$ by Theorem \ref{holcurvethm}.
Thus 
$$\ind_{\mathcal{C}_{\circ}}(C_{\mathbf{a}})=\sum_{\lambda\in(0,1)}m_{L_{\mathbf{a}}}\big(\lambda(\lambda+2)\big)+\big(m_{L_{\mathbf{a}}}(3)-6\big)$$ 
by Theorem \ref{holcurvethm}.  We deduce that $\mathbf{a}\in\mathcal{T}_{\circ}$ if and only if $m_{L_{\mathbf{a}}}(3)=6$ and $m_{L_{\mathbf{a}}}\big(\lambda(\lambda+2)\big)=0$ for all 
$\lambda\in(0,1)$.     

Suppose, for a contradiction, that $\mathcal{T}_{\circ}$ is not closed, so there exists 
$\mathbf{a}\in\overline{\mathcal{T}_{\circ}}\cap\mathcal{T}\setminus\mathcal{T}_{\circ}$.  Therefore $m_{L_{\mathbf{a}}}(3)>6$ or 
$m_{L_{\mathbf{a}}}\big(\lambda(\lambda+2)\big)>0$ for some $\lambda\in(0,1)$.  However, since the spectrum of the Laplacian varies continuously under deformations of the metric, the latter can occur  
if and only if there exists $\mathbf{a}'\in\overline{\mathcal{T}_{\circ}}$ for which $m_{L_{\mathbf{a}'}}(3)>6$.  Since this cannot happen for 
$\mathbf{a}'\in\mathcal{T}_{\circ}$ we may suppose therefore that $m_{L_{\mathbf{a}}}(3)>6$.  

Hence, there is a Lagrangian Jacobi field 
$v$ on $L_{\mathbf{a}}$ which is independent of the Lagrangian Jacobi fields corresponding to deformations of $C_{\mathbf{a}}$ in $\mathcal{C}_{\circ}$.  
For all $\mathbf{a}'$  close to $\mathbf{a}$ we have that 
$L_{\mathbf{a}'}=\exp_{v'}(L_{\mathbf{a}})$ for some Lagrangian Jacobi field $v'$.  
Moreover, since $L_{\mathbf{a}'}$ is Lagrangian and $v+v'$ is a Lagrangian Jacobi field on $L_{\mathbf{a}}$, 
we may view $v$ as a Lagrangian Jacobi field on $L_{\mathbf{a}'}$ for $\mathbf{a}'$ sufficiently close to $\mathbf{a}$.  However, 
any open neighbourhood of $\mathbf{a}$ meets $\mathcal{T}_{\circ}$ as $\mathbf{a}\in\overline{\mathcal{T}_{\circ}}$, so there exists $\mathbf{a}'\in\mathcal{T}_{\circ}$ such that $L_{\mathbf{a}'}$ has a
 Lagrangian Jacobi field independent of those corresponding to deformations of $C_{\mathbf{a}'}$ in $\mathcal{C}_\circ$.  Thus 
$m_{L_{\mathbf{a}'}}(3)>6$ for some $\mathbf{a}'\in\mathcal{T}_{\circ}$, which is our required contradiction.  
%
\end{proof}

\section{Examples of coassociative 4-folds with conical singularities}\label{exsec}

In this section we produce our examples of coassociative 4-folds with conical singularities.  We describe the 
construction of compact $\GG_2$ manifolds we require and the singular 4-dimensional submanifolds which arise following 
\cite{Kovalevsums} and \cite{Kovalevfib}.  Initially we will have a singular coassociative 4-fold $N$ in a compact 
almost $\GG_2$ manifold $M$, but the ambient $\GG_2$ structure will have torsion.  We then show that $N$ has 
conical singularities which are stable under deformations of the $\GG_2$ structure.  Finally, we deform the 
$\GG_2$ structure on $M$ so that it has no torsion and simultaneously deform $N$ to produce our CS coassociative 4-fold.

\subsection{Examples of compact {\boldmath $\GG_2$} manifolds}

Here we review the relevant material from \cite{Kovalevsums}.  The key 
ingredients will be \emph{Fano 3-folds} and \emph{$K3$ surfaces}, which we now define.

\begin{dfn}
A compact complex 3-dimensional manifold $X$ is a \emph{Fano 3-fold} 
if its first Chern class is positive.   
Equivalently, $X$ has ample anticanonical bundle.  
Fano 3-folds are simply connected and projective.

A \emph{$K3$ surface} $P$ is a simply connected, compact, complex surface with $c_1(P)=0$. 
A generic divisor in the anticanonical linear system of a Fano 3-fold 
is a smooth $K3$ surface by the work of Shokurov \cite{Shokurov}.
\end{dfn}

\begin{remark}
One definition of a Calabi--Yau 3-fold is a compact K\"ahler 3-fold with vanishing first Chern class or, equivalently, 
with trivial canonical bundle.
\end{remark}

The construction of compact $\GG_2$ manifolds in \cite{Kovalevsums} proceeds via the construction of certain 
\emph{non-compact} Calabi--Yau 3-folds.  These non-compact Riemannian manifolds are \emph{asymptotically cylindrical}.  We define these manifolds formally.
 
\begin{dfn}
Let $(Y,g)$ be a Riemannian $n$-manifold.  We say that $Y$ is \emph{asymptotically cylindrical} (with rate $\lambda$) if there 
exist constants $\lambda<0$ and $R>0$, a compact subset $K$ of $Y$, a compact Riemannian $(n\!-\!1)$-manifold $(S,g_S)$ 
and a diffeomorphism $\Psi:(R,\infty)\times S\rightarrow Y\setminus K$ satisfying
\begin{equation*}
\big|\nabla^j\big(\Psi^*(g)-g_{\text{cyl}}\big)\big|=O\left(e^{\lambda t}\right)\qquad\text{as $t\rightarrow\infty$ for all 
$j\in\N$,}
\end{equation*}
where $g_{\text{cyl}}=\d t^2+g_{S}$ is the cylindrical metric on $(0,\infty)\times S$, $\nabla$ is the Levi-Civita 
connection of $g_{\text{cyl}}$ and $|.|$ is calculated with respect to $g_{\text{cyl}}$.
\end{dfn}

We now introduce some important notation.

\begin{dfn}\label{blowupdfn}
Let $\mathcal{X}$ be a maximal deformation family of Fano 3-folds, let $X\in\mathcal{X}$ and let $P,Q$ be $K3$ surfaces in the anticanonical linear system of $X$ such that $P\cap Q$ is a 
nonsingular curve in $X$.  Let $\tilde{X}(P,Q)$ denote the blow-up of $X$ along $P\cap Q$ and let $\tilde{P}$ denote the 
proper transform of $P$ in $\tilde{X}(P,Q)$.  Finally, let $Y(X,P,Q)=\tilde{X}(P,Q)\setminus \tilde{P}$.
\end{dfn}

\noindent Since $\tilde{X}(P,Q)$ is the blow-up of $X$ along $P\cap Q$, we have a smooth map $\varpi:\tilde{X}(P,Q)\rightarrow\C\P^1$
 whose fibres are the proper transforms of the divisors in the pencil defined by $P$ and $Q$.  We may introduce a holomorphic
 coordinate $\zeta$ on $\C\P^1$ such that $\varpi^{-1}(0)=\tilde{P}$ and see that, for some open neighbourhood $U$ of 
 $0$ in the $\zeta$ coordinate, $\varpi^{-1}(U\setminus\{0\})$ is diffeomorphic to $(0,\infty)\times
 P\times\mathcal{S}^1$.  Thus we may view $Y(X,P,Q)$ as a manifold with a cylindrical `end' with cross-section $P\times\mathcal{S}^1$.
 This motivates the next result which follows from \cite[Corollary 6.43]{Kovalevsums}.

\begin{thm}\label{SU3thm}
Use the notation of Definition \ref{blowupdfn}.  There is a smooth complete metric $g_Y$ on $Y=Y(X,P,Q)$ such that 
$(Y,g_Y)$ is asymptotically cylindrical and the holonomy of $g_Y$ is $\SU(3)$.
\end{thm}

Suppose we have a pair of maximal deformation families $\mathcal{X}_1$ and $\mathcal{X}_2$ of Fano 3-folds.  Using the notation 
of Definition \ref{blowupdfn}, we have a pair of 
asymptotically cylindrical complex 3-folds $Y_1=Y_1(X_1,P_1,Q_1)$ and $Y_2=Y_2(X_2,P_2,Q_2)$ with holonomy $\SU(3)$ by Theorem 
\ref{SU3thm}.  We thus have 7-manifolds $Z_i=Y_i\times\mathcal{S}^1$, for $i=1,2$, which are asymptotically 
cylindrical to $(0,\infty)\times P_i\times\mathcal{S}^1\times\mathcal{S}^1$.  In \cite[$\S$4]{Kovalevsums} it is explained that 
if $P_1$ and $P_2$ satisfy a certain `matching condition', then one can apply a `twisted connected sum' construction to $Z_1$ and 
$Z_2$ to get a one-parameter family of compact almost $\GG_2$ manifolds $\{(M_T,\varphi_T,g_{\varphi_T})\,:\,T>T_0\}$ for some 
$T_0>0$.  Moreover, $(\varphi_T,g_{\varphi_T})$ is simply the product $\GG_2$ structure on $Z_i$, as described 
in Definition \ref{CYdfn}, away from the `interpolation region' where $Z_1$ and $Z_2$ are `glued'.  The only question is whether this `matching condition' holds for $P_1$ and $P_2$.  The answer 
\cite[Theorem 6.44]{Kovalevsums} is that there always exist $X_i\in\mathcal{X}_i$ such that $P_1$ and $P_2$ can be chosen 
which satisfy the `matching condition'.  Finally, \cite[Proposition 5.32 \& Theorem 5.34]{Kovalevsums} imply that one can always perturb 
the closed $\GG_2$ structure on $M_T$ to a torsion-free one for sufficiently large $T$.  We can summarize these observations as a theorem.

\begin{thm}\label{G2thm}
Let $\mathcal{X}_1$ and $\mathcal{X}_2$ be maximal deformation families of Fano 3-folds and recall the notation of Definitions 
\ref{pos3formdfn}, \ref{G2structdfn}, \ref{G2mflddfn} and \ref{blowupdfn}.  There exist constants $T_0>0$ and $\lambda<0$ and, 
for $i=1,2$,  $X_i\in\mathcal{X}_i$ and a $K3$ surface $P_i$ in the anticanonical linear system of $X_i$ such that, for all $T>T_0$ 
and suitable $Q_1,Q_2$ as in Definition \ref{blowupdfn}, the following hold.
\begin{itemize}
\item[\emph{(a)}] There is a compact almost $\GG_2$ manifold $(M_T,\varphi_T,g_{\varphi_T})$ which, outside some compact 
set $I_T$, is diffeomorphic to
 the disjoint union of $Y_1(X_1,P_1,Q_1)\times\mathcal{S}^1$ and $Y_2(X_2,P_2,Q_2)\times\mathcal{S}^1$ endowed with the product 
$\GG_2$ structure induced from the asymptotically cylindrical $\SU(3)$ structure given by Theorem \ref{SU3thm}.
 \item[\emph{(b)}] Let $p>4$.  There is a smooth 2-form $\eta_T$ on $M_T$, satisfying $\|\eta_T\|_{L^p_2}\leq c_p e^{\lambda T}$ and 
$\|\eta_T\|_{C^1}\leq c_p e^{\lambda T}$ for some constant $c_p>0$, such that $\varphi_T+\d\eta_T\in C^{\infty}(\Lambda^3_+T^*M_T)$ and 
the metric $g_{\varphi_T+\d\eta_T}$ on $M_T$ has holonomy $\GG_2$.
\end{itemize}
\end{thm}

\noindent The compact set $I_T$ is the interpolation region between $Z_1$ and $Z_2$.

\begin{note}
Of course, the $M_T$ are topologically the same for all $T$, 
so we can view $(\varphi_T,g_{\varphi_T})$ as a one-parameter family of closed $\GG_2$ structures on a 7-manifold $M$.
\end{note}

We now have the following important result.

\begin{prop}\label{singfibprop}
Use the notation of Theorem \ref{G2thm}.  Outside $I_T$, $M_T$ is fibered by $K3$ surfaces which are coassociative with respect to 
$\varphi_T$.  Moreover, for generic $Q_1$ and $Q_2$, there exist coassociative $K3$ surfaces in $(M_T,\varphi_T,g_{\varphi_T})$
 whose singularities are isolated and are ordinary double points.
\end{prop}

\begin{proof}
After Definition \ref{blowupdfn} we noted that we have a fibration $\varpi:\tilde{X}(P,Q)\rightarrow\C\P^1$ whose fibres are
$K3$ surfaces.  Necessarily some of these fibres will be singular and the generic singularity is an ordinary double point.  Recall 
that we have the freedom to choose any smooth $K3$ surface $Q$ in the anticanonical linear 
system of $X$ which meets $P$ in a nonsingular curve.  Therefore, through 
generic choice of $Q$ we can be assured that there are fibres other than the exceptional divisor whose only singularities
 are ordinary double points.  
Thus, for generic $Q_1$ and $Q_2$, we have a fibration 
$$\varpi_T:M_T\setminus I_T\cong (Y_1\times\mathcal{S}^1)\sqcup (Y_2\times\mathcal{S}^1)\rightarrow
 (\C\P^1\sqcup\C\P^1)\times\mathcal{S}^1$$ with fibres which are $K3$ surfaces in 
$Y_i\times\{x\}$ for some $i$ and some $x\in\mathcal{S}^1$.  Moreover, there are some fibres of $\varpi_T$ which 
only have ordinary double point singularities.  Since $Y_1$ and 
$Y_2$ are Calabi--Yau manifolds and the almost $\GG_2$ structure on $M_T\setminus I_T$ agrees with the product $\GG_2$ structure 
on $Y_1\times\mathcal{S}^1$ and $Y_2\times\mathcal{S}^1$, a simple generalization of Corollary \ref{coasscor} leads us to deduce that
 the fibres of $\varpi_T$ are coassociative with respect to $\varphi_T$.  
\end{proof}

\begin{remark}
By studying the $\GG_2$ structure on $I_T$, it is shown in \cite{Kovalevfib} that one can extend the coassociative fibration $\varpi_T$
through $I_T$.
\end{remark}

For convenience we introduce the following notation.

\begin{dfn}\label{singfibdfn}
Use the notation of Theorem \ref{G2thm} and let $Q_1$ and $Q_2$ be generic so that Proposition \ref{singfibprop} applies.  
Let $\Gamma_T$ denote the set of coassociative $K3$ surfaces in $(M_T,\varphi_T,g_{\varphi_T})$ which have isolated ordinary 
double point singularities.  By Proposition \ref{singfibprop}, there exist $N\in\Gamma_T$ 
such that $N\subseteq M_T\setminus I_T$.
\end{dfn} 

Our aim now is to show that some of the singular coassociative 4-folds in $\Gamma_T$ are `stable' 
under the deformation from the closed $\GG_2$ structure $\varphi_T$ to the torsion-free $\GG_2$ structure $\varphi_T
+\d\eta_T$ given in Theorem \ref{G2thm}(b).  

\subsection{Stable coassociative conical singularities}

We begin with a crucial result which allows us to implement our stability theory.

\begin{prop}\label{csfibprop1}
Use the notation of Theorem \ref{G2thm} and Definition \ref{singfibdfn}.
If $N\in\Gamma_T$, 
then $N$ is a CS coassociative 4-fold in
$(M_T,\varphi_T,g_{\varphi_T})$ in the sense of Definition \ref{csdfn}.  Moreover, the singularities of $N$ have cones in 
the family $\mathcal{C}_{\circ}$ given in Corollary \ref{cxso3cor}.

\end{prop}

\begin{proof}
First observe that $N$ is clearly a connected coassociative integral current with $\partial N=\emptyset$.  Second, if $z$ is a
singular point of $N$ then it is an ordinary double point of a complex surface, so the tangent cone at $z$  has 
multiplicity one and is modelled on a cone in $\mathcal{C}_{\circ}$ by definition.  
Since cones in $\mathcal{C}_{\circ}$ are Jacobi integrable by Corollary \ref{Cacor}, we may apply Corollary \ref{cosingcor} 
to deduce the result.
\end{proof}

\begin{prop}\label{stabprop}
Let $N$ be a CS coassociative 4-fold in an almost $\GG_2$ manifold.  Suppose that the singularities of $N$ are $z_1,\ldots,z_s$ 
with rate $\mu$ and cones $C_1,\ldots,C_s$ such that:
\begin{itemize}
\item[\emph{(i)}] $C_i$ is in the family $\mathcal{C}_{\circ}$ given in Corollary \ref{cxso3cor} for all $i$; and 
\item[\emph{(ii)}] $(1,\mu]\cap\mathcal{D}=\emptyset$, where $\mathcal{D}$ is given in Definition \ref{Ddfn}.
\end{itemize}  
If $\mathfrak{C}=\mathcal{C}_{\circ}^s$, then
$\mathcal{O}(N,\mu,\mathfrak{C})=\{0\}$, in the notation of Theorem \ref{maindefthm}.
\end{prop}

\begin{proof}
This follows immediately from Definition \ref{stabdfn} and Corollary \ref{Cacor}.
\end{proof}

Theorem \ref{mainthm} now follows from our final result.

\begin{thm}
Use the notation of Theorem \ref{G2thm} and Definition \ref{singfibdfn}.  Let $T>T_0$ and let $N\in\Gamma_T$ be such that $N\subseteq M_T\setminus I_T$.  Making $T_0$ larger if necessary,  
there exists a CS deformation $N^{\prime}$ of $N$ which is coassociative with respect to $\varphi_T+\d\eta_T$.
\end{thm}

\begin{proof}
Notice that there are $N\in\Gamma_T$ which do not lie in $I_T$ and that we are free to make the 
rate $\mu$ at the singularities of $N$ lower if necessary to satisfy $(1,\mu]\cap\mathcal{D}=\emptyset$.  Hence,
Propositions \ref{csfibprop1} and \ref{stabprop} imply $N$ is a CS coassociative 4-fold with respect to $\varphi_T$ and that
 $\mathcal{O}(N,\mu,\mathfrak{C})=\{0\}$.   

Recall that $M_{T^{\prime}}$ is diffeomorphic to some 7-manifold $M$ for all $T^{\prime}\geq T_0$.  Let $\mathcal{F}=\{(\varphi_{t}+s\d\eta_{t},g_{\varphi_{t}+s\d\eta_{t}})\,:\,s\in\R,\,t>T_0\}$.  Then $\mathcal{F}$ is a smooth 
2-dimensional family of closed $\GG_2$ structures on $M$.  Moreover, for each $T>T_0$ we may parameterize $\mathcal{F}$ by 
$(u,v)\in B(0;1)$ via
$$(u,v)\in B(0;1)\subseteq\R^2\mapsto \varphi^{(u,v)}=\varphi_{t_T(u,v)}+s_T(u,v)\d\eta_{t_T(u,v)}$$
where 
\begin{align*}
s_T(u,v)&=\frac{(T-T_0)v}{(1-u)^2+v^2}\\ \intertext{and}t_T(u,v)&=\frac{T_0(1-u)^2+T_0v^2+(T-T_0)(1-u)}{(1-u)^2+v^2}.\end{align*}
With this parameterization of $\mathcal{F}$ we see that $\varphi^{(0,0)}=\varphi_T$.  

Notice that, since $N\cap I_T=\emptyset$ there exists $\tau>0$ such that $T>T_0+\tau$ and $N\cap I_{T^{\prime}}=\emptyset$ for all 
$T^{\prime}\in (T-\tau,T+\tau)$.  As the $\GG_2$ structure outside $I_{T^{\prime}}$ is the same for all 
$T^{\prime}\in(T-\tau,T+\tau)$ by Theorem \ref{G2thm}(a), $N$ is coassociative with respect to $\varphi_{T^{\prime}}$ 
for all such $T^{\prime}$.  We therefore see that $[\varphi^{(u,v)}]=0$ in $H^3_{\text{cs}}(\hat{N})$ for all
 $(u,v)$ sufficiently near $(0,0)$.  Hence, by Theorem \ref{G2deformthm}, there exists some $\delta_N>0$ 
such that for all $(u,v)\in B(0;\delta_N)$ there exists a
  CS deformation $N^{(u,v)}$ of $N$ which is coassociative with respect to $\varphi^{(u,v)}$.  

To complete the proof, we need to show that
$\varphi^{(u,v)}=\varphi_T+\d\eta_T$ for some $(u,v)\in B(0;\delta_N)$.  As discussed after Theorem \ref{G2deformthm}, this occurs if 
$\|\d\eta_T\|_{C^1}<\epsilon$ and $\|\d\eta_T\|_{L^p_2}<\epsilon$ for some $p>4$, where $\epsilon>0$ is determined by the geometry 
near $N$ with respect to $g_{\varphi_T}$.  By Theorem \ref{G2thm}(b), $\|\d\eta_T\|_{C^1}$ and $\|\d\eta_T\|_{L^p_2}$ are of 
order $O(e^{\lambda T})$ for some $\lambda<0$.  Moreover, by Theorem \ref{G2thm}(a), the geometry near $N$ is not changing as $T$ 
varies, so the constant $\epsilon$ can be chosen to be independent of $T$.  Thus, we can ensure that the relevant norms of 
$\d\eta_T$ are sufficiently small by making $T_0$ larger.
\end{proof}

\begin{ack}
The author particularly thanks Robert Bryant, Dominic Joyce and Alexei Kovalev for useful conversations, 
and University College, Oxford, for hospitality 
during the course of this project.  The author also thanks the referee for helpful comments and suggestions.  The author is currently 
supported by an EPSRC Career Acceleration Fellowship.
\end{ack}



\begin{thebibliography}{99}

\bibitem{AdamsSimons} D.~Adams and L.~Simon, {\it Rates of Asymptotic Convergence Near Isolated
Singularities of Geometric Extrema}, Indiana Univ. Math. J. {\bf 37} (1988), 225--254.

\bibitem{Bartnik} R.~Bartnik, {\it The Mass of an Asymptotically
Flat Manifold}, Comm. Pure Appl. Math. {\bf 39} (1986), 661--693.

\bibitem{Berard} L.~B\'erard Bergery and J.-P.~Bourguignon, {\it Laplacians and Riemannian Submersions with Totally Geodesic Fibres}, Illinois J. Math. {\bf 26} (1982), 181--200.

\bibitem{BryantOct} R. L. Bryant, {\it Submanifolds and Special Structures on the Octonions}, J. Diff. Geom. 
{\bf 17} (1982), 185--232.

\bibitem{BryantG2struct} R.~L.~Bryant, {\it Some Remarks on
$\GG_2$-Structures}, Proceedings of G\"okova Geometry-Topology Conference 2005, edited by S.~Akbulut, T.~\"Onder and R.~J.~Stern, 
International Press, 2006.

\bibitem{CHNP} A. Corti, M. Haskins, J. Nordstr\"om and T. Pacini, {\it $\GG_2$ Manifolds and Associative Submanifolds via Semi-Fano 3-folds}, 
preprint, arXiv:1207.4470.

\bibitem{Dillen} F.~Dillen, L.~Verstraelen and L.~Vrancken, {\it Classification of Totally Real 3-Dimensional Submanifolds of
$\mathcal{S}^6(1)$ with $K\geq1/16$}, J. Math. Soc. Japan
\textbf{42} (1990), 565--584.

\bibitem{Etnyre} J.~Etnyre and R.~Ghrist, {\it Contact Topology and Hydrodynamics: I. Beltrami Fields and the Seifert Conjecture}, 
 Nonlinearity {\bf 13} (2000), 441--458.

\bibitem{Foxcoass} D. Fox, {\it Coassociative Cones that are Ruled by 2-Planes}, Asian J. Math. {\bf 11} (2007), 535--554.


\bibitem{HarLaw} R.~Harvey and H.~B.~Lawson, {\it Calibrated Geometries}, Acta Math. {\bf 148} (1982), 47--152.

\bibitem{Haskins} M.~Haskins, {\it The Geometric Complexity of Special Lagrangian $T^2$-Cones}, Invent. Math. {\bf 157} (2004), 11--70.

\bibitem{Ikeda} A.~Ikeda, {\it On the Spectrum of Homogeneous Spherical Space Forms}, Kodai Math. J. {\bf 18} (1995), 57--67.

\bibitem{Joyexcept} D.~D.~Joyce, {\it Compact Manifolds with Special Holonomy},
OUP, Oxford, 2000.

\bibitem{JoyceSLsing5} D.~D.~Joyce, {\it Special Lagrangian Submanifolds
with Isolated Conical Singularities. V. Survey and Applications}, J. Diff. Geom. {\bf 63} (2003), 299--347.

\bibitem{JoyceSLsing1} D.~D.~Joyce, {\it Special Lagrangian Submanifolds
with Isolated Conical Singularities. I. Regularity}, Ann. Global
Ann. Geom. {\bf 25} (2004), 201--251.

\bibitem{JoyceSLsing2} D.~D.~Joyce, {\it Special Lagrangian Submanifolds
with Isolated Conical Singularities. II. Moduli Spaces}, Ann. Global
Ann. Geom. {\bf 25} (2004), 301--352.

\bibitem{JoyceSLsing3} D.~D.~Joyce, {\it Special Lagrangian Submanifolds
with Isolated Conical Singularities. III. Desingularization, The Unobstructed Case}, Ann. Global
Ann. Geom. {\bf 26} (2004), 1--58.

\bibitem{JoyceSLsing4} D.~D.~Joyce, {\it Special Lagrangian Submanifolds
with Isolated Conical Singularities. IV. Desingularization, Obstructions and Families}, Ann. Global
Ann. Geom. {\bf 26} (2004), 117--174.


\bibitem{JoyceRiem} D.~D.~Joyce, {\it Riemannian Holonomy Groups and Calibrated Geometry}, 
Oxford Graduate Texts in Mathematics {\bf 12}, OUP, Oxford, 2007. 

\bibitem{Kovalevsums} A.~G.~Kovalev, {\it Twisted Connected Sums and Special
Riemannian Holonomy}, J. Reine Angew. Math. {\bf 565} (2003), 125--160.

\bibitem{Kovalevfib} A.~G.~Kovalev, {\it Coassociative $K3$ Fibrations of Compact
$\GG_2$-Manifolds}, preprint, arXiv:math/0511150.

\bibitem{LO} H.~B.~Lawson and R.~Osserman, {\it Non-existence, Non-uniqueness and Irregularity
of Solutions to the Minimal Surface System}, Acta Math. {\bf 139} (1977), 1--17.


\bibitem{LockhartMcOwen} R.~B.~Lockhart and R.~C.~McOwen, {\it
Elliptic Differential Operators on Noncompact Manifolds}, Ann. Sc.
Norm. Super. Pisa Cl. Sci. {\bf 12} (1985), 409--447.

\bibitem{Lotay2r} J. D. Lotay, {\it 2-Ruled Calibrated 4-folds in $\R^7$ and $\R^8$}, J. London Math. Soc. {\bf 74} (2006), 219--243.

\bibitem{Lotaycs} J.~D.~Lotay, {\it Coassociative 4-folds with Conical
Singularities}, Comm. Anal. Geom. {\bf 15} (2007), 891--946.

\bibitem{Lotaydesing} J.~D.~Lotay, {\it Desingularization of Coassociative 4-folds with Conical Singularities}, Geom.~Funct.~Anal. 
{\bf 18} (2008), 2055--2100.

\bibitem{Lotayac} J.~D.~Lotay, {\it Deformation Theory of Asymptotically Conical Coassociative 4-folds}, Proc. London Math. Soc. {\bf 99} (2009), 386--424.

\bibitem{LotayLag} J.~D.~Lotay, {\it Ruled Lagrangian Submanifolds of the 6-Sphere}, Trans. Amer. Math. Soc.
{\bf 363} (2011), 2305--2339.

\bibitem{Mashimo} K.~Mashimo, {\it Homogeneous Totally Real Submanifolds of $\mathcal{S}^6$}, Tsukuba J. Math. {\bf 9} (1985),
185--202.

\bibitem{Mazya} V.~G.~Maz'ya and B.~Plamenevskij, {\it Elliptic Boundary Value Problems}, Amer. Math. Soc. Transl. {\bf 123} (1984), 
1--56.

\bibitem{McLean} R.~C.~McLean, {\it Deformations of Calibrated Submanifolds}, Comm. Anal. Geom. {\bf 6} (1998), 705--747.

\bibitem{Morgan} F.~Morgan, {\it Geometric Measure Theory: A Beginner's Guide}, Fourth Edition, Academic Press, San Diego, 2009. 

\bibitem{Nagy} P.-A.~Nagy, {\it Un Principe de S\'eparation des Variables pour le Spectre du Laplacien des Formes Differentielles 
et Applications}, PhD thesis, Universit\'e de Savoie, 2001.

\bibitem{Ohnita} Y.~Ohnita, {\it Stability and Rigidity of Special Lagrangian Cones over Certain Minimal Legendrian Orbits}, 
Osaka J. Math. {\bf 44} (2007), 305--334.

\bibitem{Salamon} S.~Salamon, {\it Riemannian Geometry and Holonomy Groups}, Pitman Research Notes in Mathematics \textbf{201},
 Longman, Harlow, 1989.

\bibitem{Shokurov} V.~V.~Shokurov, {\it Smoothness of a General Anticanonical Divisor on a Fano Variety}, Izv. Akad. Nauk SSSR Ser.
Mat. {\bf 43} (1979), 430--441; English transl.: Math. USSR Izb {\bf 14} (1980), 395--405.

\bibitem{Simon} L.~Simon, {\it Isolated Singularities of Extrema of Geometric Variational Problems}, in Harmonic Mappings and 
Minimal Immersions, edited by E.~Giusti, Lecture Notes in Mathematics {\bf 1161}, Springer--Verlag, Berlin, 1985, 206--277.

\bibitem{SYZ} A.~Strominger, S.~T.~Yau, and E.~Zaslow,  {\it Mirror Symmetry is T-Duality}, Nucl. Phys. B {\bf 479} (1996), 243--259.

\bibitem{Tanno} S.~Tanno, {\it The First Eigenvalue of the Laplacian on Spheres}, T\^ohoku Math. J. {\bf 31} (1979), 179--185.

\bibitem{Urakawa} H.~Urakawa, {\it On the Least Eigenvalue of the Laplacian for Compact Group Manifolds}, J. Math. Soc.
Japan {\bf 31} (1979), 209--226.

\bibitem{Vrancken} L.~Vrancken, {\it Killing Vector Fields and Lagrangian Submanifolds of the Nearly K\"ahler $\mathcal{S}^6$}, 
J. Math. Pures Appl. {\bf 77} (1998), 631--645.

\bibitem{Watson} B.~Watson, {\it Manifold Maps Commuting with the Laplacian}, J. Diff. Geom. {\bf 8} (1973), 85--94.

\end{thebibliography}
\end{document}